\documentclass[24pt,twoside]{article}
\usepackage[centertags]{amsmath}
\usepackage{amsfonts}
\usepackage{amssymb}
\usepackage{amsthm}
\usepackage{newlfont}
\usepackage{color}
\usepackage[svgnames,dvipsnames]{xcolor}
\usepackage{graphicx}

\date{\color{green}  2013 November}
\def \n {\noindent}
\usepackage{fancyhdr}
\usepackage[T1]{fontenc}
\usepackage[french]{babel}
\usepackage{array,multirow,makecell}
\setcellgapes{1pt}
\makegapedcells
\newcolumntype{R}[1]{>{\raggedleft\arraybackslash }b{#1}}
\newcolumntype{L}[1]{>{\raggedright\arraybackslash }b{#1}}
\newcolumntype{C}[1]{>{\centering\arraybackslash }b{#1}}

\begin{document}
 \n \fcolorbox{black}{yellow}{
\begin{minipage}{15.6cm}
\quad\\

\begin{center}
\n {\bf {\color{red}{\large Application of the criterion of Li-Wang to a five dimensional epidemic model of COVID-19. Part I}}}\\
\quad\\
{\bf {\color{blue}{\large Abdelkader Intissar }}}\\
\quad\\
{\color{blue}$^{(1)}$} \quad\quad \quad{\color{blue}$^{(2)}$} \\
 \includegraphics[scale=0.33]{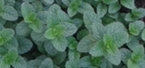} \\

\end{center}
\end{minipage}}\\
\quad\\

\n {\color{blue}$^{(1)}$}  Equipe d'Analyse Spectrale, UMR-CNRS n: 6134, Université de Corse, Quartier Grossetti, 20250 Corté, France.\\
\n {\color{blue}intissar@univ-corse.fr}\\

\n {\color{blue}$^{(2)}$} Le Prador, 129 rue du commandant Rolland, 13008 Marseille, France.\\
\n {\color{blue}abdelkader.intissar@orange.fr}\\

\n {\it {\bf {\color{red}{\large Abstract}}}\\

\n The dynamics of many epidemic models for infectious diseases that spread in a single host population demonstrate a threshold phenomenon.\\

\n  If the basic reproduction number $R_{0}$ is below unity, the disease-free equilibrium $P_{0}$ is globally stable in the feasible region and the disease always dies out.\\

\n  If $R_{0} > 1$, a unique endemic equilibrium $P^{*}$ is globally asymptotically stable in the interior of the feasible region and the disease will persist at the endemic equilibrium if it is initially present. \\

\n In this paper (Part I), we reinvestigate the study of the stability or the non stability of a mathematical Covid-19 model constructed by Nita H. Shah, Ankush H. Suthar and Ekta N. Jayswal in  {\color{blue} https://doi.org/10.1101/2020.04.04.20053173.}\\ We use a criterion of Li-Wang for stability of matrices {\bf{\color{blue}[Li-Wang]}} on the second additive compound matrix associated to their model.\\

\n In second paper (Part II), In order to control the Covid-19 system, i.e., force the trajectories to go to the equilibria we will add some control parameters with uncertain parameters to stabilize the five-dimensional  Covid-19 system studied in this paper. Based on compound matrices theory, we apply  in {\bf {\color{blue}[Intissar]}} again the criterion of Li-Wang to study the stability of equilibrium points of Covid-19 system with uncertain parameters. In  this part II, all sophisticated technical calculations including those in part I are given in appendices.}\\

\n {\bf Keywords :} Epidemic models, Endemic equilibrium, Stability of matrices, Compound matrices,  Dynamical systems, Covid-19 model.\\
\newpage

\begin{center}
\n  {\color{red} $ \S  1$ { \bf Introduction and preliminary results}}\\
\end{center}

\n A mathematical Covid-19 model is constructed  in {\bf {\color{blue} [Shah et al]}}  by Nita H. Shah, Ankush H. Suthar and Ekta N. Jayswal to study human to human transmission of the  Covid-19. \\

\n The model consists all possible human to human transmission of the virus. \\

\n The Covid-2019 is highly contagious in nature and infected cases are seen in most of the countries around the world, hence in the model the susceptible population class is ignored and whole population is divided in five compartments:\\

\n (1) class of exposed individuals E(t) (individuals surrounded by infection by not yet infected),\\

\n (2) class of infected individuals by Covid-19 I(t), \\
 
\n (3)  class of critically infected individuals by Covid-19 C(t),\\

\n (4) class of hospitalised individuals H(t), \\

\n and \\

\n (5) class of dead individuals due to Covid-19 D(t).\\ 

\n Human to human transmission dynamics of Covid-19 is describe graphically in\\

\begin{center}
\includegraphics[scale=0.82]{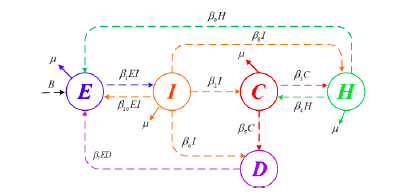}
\end{center}

\n  {\color{blue}$\bullet$} {\bf {\color{red}Table of parameters used in the model is described as follow:}}\\

\n B : Birth rate of class of exposed individuals : {\color{blue}0.80} Calculated\\

\n $\mu$  Natural death rate : {\color{blue}0.01} Assumed\\

\n $\beta_{1}$ :  Transmission rate of individuals moving from exposed to infected class: {\color{blue}0.55} Calculated\\

\n $\beta_{2}$ :  Rate at which infected individuals goes into sever condition or in critical condition : {\color{blue}0.40} Calculated\\

\n $\beta_{3}$ :  Rate at which critically infected individuals get hospitalized : {\color{blue} 0.60} Calculated\\

\n $\beta_{4}$ :  Rate by which hospitalized individuals not recovered and remain in critical condition : {\color{blue} 0.80} Calculated\\

\n $\beta_{5}$ :  Mortality rate of critically infected individuals : {\color{blue} 0.34} Calculated\\

\n $\beta_{6}$ :  Mortality rate of infected individuals : {\color{blue} 0.30} Calculated\\

\n $\beta_{7}$ :  Rate by which infected dead body spreads infection: {\color{blue}0.35} Assumed\\

\n $\beta_{8}$ : Rate at which infected individuals get hospitalized :  {\color{blue}0.30} Calculated\\

\n $\beta_{9}$ : Rate at which hospitalised individuals get recovered and become exposed again : {\color{blue}0.35} Assumed \\

\n $\beta_{10}$ :  Rate at which infected individuals recovered themselves due to strong immunity and again become exposed \\

\n Using the above representation, dynamical system of set of nonlinear differential for the model is formulated as follow: \\

\n  $\left \{ \begin{array} {c} \displaystyle{\frac{dE}{dt} = B - \beta_{1}EI  + \beta_{7}ED + \beta_{9}H  + \beta_{10}EI  - \mu E}\quad\\
\quad\\
\displaystyle{\frac{dI}{dt} =  \beta_{1}EI  - \beta_{2}I - \beta_{6}I - \beta_{8}I -\beta_{10}EI  - \mu I} \quad \quad \quad\\
\quad\\
\displaystyle{\frac{dC}{dt} =  \beta_{2}I  - \beta_{5}C - \beta_{3}C + \beta_{4}H  - \mu C} \quad \quad \quad \quad \quad \quad\\
\quad\\
\displaystyle{\frac{dH}{dt} =  \beta_{3}C  - \beta_{4}H + \beta_{8}I - \beta_{9}H  - \mu H} \quad \quad \quad \quad \quad \quad\\
\quad\\
\displaystyle{\frac{dD}{dt} =  \beta_{5}C  + \beta_{6}I -  \beta_{7}DE} \quad \quad \quad \quad \quad \quad \quad \quad \quad \quad \quad\\
\end{array} \right.$ \hfill {} {\bf{\color{blue} (covid-19)}}\\
\quad\\

\n {\bf {\color{red} Remark 1.1}}\\

\n (i) All of the parameters in {\bf{\color{blue} (covid-19)}} are assumed to be nonnegative.\\

\n (ii) $ \displaystyle{\frac{dE}{dt} + \frac{dI}{dt} + \frac{dC}{dt} + \frac{dH}{dt}  + \frac{dD}{dt} = B - \mu(E + I + C + H + D)}$ \hfill { } {\bf{\color{blue}(1.1)}}\\

\n (iii) $ \displaystyle{\frac{dE}{dt} + \frac{dI}{dt} + \frac{dC}{dt} + \frac{dH}{dt}  + \frac{dD}{dt} \leq 0 \iff  E + I + C + H + D \leq \frac{B}{\mu}}$ \hfill { } {\bf{\color{blue}(1.2)}}\\ 

\n (iv) A detailed description of the model can be found in {\bf {\color{blue} [Shah et al ]}} and there in references.:\\

\n (v)  For other mathematical systems of epidemic models, we can consult these references {\bf {\color{blue} [Li et al]}}, {\bf {\color{blue} [Beretta et al]}} and {\bf {\color{blue} [Sun et al]}}.\\

\n {\bf {\color{red} Theorem 1.2}}\\ 

\n (i) The positive orthant $\mathbb{R}_{+}^{5}$ is positively invariant under the flow of {\bf{\color{blue} (covid-19)}}. Precisely, if $E(0) > 0$; $I(0) > 0$; $C(0) > 0$; $H(0) > 0$; $D(0)  > 0$ then $\forall t > 0$; $E(t) > 0$; $I(t) > 0$ $C(t) > 0$; $H(t) > 0$; $D(t)  > 0$.\\

\n {\bf {\color{red} Proof}}\\

\n (i) {\color{red}$\bullet$} Let's suppose $I(0) > 0$, then from the  second equation of {\bf{\color{blue} (covid-19)}}), if \\

\n  $\displaystyle{ \chi(t) = (\beta_{10} - \beta_{1})E + \beta_{2} + \beta_{6} + \beta_{8} +\mu}$ then the integration from $0$ to $t > 0$ gives :\\

\n $\displaystyle{ I(t) = I(0) e^{-\int_{0}^{t}\chi(s)ds}}$.\\

\n  Therefore $I(t) > 0; \, \forall  t \geq 0.$\\

\n {\color{red}$\bullet$} Consider the following sub-equations related to the variables $C$ and $H$ :\\

\n  $\left \{ \begin{array} {c}\displaystyle{\frac{dC(t)}{dt} =  \beta_{2}I(t)  - \beta_{5}C(t) - \beta_{3}C(t) + \beta_{4}H(t)  - \mu C(t)}\\
\quad\\
\displaystyle{\frac{dH(t)}{dt} =  \beta_{3}C(t)  - \beta_{4}H(t) + \beta_{8}I(t) - \beta_{9}H(t)  - \mu H(t)} \\
\quad \\
C(0) > 0 \, and \,  H(0) > 0 \quad \quad   \quad \quad  \quad \quad  \quad \quad  \quad \quad  \quad \quad  \quad \quad \\
\end{array} \right.$ \hfill {} {\bf{\color{blue} (1.3)}}\\
\quad\\

\n The system {\bf{\color{blue} (1.3)}} takes the matrix form:\\

\n $\displaystyle{\frac{dU}{dt} = MU(t) + F(t) }$ \hfill {} {\bf{\color{blue} (1.4)}}\\

\n where\\

\n  $M = $ $\left ( \begin{array} {cc} - (\beta_{5} + \beta_{3}  + \mu)& \beta_{4}\\
\quad\\
 \beta_{3}&  - (\beta_{4} + \beta_{9}  + \mu)\\
 \end{array} \right)$ \hfill {} {\bf{\color{blue} (1.5)}}\\
 
 \n and\\

\n  $U(t) =$ $\left ( \begin{array} {c} C(t)\\
\quad\\
H(t)\\
\end{array} \right)$ ;  \n  $U(0) =$ $\left ( \begin{array} {c} C(0)\\
\quad\\
H(0)\\
\end{array} \right)$ $ > 0$ and  $F(t) =$ $\left ( \begin{array} {c} \beta_{2}I(t) \\
\quad\\
 \beta_{8}I(t) \\
\end{array} \right)$  $> 0$.\\ 

\n  One can turn the Cauchy problem {\bf{\color{blue} (1.3)}} into an integral equation by using the following so called Duhamel formula :\\

\n $\displaystyle{U(t) = e^{tM}U(0) + \int_{0}^{t}e^{(t -s)M}F(s)ds}$\\

\n Therefore $C(t) > 0; \, \forall  t \geq 0.$ and  $H(t) > 0; \, \forall  t \geq 0.$\\

\n We can observe  also that  $M$ is a Metzler matrix (a matrix  $A=(a_{ij} 1 \leq i, j \leq n$ is  a Metzler matrix  if all of its elements are non-negative except for those on the main diagonal, which are unconstrained.)  That is, a Metzler matrix is any matrix $A$ which satisfies $\displaystyle{ A=(a_{ij});\quad a_{ij}\geq 0,\quad i\neq j.}$\\

\n Thus, {\bf{\color{blue} (1.3)}} is a monotone system. It follows that, $\mathbb{R}_{+}^{2}$ is invariant under the flow of {\bf{\color{blue} (1.3)}}.\\

\n {\color{red}$\bullet$} Let's suppose $E(0) > 0$, then from the  first equation of {\bf{\color{blue} (covid-19)}}), if\\

\n $\displaystyle{ \chi(t) = (\beta_{1} - \beta_{10})I - \beta_{7}D + \mu}$ and  $\displaystyle{\pi(t) = \beta_{9}H(t) + B}$ which is $ > 0$.\\

\n the integration from $0$ to $t > 0$ gives :\\

\n $\displaystyle{ E(t) = E(0) e^{-\int_{0}^{t}\chi(s)ds}  + e^{-\int_{0}^{t}\chi(s)ds}\int_{0}^{t} \pi(u)e^{\int_{0}^{u}\chi(w)dw} du}$.\\

\n Therefore $E(t) > 0; \, \forall  t \geq 0.$\\

\n {\color{red}$\bullet$} Let's suppose $D(0) > 0$, then from the  $5^{th}$ equation of {\bf{\color{blue} (covid-19)}}), if\\

\n $\displaystyle{ \chi(t) = \beta_{7}E}$ and  $\displaystyle{\pi(t) = \beta_{6}I(t) + \beta_{5}C(t)}$ which is $ > 0$.\\

\n the integration from $0$ to $t > 0$ gives :\\

\n $\displaystyle{ D(t) = D(0) e^{-\int_{0}^{t}\chi(s)ds}  + e^{-\int_{0}^{t}\chi(s)ds}\int_{0}^{t} \pi(u)e^{\int_{0}^{u}\chi(w)dw} du}$.\\

\n Therefore $D(t) > 0; \, \forall  t \geq 0.$\\
 \begin{center}
\includegraphics[scale= 0.7]{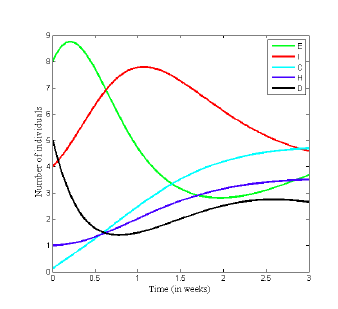}\\
\end{center}
\n {\bf {\color{red} Remark 1.3}}\\

\n (i) The theorem 1.2 ensures the existence and uniqueness of global (in time) solution of system {\bf{\color{blue} (covid-19)}}.\\

\n (ii) Let {\color{red}$\Lambda$} be the domain $\displaystyle{{\color{red}\Lambda} = \{(E, B, C, H, D) \in \mathbb{R}_{+}^{5} ; E + I + C + H + D \leq \frac{B}{\mu}\}}$ then this domain is positively invariant, and all the solutions of the system {\bf{\color{blue} (covid-19)}} are remain in this domain.\\

\n Consider the following n-dimensional system:\\

\n $\displaystyle{x'(t) = f(x(t)) ; t \geq 0}$  \hfill {} {\bf{\color{blue} (1.6)}}\\

\n where  $\displaystyle{f : \Omega \subset \mathbb{R}^{n} \longrightarrow  \mathbb{R}^{n}}$ is $\mathcal{C}^{1}$-function. \\

\n {\bf {\color{red} Definition 1.4}}\\

\n {\color{red} $\bullet$} We say that $x^{*}$ is an equilibrium point of {\bf{\color{blue} (1.6)}} if $f(x^{*}) = 0$.\\

\n {\color{red} $\bullet$} We will say that  an equilibrium point $x^{*}$ is {\color{red} stable} if \\

\n $\displaystyle{ \forall \, \epsilon  > 0 , \exists \delta > 0 \, such \, that \, \mid\mid x - x^{*} \mid\mid  < \delta \, and \, t   > 0 \Longrightarrow \mid\mid \phi_{t}(x) - x^{*} \mid\mid < \epsilon}$\\

\n where  $\phi_{t}(x) $ is a solution  of {\bf{\color{blue} (1.6)}}\\

\n {\color{red} $\bullet$} We will say that  an equilibrium point $x^{*}$ is {\color{red} asymptotically stable} if  for each neighborhood $\mathbb{U}$ of $x^{*}$ there exists a neighborhood $\mathbb{W}$ such that $x^{*}  \in \mathbb{W} \subset \mathbb{U}$ and $x(0) \in \mathbb{W}$ implies that the solution $\phi_{t}(x)$ satisfies $\phi_{t}(x) \in  \mathbb{U}$ for all $t > 0$, and that $\phi_{t}(x) \longrightarrow  x^{*}$as $t \longrightarrow + \infty $.\\

\n In particular, a system is called {\color{red}asymptotically stable} around its equilibrium point at {\color{red}the origin} if it satisfies the following two conditions:\\

\n 1. Given any $\epsilon > 0; \exists \delta_{1} > 0$ such that if $ \mid\mid x(0)\mid\mid  < \delta_{1}$, then \\

\n  $ \mid\mid \phi_{t}(x) \mid\mid < \epsilon, \forall t > 0$.\\

\n 2. $\exists \delta_{1} > 0$ such that if $ \mid\mid x(0) \mid\mid  < \delta_{2}$, then $\phi_{t}(x) \longrightarrow 0$ as $t \longrightarrow \infty$. \hfill { } $\blacklozenge$\\

\n The first condition requires that the state trajectory can be confined to an arbitrarily small ``ball" centered at the equilibrium point and of radius $\epsilon$, when released from an arbitrary initial condition in a ball of sufficiently small (but positive) radius $\delta_{1}$. This is called stability in the sense of Lyapunov {\color{red}(i.s.L.)}.\\

\n It is possible to have stability in the sense of Lyapunov without having asymptotic stability, in which case we refer to the equilibrium point as marginally stable. Nonlinear systems also exist that satisfy the second requirement without being stable i.s.L. An equilibrium point that is not stable i.s.L. is termed unstable.\\

\n {\color{blue} $\bullet$} {\color{red} {\bf Linear stability analysis for systems of ordinary differential equations}}\\

\n Consider the  n-dimensional dynamical system {\bf{\color{blue}(1.6)}} written in the following form:\\

\n $\displaystyle{\frac{dx_{i}}{dt} = f_{i}(x(t)) } $; \hfill { } {\bf{\color{blue}(1.7)}}\\

\n $\displaystyle{ x(t) =(x_{1}(t),..,x_{i}(t), .., x_{n}(t) , 1 \leq i \leq n } $ and  $\displaystyle{ 0 \leq t < +\infty}$ \\

\n where \\

\n $\displaystyle{ x(0) =(x_{1}(0),..,x_{i}(0), .., x_{n}(0)) = x_{0}}$ is fixed  \\

\n and  \\

\n $\displaystyle{f_{i} : \mathbb{R}^{n} \longrightarrow \mathbb{R}}$ are $\mathcal{C}^{1}$-functions which are given.\\

\n and suppose that $\displaystyle{x^{*} = (x_{1}^{*}, ..., x_{i}^{*}, ...., x_{n}^{*})}$ is a steady state, that is,  $\displaystyle{f_{i}(x^{*}) = 0}$.\\

\n The question of interest is whether the steady state is {\color{red}stable} or {\color{red}unstable}. Consider a small perturbation from the steady state by letting  $\displaystyle{x_{i} = x_{i}^{*} + u_{i} , 1 \leq i \leq n}$ where both $u_{i} , 1 \leq i \n$ are understood to be small. The question of interest translates into the following: will  $\displaystyle{u_{i} , 1 \leq i \n}$ where both grow (so that $\displaystyle{x_{i} , 1 \leq i \leq n}$ move away from the steady state), or will they decay to zero (so that $\displaystyle{x_{i}, , 1 \leq i \leq n}$ move towards the steady state)?\\

\n In the former case, we say that the steady state is unstable, in the latter it is stable.To see whether the perturbation grows or decays, we need to derive differential equations for $u_{i}$ , $1 \leq i \n$  We do so as follows:\\

\n $\displaystyle{\frac{du_{i}}{dt} = \frac{dx_{i}}{dt},  1 \leq i \leq n}$ (since $x_{i}^{*}$ is constant  $1 \leq i \leq n$) \\

\n $\displaystyle{= f_{i}(x) }$ (by definition)\\

\n $\displaystyle{= f_{i}(x^{*} + u), u = (u_{1}, ..., u_{i}, ....u_{n})}$ (substitution)\\

\n $\displaystyle{= f_{i}(x^{*})  + \sum_{j=1}^{n} \frac{\partial f_{i}}{\partial x_{j}}(x^{*}) + ....}$  (Taylor series expansion)\\

\n $\displaystyle{=   \sum_{j=1}^{n} \frac{\partial f_{i}}{\partial x_{j}}(x^{*})  + .... }$ (since $\displaystyle{f_{i}(x^{*}) = 0}$)\\

\n The $....$  denote higher order terms, Since $u_{i}; 1\leq i \leq n$  are assumed to be small, these higher order terms are {\color{red}extremely small}.\\

\n The above linear system for $u_{i}; 1\leq i \leq n$  has the trivial steady state $u_{i} = 0; 1\leq i \leq n$, and the stability of this trivial steady state is determined by the eigenvalues of the matrix, as follows: \\

\n If we can safely neglect the higher order terms, we obtain the following linear system of equations governing the evolution of the perturbations $u_{i}, 1 \leq i \leq n$:\\

\n $\left (\begin{array} {c}\displaystyle{\frac{du_{1}}{dt}}\\
\quad\\
\displaystyle{\frac{du_{2}}{dt}}\\
\quad\\
.\\
\quad\\
.\\
\quad\\
.\\
\quad
.\\
\quad\\
\quad\\
\displaystyle{\frac{du_{n}}{dt}}\\
\end{array}\right)$ $=$ $\left (\begin{array} {cccccccc}\displaystyle{\frac{\partial f_{1}}{\partial x_{1}}(x^{*})} &.& .&\displaystyle{\frac{\partial f_{1}}{\partial x_{j}}(x^{*})} &.&.&.&\displaystyle{\frac{\partial f_{1}}{\partial x_{n}}(x^{*})} \\    
\quad\\
\displaystyle{\frac{\partial f_{2}}{\partial x_{1}}(x^{*})} &.& .&\displaystyle{\frac{\partial f_{2}}{\partial x_{j}}(x^{*})} &.&.&.&\displaystyle{\frac{\partial f_{2}}{\partial x_{n}}(x^{*})} \\   
\quad\\
.\\
\quad
.\\
\quad\\
\displaystyle{\frac{\partial f_{i}}{\partial x_{1}}(x^{*})} &.& .&\displaystyle{\frac{\partial f_{i}}{\partial x_{j}}(x^{*})} &.&.&.&\displaystyle{\frac{\partial f_{i}}{\partial x_{n}}(x^{*})} \\   
\quad\\
.\\
\quad\\                            
\displaystyle{\frac{\partial f_{n}}{\partial x_{1}}(x^{*})} &.& .&\displaystyle{\frac{\partial f_{n}}{\partial x_{j}}(x^{*})} &.&.&.&\displaystyle{\frac{\partial f_{n}}{\partial x_{n}}(x^{*})} \\   
\end{array}\right)$  $\left (\begin{array} {c}\displaystyle{u_{1}}\\
\quad\\
\displaystyle{u_{2}}\\
\quad\\
.\\
\quad\\
.\\
\quad\\
.\\
\quad
.\\
\quad\\
.\\
\quad\\
.
\quad\\
\quad\\
\displaystyle{u_{n}}\\
\end{array}\right )$\\
\quad\\

\n We refer to the matrix as the Jacobian matrix of the original system at the steady state $x^{*}$.\\

\n $\displaystyle{\frac{du}{dt} = \mathbb{J}_{x^{*}}u }$ where $\displaystyle{\mathbb{J}_{x^{*}} = (\frac{\partial f_{i}(x^{*})}{\partial x_{j}})_{1 \leq i, j \leq n}}$.  \hfill { } {\bf{\color{blue}(1.8)}}\\

\n {\bf {\color{red} Theorem 1.5}}\\

\n if the eigenvalues of the Jacobian matrix all have real parts less than zero, then the steady state is stable.\\

\n if the eigenvalues of the Jacobian matrix all have real parts  $< 0$, then  the steady state is asymptotically stable.\\

\n If at least one of the eigenvalues of the Jacobian matrix has real part greater than zero, then the steady state is unstable.\\

\n  Otherwise there is no conclusion (then we have a borderline case between stability and instability; such cases require an investigation of the higher
order terms we neglected, and this requires more sophisticated mathematical machinery discussed in advanced courses on ordinary differential equations). \hfill { } $\blacklozenge$\\

\n {\bf {\color{red} Definition 1.6}}\\

\n An equilibrium point $x^{*}$ is said {\color{red}hyperbolic} if all eigenvalues of the Jacobian matrix have real parts  $\neq 0$.\\

\n {\bf {\color{red} Remark 1.7}}\\

\n A hyperbolic equilibrium point $x^{*}$ is asymptotically stable if the eigenvalues of the Jacobian matrix all have real parts  $< 0$ or  otherwise  it is unstable.\\

\n Let $A$ be the Jacobian matrix, assume  that it is a real hyperbolic matrix, i.e. $\Re e \lambda \neq 0$ for for all eigenvalues  $\lambda$ of $A$, then \\

\n There is a linear change of variables [good coordinates $(x_{{\color{red}s}}, x_{{\color{red}u}})$] that induces a splitting into {\color{red}stable} and  {\color{red}unstable spaces} $\displaystyle{\mathbb{R}^{n} = \mathcal{E}_{s} \oplus \mathcal{E}_{u}}$ so that in the new variables \\

\n $\displaystyle{A = }$ $\left ( \begin{array} {cc}A_{s}&0\\
\quad\\
0&A_{u}\\
\end{array} \right)$
\quad\\

\n and a constant $\alpha >  0$ so that for $t  \geq 0$,\\

\n $\left \{ \begin{array} {c} \displaystyle{\mid\mid e^{tA}x_{{\color{red}s}}\mid\mid \quad \leq e^{-\alpha t}\mid\mid x_{{\color{red}s}}\mid\mid} \quad\\
\quad\\
\displaystyle{\mid\mid e^{- tA}x_{{\color{red}u}}\mid\mid \quad \leq e^{-\alpha t}\mid\mid x_{{\color{red}u}}\mid\mid}\\
\end{array} \right.$
\quad\\

\n We have written $\displaystyle{x_{{\color{red}s}} = P_{{\color{red}s}}x , x_{{\color{red}u}} = P_{{\color{red}u}}x }$ where $\displaystyle{P_{{\color{red}s}} : \mathbb{R}^{n} \longrightarrow \mathcal{E}_{{\color{red}s}}}$ and $\displaystyle{P_{{\color{red}u}} : \mathbb{R}^{n} \longrightarrow \mathcal{E}_{{\color{red}u}}}$ are the orthogonal projections.\\

\n Last but not least, there is a theorem (the Hartman- Grobman Theorem) that guarantees that the stability of the steady state $x^{*}$ of the original system is the
same as the stability of the trivial steady state $0$ of the {\color{red}linearized system}.\\

\n Let $x^{*}$ be an equilibrium point  of nonlinear system {\bf{\color{blue}(1.6)}} then by applying a translation, we can always assume $0$ is a equilibrium point of {\bf{\color{blue}(1.6)}}.\\

\n {\color{red}$\bullet$} Poincaré in his dissertation showed that if $f$ is analytic at the equilibrium point $x^{*}$, and the eigenvalues of $\mathbb{J}_{x^{*}}$ are nonresonant, then
there is a formal power series of change of variable to change {\bf{\color{blue}(1.6)}} to a {\color{red}linear system} {\bf{\color{blue}[Poincaré]}} and {\bf{\color{blue}[Arnold]}}.\\

\n {\color{red}$\bullet$} Hartman and Grobman showed that if $f$ is continuously differentiable, then there is a neighborhood of a hyperbolic equilibrium point and a homeomorphism on this neighborhood, such that the system in this neighborhood is changed to a {\color{red}linear system} under such a homeomorphism  {\bf{\color{blue}[Grobman]}}, {\bf{\color{blue}[Hartman1]}}, {\bf{\color{blue}[Hartman 2]}} ,  {\bf{\color{blue}[Perko,]}}  and {\bf{\color{blue}[X-Wang]}}.\\

\n {\bf {\color{red} Theorem 1.8}} (Hartman-Grobman theorem)  \\

\n Let $\Omega$ be an open set of $\mathbb{R}^{n}$ containing the origin, $\displaystyle{f : \Omega \longrightarrow \mathbb{R}^{n}}$ be a $\mathcal{C}^{1}$- function on $\Omega$, $0 $ be a hyperbolic equilibrium point of the system {\color{blue}{\bf(1.6)}}, and $\displaystyle{ \mathbb{U}_{r}  = \{ x ; \mid\mid x \mid\mid < r\}}$ be the neighborhood of the origin of radius $r$. For any $r, \epsilon > 0$ such that $\displaystyle{ \overline{\mathbb{U}}_{r+\epsilon} \subset \Omega}$, there exists a transformation $y = H(x)$, $H(0) = 0$ and $H$ is a homeomorphism in a neighborhood of $0$, such that the system {\bf{\color{blue}(1.6)}} is changed into the {\color{red}linear system}\\

\n $\displaystyle{y'(t) = \mathbb{A}y,  \quad \mathbb{A} = (\frac{\partial f_{i}(0)}{\partial x_{j}})_{1 \leq i, j \leq n}}$ in $\displaystyle{ \mathbb{U}_{r}}$.\\

\n {\color{red}{\bf Proof}} \\

\n see  {\color{blue}http://www.math.utah.edu/~treiberg/M6414HartmanGrobman.pdf}.\\

\n Thus, the procedure to determine stability of $x^{*}$ is as follows:\\

\n 1. Compute all partial derivatives of the right-hand-side of the original system of differential equations, and construct the Jacobian matrix.\\

\n 2. Evaluate the Jacobian matrix at the steady state.\\

\n 3. Compute eigenvalues.\\

\n 4. Conclude stability or instability based on the real parts of the eigenvalues.\\

\n {\bf {\color{red} Definition 1.9}} (Liapunov function)\\

\n Let $x^{*}$ be an equilibrium point of {\bf {\color{blue}(1.6)}}, $\displaystyle{\mathbb{U}  \subset \Omega}$ be a neighborhood of $x^{*}$ and $\displaystyle{ L : \mathbb{U} \longrightarrow \mathbb{R}}$ be a continuous function. We say that $ L$ is Liapunov function for {\bf {\color{blue}(1.6)}} at $x^{*}$ if \\

\n (1) $\displaystyle{L(x^{*}) = 0}$ and for every $x \neq x^{*}$ we have $\displaystyle{L(x) >  0}$;\\

\n (2) The function $\displaystyle{t \longrightarrow L(\phi_{t}(x)) }$ is decreasing.\\

\n  We say that $ L$ is strictly Liapunov function for {\bf {\color{blue}(1.6)}} at $x^{*}$ if  $L$ satisfy (1) and   \\

\n (3) the function $\displaystyle{t \longrightarrow L(\phi_{t}(x))}$ is strictly decreasing.\\

\n {\bf {\color{red} Remark 1.10}} \\

\n If $L$ is $\mathcal{C}^{1}$ function then we can replace :\\

\n - The condition (2)  by  $\displaystyle{\forall \, x \in \mathbb{U},  < \nabla L(x) , f(x) > \leq 0}$.\\

\n and\\

\n - The condition (3)  by  $\displaystyle{\forall \, x \in \mathbb{U}, < \nabla L(x) , f(x) > < 0}$.\\

\n{\bf {\color{red} Theorem 1.11}} \\

\n If {\bf {\color{blue}(1.6)}} admits a Liapunov function at an equilibrium point $x^{*}$, then $x^{*}$  is stable and if  the Liapunov function is strictly decreasing then $x^{*}$ is asymptotically stable. \hfill { } $\blacklozenge$\\

\n We outline in the next section the Li-Wang's stability criterion {\bf {\color{blue} [Li-Wang]}}  for real matrices and we recall of some spectral properties of M-matrices. In section 3 we give some preliminary definitions and some lemmas for linear stability of {\bf{\color{blue} (covid-19)}} system. In section 4, we present a study of stability of equilibrium points of {\bf{\color{blue} (covid-19)}} system by using  the $\mathcal{R}_{0 }$ criterion and  Li-Wang criterion on {\color{red}second additive compound matrix} associated to Jacobian matrix of system {\bf{\color{blue}(covid-19)}}.\\

\begin{center}
{\color{red}$\S \, 2$ {\bf On Li-Wang's stability criterion of real matrix.}}\\
\end{center}

\n {\bf{\color{red} Definition 2.1}}\\

\n Let $\mathbb{A}$ be an  $n\times n$  matrix and let $\sigma (\mathbb{A})$ be its spectrum. The stability modulus of $\mathbb{A}$ is defined by $\displaystyle{s(\mathbb{A}) = Max\{\mathcal{R}e \lambda ; \lambda \in \sigma (\mathbb{A})\}}$ i.e. s($\mathbb{A}$) is the maximum real part of the eigenvalues of $\mathbb{A}$ called also  the spectral abscissa.\\

\n  $\mathbb{A}$ is said to be stable if $s(\mathbb{A}) < 0$. \hfill { } $\blacklozenge$ \\

\n  The stability of a matrix is related to the Routh-Hurwitz problem on the number of zeros of a polynomial that have negative real parts. Routh-Hurwitz discovered necessary and sufficient conditions for all of the zeros to have negative real parts, which are known today as the Routh-Hurwitz conditions. A good and concise account of the Routh-Hurwitz problem can be found in {\bf{\color{blue}[Banks et al]}}.\\

\n  The Li-Wang criterion offer an alternative to the well-known Routh-Hurwitz. It based on Lozinski$\check{i}$ measures and seconde additive compound matrix. For detailed discussions on compound matrices, the reader is referred to {\bf{\color{blue} [Li-Wang]}}  and for additive compound matrices to {\bf{\color{blue}[Fiedle]}}.\\

\n {\color{red}$\bullet$} In  {\color{blue} [Li-Wang]} a necessary and sufficient condition for the stability of an $n\times n$ matrix with real entries is derived (Li-Wang criterion) by using a simple spectral property of additive compound matrices.\\

\n {\color{red}$\bullet$} A survey is given of a connection between compound matrices and ordinary differential equations  by James S; Muldowney in  {\bf{\color{blue} [Muldowney]} }\\

\n And for an application of Li-Wang criterion, we can consult  {\color{blue}[Diekmann et al]} , {\color{blue}[Intissar et al]} and  {\color{blue}[Khanh]}.\\

 \n Now, let $\mathbb{M}_{n}(\mathbb{K})$ be the linear space of $n\times n$ matrices with entries in $\mathbb{K}$,\\
 
 \n  where $\mathbb{K}$ $=$ $\mathbb{R}$ or $\mathbb{C}$.\\
 
 \n {\bf{\color{red} Definition 2.2}}\\

 \n {\color{red}$\bullet$} Let $\wedge$ denote the exterior product in $\mathbb{K}^{n}$, and let $1 \leq k \leq n$ be an integer. With respect to the canonical basis in the $kth$ exterior product space $\displaystyle{\wedge^{k}\mathbb{K}^{n}}$, the $kth$ additive compound matrix $\mathbb{A}^{[k]}$  of $\mathbb{A}$ is a linear operator on $\displaystyle{\wedge^{k}\mathbb{K}^{n}}$ whose definition on a decomposable element $\displaystyle{x_{1}\wedge x_{2}\wedge ......\wedge x_{k}}$ is\\

\n $\displaystyle{\mathbb{A}^{[k]}x_{1}\wedge x_{2}\wedge ......\wedge x_{k} = \sum_{i=1}^{k}x_{1}\wedge x_{2}\wedge ...\wedge \mathbb{A}x_{i}\wedge...\wedge x_{k}}$ \hfill { } {\color{blue} (2.1)} \\

\n {\color{red}$\bullet$} Let $\displaystyle{\mathbb{A} = (a_{ij})_{1 \leq i,j \leq n}}$ and for any integer $i = 1,...,C_{n}^{k}$, let $\displaystyle{((i)) = (i_{1}, i_{2}, ....,i_{k})}$ be the $ith$ member in the lexicographic ordering of integer $k$-tuples such that $\displaystyle{1 \leq i_{1} < i_{2} < .... < i_{1} \leq n }$  where $\displaystyle{C_{n}^{k}= \frac{n!}{k!(n-k)!}}$. Then \\

\n {\color{red}$\bullet$} The entry  in the $ith$ row and the $jth$ column of $\displaystyle{\mathbb{A}^{[k]} = (\hat{a}_{ij})_{1 \leq i,j \leq C_{n}^{k}}}$ is\\

\n $\hat{a}_{i,j} = $ $ \left\{\begin{array}{c}{\color{red}a_{i_{1},i_{1}} + ..... a_{i_{k},i_{k}}} \quad if \quad ((i)) = ((j)) ; 1 \leq i_{1} < i_{2} \leq n \quad \quad \quad \quad \quad \quad \quad \quad \quad \quad\\
\quad\\
{\color{red}(-1)^{r+s}a_{j_{r},i_{s}}} \quad if \quad exactly \quad one \quad entry\quad of \quad i_{s} \quad does\quad not \quad \quad \quad \quad \quad\\

 occur \quad in\quad ((j))\quad and \quad j_{r} \quad does\quad not\quad occur \quad in \quad((i)),\quad \quad\quad \quad\quad \\
 \quad\\

 {\color{red}0} \quad if\quad ((i)) \quad differs\quad  from \quad ((j))\quad in \quad two\quad or \quad more \quad entries.
 \end{array} \right.$ \hfill { } {\color{blue}(2..2)}\\
.\hfill { } $\blacklozenge$ \\

\n {\color{red}$\bullet$} Let $||. ||$  denote a vector norm in $\mathbb{K}^{n}$ and the operator norm it induces in $\mathbb{M}_{n}(\mathbb{K})$. \\

\n {\color{red} $\circ_{1}$} The Lozinski$\breve{i}$  measure  $\mu$ (also known as logarithmic norm $\mid\mid . \mid\mid_{log}$) on $\mathbb{M}_{n}(\mathbb{K})$ with respect to $||. ||$ is defined by (see {\bf{\color{blue}[Coppel]}}, p. 41)\\

\n For $\mathbb{A} \in \mathbb{M}_{n}(\mathbb{K})$,\\

\n $\displaystyle{\mu(A) : = _{_{_{h \longrightarrow 0^{+}}}}\!\!\!\!\!\!\!\!\!\!\!Lim \frac{||\mathbb{I} + h\mathbb{A} || - 1}{h}}$ $\hfill { } {\color{blue}(2.3)}$ \\

\n {\color{red} $\circ_{2}$}  By the logarithmic norm of a matrix $A$ we mean the real number defined by the formula :\\

\n $\displaystyle{ \mid\mid A \mid\mid_{log} : = _{_{_{t \longrightarrow 0^{+}}}}\!\!\!\!\!\!\!\!\!\!\!Lim \frac{ln ||\mathbb{I} + h\mathbb{A} || - \mid\mid I \mid\mid}{t}}$ $\hfill { } {\color{blue}(2.3)_{bis}}$ \\

\n {\color{red}$\bullet$} The existence of a limit in {\color{blue}$(2.3)_{bis}$} is established on the basis of the convexity of the function $I + tA $ (see [{\color{blue}{\bf [Bylov et al ]}}, Supplement I, Sec. 2], whence we also borrow the notation for the logarithmic norm).\\

\n {\color{red}$\bullet$} The logarithmic norm of a matrix for an arbitrary norm was introduced by the Leningrad mathematician Lozinskii {\color{blue}{\bf [Lozinski]}} and the Swedish mathematician Dahlquist {\color{blue}{\bf [Dahlquist]}} in their papers on the numerical integration of ordinary differential equations. For linear bounded operators in Banach spaces, a similar notion was introduced Daletskii and Krein in their book [{\color{blue}{\bf [Daletskii-Krein]}}, Problems and supplement to Chap. I].\hfill { } $\blacklozenge$\\

\n {\color{red}$\bullet$} Let $A = (a_{ij})$ be a real or complex square $n \times n$ matrix, and let $\lambda_{1}, \lambda_{2}, ...., \lambda_{n}$ be the complete set of its eigenvalues denoted by $\sigma(A)$ (the spectrum of the matrix $A$). The maximal real part of these eigenvalues is denoted  by $s(A)$ i.e. $\displaystyle{s(A) = max_{1\leq i \leq n}\Re e \lambda_{i}}$. (spectral abscissa).\\
\n The term ``spectral abscissa'' (by analogy with the spectral radius $\displaystyle{\rho(A) = \lim \, \mid\mid A^{n} \mid\mid^{\frac{1}{n}}}$ as $n \longrightarrow +\infty$ of a matrix $A$) and the notation for it were proposed in [{\color{blue}{\bf [Perov]}}, p. 23].\\

\n {\color{red}$\bullet$}  What the best  upper and lower bounds for $\displaystyle{\mid\mid e^{tA} \mid\mid ?, 0 \leq t < +\infty }$ where $\displaystyle{e^{tA} = I + \sum_{k=1}^{n}\frac{t^{k}A^{k}}{k!}}$.\\

\n It follows from the definition of $e^{tA}$ that $\displaystyle{e^{-t\mid\mid A\mid\mid} \leq \mid\mid e^{tA}\mid\mid \leq e^{t\mid\mid A\mid\mid}, 0\leq t < + \infty}$ but $- \mid\mid A \mid\mid$ and  $\mid\mid A \mid\mid$  are not the best constants.\\

\n Now, let $\alpha$ and $\beta$ the best constants in the estimate :\\

\n $\displaystyle{e^{t \alpha} \leq \mid\mid e^{tA} \mid\mid\leq e^{\beta t} , 0\leq t < + \infty}$ \hfill { } {\color{blue}(2.4)}\\

\n the existence of such constants is beyond doubt.\\

\n {\bf{\color{red}Theorem 2.3}}\\

\n (i) Let $\alpha$  be the best constant in estimate {\color{blue}(2.4)} from below. Then\\

\n $\displaystyle{\alpha = inf_{0 < t}\frac{ln(\mid\mid e^{tA} \mid\mid)}{t} = lim \, \frac{ln(\mid\mid e^{tA} \mid\mid)}{t} = max_{1 \leq i \leq n}\Re e \lambda_{i}}$ as $t \longrightarrow + \infty$. \hfill { } {\color{blue}(2.5)}\\

\n (ii) Let $\beta$  be the best constant in estimate {\color{blue}(2. 4)} from below. Then\\

\n $\displaystyle{\beta = sup_{0 < t}\frac{ln(\mid\mid e^{tA} \mid\mid)}{t} = lim \ \frac{ln \mid\mid I + tA \mid\mid}{t} = lim \ \frac{ln \mid\mid I + tA \mid\mid - \mid\mid I \mid\mid}{t}}$ as $t \longrightarrow 0^{+}$. \hfill { } {\color{blue}(2.6)}\\

\n {\bf{\color{red}Proof}}\\

\n (i) For the proof, see [{\color{blue}{\bf [Daletskii-Krein]}}, Chap. I, Theorem 4.1]..\\

\n We see from the last equality in {\color{blue}(2.5)} that $\alpha$ is the spectral abscissa of the matrix $A$: $\alpha = s(A)$. Let us stress that the spectral abscissa is independent of the choice of the norm.\\

\n (ii) We see from the last equality in {\color{blue}(2.6)} that $\beta$ is the logarithmic norm of the matrix $A$: $\beta = \mid\mid A\mid\mid_{log}$.\\

\n Consider the logarithmic function on the positive semi-axis. In view of its continuous differentiability, it locally satisfies the Lipschitz condition. Therefore, for any $\epsilon  > 0$, we can indicate a $\delta = \delta_{\epsilon} > 0$such that \\

\n $\displaystyle{|ln u - ln v| \leq  (1 + \epsilon)|u - v|}$  for $\displaystyle{|u - 1| < \delta, |v - 1| < \delta}$.\\

\n Therefore, under the conditions $\displaystyle{\mid \mid\mid e^{tA} \mid\mid - 1\mid < \delta}$ and $\displaystyle{\mid \mid\mid I + tA \mid\mid - 1\mid < \delta}$, , we have\\

\n $\displaystyle{\mid ln \mid\mid exp(tA) \mid\mid-  ln\mid\mid I + tA \mid\mid \leq  (1 + \epsilon) \mid  \mid\mid e^{tA}  \mid\mid -  \mid\mid I + tA  \mid\mid \mid }$\\

\n $\displaystyle{\mid\mid e^{tA} - (I + tA \mid\mid \leq \sum_{k=2}^{n}\frac{t^{k}\mid\mid A \mid\mid^{k}}{k!} = e^{t\mid\mid A \mid\mid} - 1 - t \mid\mid A \mid\mid}$.\\

\n Therefore,\\

\n $\displaystyle{lim_{t\longrightarrow 0^{+}} \, \frac{ln \mid\mid e^{tA} \mid\mid}{t} = lim_{t\longrightarrow 0^{+}} \, \frac{ln \mid\mid I + tA\mid\mid}{t}}$,\hfill { } {\color{blue}(2.7)}\\

\n provided that at least one of the limits in {\color{blue}(2.7)} exists.\\

\n Further, setting $\displaystyle{\epsilon(t) =  \mid\mid I + tA  \mid\mid-  \mid\mid I  \mid\mid}$, we can write \\

\n $\left \{ \begin{array} {c} \displaystyle{\frac{ln \mid\mid I + tA \mid\mid}{t} = \frac{ \mid\mid I + tA \mid\mid - \mid\mid I \mid\mid}{t} \quad \quad \quad  \quad \quad \quad if \, \epsilon(t) = 0}\\
\quad\\
\displaystyle{\frac{ln \mid\mid I + tA \mid\mid}{t} = \frac{ ln(1 + \epsilon(t))}{\epsilon(t)}\frac{\mid\mid I + tA \mid\mid - \mid\mid I \mid\mid}{t} \quad if \, \epsilon(t) \neq 0}\\
\end{array} \right.$\\

\n whence, using the well-known relation $\displaystyle{\frac{ln(1 + x)}{x } \longrightarrow 1}$  as $x \longrightarrow  0$, we obtain\\

\n $\displaystyle{ lim_{t\longrightarrow 0^{+}} \frac{ln \mid\mid I + tA\mid\mid}{t} =   lim_{t\longrightarrow 0^{+}} \frac{ \mid\mid I + tA\mid\mid - \mid\mid I \mid\mid}{t}}$ \hfill { } {\color{blue}(2.8)}\\ 

\n provided that at least one of the limits in {\color{blue}(2.8)} exists. As we have already said above, the last limit exists and serves to define the logarithmic norm. It remains to prove that \\

\n $\displaystyle{sup_{0 < t}\frac{ln \mid\mid e^{tA} \mid\mid}{t} =  lim_{t\longrightarrow 0^{+}} \frac{ln \mid\mid  e^{tA} \mid\mid}{t}}$ \hfill { } {\color{blue}(2.9)}\\ 

\n In {\color{blue}(2.9)}, the quantity on the left exists, is finite and is equal to $\beta$; as proved above, the limit on the right exists, is finite and will be denoted by $b$. The definition of the number $\beta$ implies the inequality $\beta \geq  b$.\\

\n Suppose for the time being that the written inequality is strict: $ \beta > b$. For a sufficiently small $\epsilon > 0$, we can write $\beta  - \epsilon \geq b + \epsilon$ . From the obtained $\epsilon > 0$, we then find a $\delta = \delta_{\epsilon}$ such that\\

\n $\displaystyle{ \mid\mid  e^{tA} \mid\mid \leq \mid\mid  e^{t(b + \epsilon)} \mid\mid }$ for $0 < t \leq \delta$.\\

\n After this, consider an arbitrary fixed $t > 0$. Let us choose a natural number $k$ so that $\displaystyle{0 < \frac{t}{k} \leq \delta}$.\\
After this, we estimate \\

\n $\displaystyle{ \mid\mid  e^{tA} \mid\mid = \mid\mid  e^{k\frac{t}{k}A} \mid\mid \leq \mid\mid  e^{\frac{t}{k}A} \mid\mid^{k} \leq e^{\frac{t}{k}(b + \epsilon)k} = e^{t(b + \epsilon)} \leq e^{t(\beta - \epsilon)}}$ .\\

\n Thus, \\

\n $\displaystyle{ \mid\mid  e^{tA} \mid\mid \leq e^{t(\beta - \epsilon)}}$ for $0 < t < +\infty$,\\

\n and this explicitly contradicts the definition of the number $\beta$.\\

\n This Theorem implies the important inequality :\\

\n $\displaystyle{\alpha \, {\color{red} =} s(A) \, \leq \,  \mid\mid A \mid\mid_{log} \,{\color{red} =}\,  \beta} $.\\

\n For every $A, B \in M_{n}(\mathbb{C})$, $\alpha \geq 0$ , and $\xi \in \mathbb{C}$ the following relations hold:\\

\n {\color{red}$\bullet_{1}$} \quad $\displaystyle{\mu(\alpha A + \xi I) = \alpha\mu(A) + \Re e\xi}$.\\ 

\n {\color{red}$\bullet_{2}$} \quad $\displaystyle{-\mid\mid A \mid\mid \leq -\mu(-A) \leq \mu(A) \leq \mid\mid A \mid\mid}$.\\

\n {\color{red}$\bullet_{3}$} \quad $\displaystyle{\mu(A) + \mu(-A) \geq 0}$\\

\n {\color{red}$\bullet_{4}$} \quad $\displaystyle{\mu(A + B) \leq \mu(A) +\mu(B)}$.\\

\n {\color{red}$\bullet_{5}$} \quad $\displaystyle{- \mu(-A) \leq \Re e \lambda \leq \mu(A)}$ for $\lambda \in \sigma(A)$.\\

\n In the partial case for the Holder vector $p$-norm defined by \\

\n $\displaystyle{\mid\mid x \mid\mid_{p} = (\sum_{i=1}^{n}\mid x_{i} \mid^{p})^{\frac{1}{p}}}$ and $\displaystyle{\mid\mid x \mid\mid_{\infty} = max_{1 \leq i \leq n}\{ \mid x_{i} \mid \}}$ \\

\n then  the corresponding matrix measure can be calculated explicitly in the cases :\\

\n {\color{blue} $\bullet$} {\bf{\color{red}Example 1}}\\

\n {\color{red}$\bullet$ (a)} Let $n \in \mathbb{N}$, $X = (x_{1}, ...., x_{i},..., x_{n}) \in \mathbb{R}^{n}$ with vector norm  $\displaystyle{||X|| = \sum_{i=1}^{n}|x_{i}|}$ and \\ $\mathbb{A} =(a_{i,j}) \in M_{n}(\mathbb{R})$  then $\displaystyle{\mu(\mathbb{A}) =\quad _{_{_{j}}}\!\!\!\!\!\!\!sup(a_{j,j} + \sum_{i,i\neq j}^{n}|a_{i,j}|)}$ is Lozinski$\breve{i}$ norm on $M_{n}(\mathbb{R})$.\\

\n  {\color{red}$\bullet$ (b)} The Lozinski measures of complex matrix $ \mathbb{A} =(a_{i,j}) \in M_{n}(\mathbb{C})$ with respect to the three common norms $\displaystyle{\mid\mid x \mid\mid_{\infty} = sup_{i}\mid x_{i}\mid }$, $\displaystyle{\mid\mid x \mid\mid_{1} =  \sum_{i=1}^{n}\mid x_{i} \mid }$  and $\displaystyle{\mid\mid x \mid\mid_{2} = \sqrt{\sum_{i}^{n}\mid x_{i}\mid^{2} }}$ are\\

\n $\displaystyle{\mu_{\infty}(A) = sup_{i} (\Re e a_{ii }+ \sum_{k, k\neq i}\mid a_{ik}\mid }$, $\displaystyle{\mu_{1}(A) = sup_{k} (\Re e a_{kk }+ \sum_{i, i\neq k}\mid a_{ik}\mid }$ and $\displaystyle{\mu_{2}(A) = s(\frac{A + A^{*}}{2})}$ respectively , where $A^{*}$  denotes the Hermitian adjoint of $A$. \\

\n If $A$ is {\color{red}real symmetric}, then $\mu_{2}(A) = s(A)$. \\

\n For a {\color{red} real matrix} A, conditions $\mu_{\infty}(A)  < 0$ or $\mu_{1}(A) < 0$ can be interpreted as $a_{ii} < 0$ \\

\n for $i =1, . . . , n$, and $A$ is {\color{red}diagonally dominant } in rows or in columns, respectively.  \hfill { } $\blacklozenge$\\\

\n {\color{blue} $\bullet$} {\bf {\color{red}Some upper and lower bounds for the determinant of $n\times n $ matrix $A$ with positive diagonal elements}}\\

\n Let $A  = (a_{ij})_{1\leq i, j \leq n}$ be a real matrix satisfying :\\

\n $\displaystyle{a_{ii} \geq \sum_{j\neq i}\mid a_{ij}\mid , \quad i = 1, 2, ... , n}$  \hfill { }  {\color{blue} (2.4)}\\

\n Then we have the following result :\\

\n {\color{red}{\bf Theorem 2.4}}\\

\n If $A  = (a_{ij})_{1\leq i, j \leq n}$ has elements satisfying  {\color{blue} (2.4)}, it is possible to define $l_{i}$ and $r_{i}$, such that :\\

\n $ \left \{ \begin{array} {c} \displaystyle{ a_{ii} = l_{i} + r_{i}  \quad, 1 \leq i \leq n}\\
\quad\\
\displaystyle{ l_{i} \geq \sum_{j < i }\mid a_{ij} \mid  \quad, 1 \leq i \leq n}\\
\quad\\
\displaystyle{ r_{i} \geq \sum_{j > i }\mid a_{ij} \mid  \quad, 1 \leq i \leq n}\\
\end{array} \right. $ \hfill { } {\color{blue} (2.5)}\\
\quad\\
Then, for any choice of $l_{i}$ and $r_{i}$, satisfying {\color{blue} (2.5)} we have\\

\n $\displaystyle{\sum_{k=0}^{n}(\prod_{i=1}^{k}l_{i}\prod_{i= k+1}^{n}r_{i}) \leq det A \leq \sum_{k=0}^{n}(\prod_{i=1}^{k-1}(l_{i} + 2r_{i})l_{i}\prod_{i=k+1}^{n}r_{i})}$ \hfill { } {\color{blue} (2.6)}\\

\n where an empty product is defined to be $1$ and $det A$ denotes determinant of $A$.\\

\n {\bf {\color{red} Proof}}\\

\n To prove this result, we need the following bound given by Price {\bf {\color{blue}[Price]}} :\\

\n If {\color{blue} (2.4)} holds  then\\

\n $\displaystyle{\prod_{i=1}^{n}(a_{ii} - r_{i}) \leq det A  \leq \prod_{i=1}^{n}(a_{ii} + r_{i})}$ \hfill { } {\color{blue} (2.7)}\\

\n where $\displaystyle{r_{i} = \sum_{j > i}\mid a_{ij} \mid}$.\\

\n Let $\mathbb{D}_{n}$ represent det $A$ then we procceed by induction on $n$:\\

\n {\color{red}(a)} For $n = 2$, let $A = $ $\left ( \begin{array} {cc} a_{11}&a_{12}\\
\quad\\
a_{21}&a_{22}\\
\end{array} \right )$  has elements satisfying  {\color{blue} (2.4)} then $ r_{1} = \mid a_{12}\mid $ and $r_{2} = 0$ then by observing that :\\

\n {\color{red}$\bullet_{1}$} $\displaystyle{\mid a_{12}\mid a_{22} \geq \mid a_{12}\mid \mid a_{21}\mid = \mid a_{12}a_{21}\mid \geq - a_{12} a_{21}}$.\\

\n and\\

\n {\color{red}$\bullet_{2}$} $\displaystyle{\mid a_{12}\mid a_{22} \geq \mid a_{12}\mid \mid a_{21}\mid = \mid a_{12}a_{21}\mid \geq  a_{12} a_{21}}$.\\

\n we deduce that the Price's theorem holds.\\

\n Now let $\mathbb{D}_{2} = $ $\left \vert \begin{array} {cc} a_{11}&a_{12}\\
\quad\\
a_{21}&a_{22}\\
\end{array} \right \vert$  $= $ $\left \vert \begin{array} {cc} l_{1} + r_{1}&a_{12}\\
\quad\\
a_{21}&l_{2} + r_{2}\\
\end{array} \right \vert$ and expanding it by diagonal elements in the following form :\\

\n $\mathbb{D}_{2} = $ $\left \vert \begin{array} {cc} l_{1}&a_{12}\\
\quad\\
0&l_{2}\\
\end{array} \right \vert$ + $\left \vert \begin{array} {cc} l_{1} &0\\
\quad\\
0& r_{2}\\
\end{array} \right \vert$ + $\left \vert \begin{array} {cc} r_{1}&0\\
\quad\\
a_{21}& r_{2}\\
\end{array} \right \vert$ + $\left \vert \begin{array} {cc} r_{1}&a_{12}\\
\quad\\
a_{21}&l_{2}\\
\end{array} \right \vert$.\\

\n Therefore\\

\n $\displaystyle{l_{1}l_{2} +l_{1}r_{2} + r_{1}r_{2} \leq \mathbb{D}_{2} \leq r_{1}r_{2} +l_{1}r_{2} + (l_{1} + 2 r_{1})l_{2}}$,\\

\n since\\

\n $\displaystyle{ 0 \leq }$ $\left \vert \begin{array} {cc} r_{1}&a_{12}\\
\quad\\
a_{21}&l_{2}\\
\end{array} \right \vert$ $\displaystyle{ \leq (r_{1} + a_{12})l_{2} < 2r_{1}l_{2}}$ by  {\bf {\color{blue}(2.5)}} and {\bf {\color{blue}(2.7) }}.\\

\n {\color{red}(b)} Assume that for any matrix of order $n-1$ with elements satisfying  {\bf {\color{blue}(2.5)}}, \\

\n $\displaystyle{\sum_{k=0}^{n-1}(\prod_{i=1}^{k}l_{i}\prod_{i= k+1}^{n_1}r_{i}) \leq \mathbb{D}_{n-1} \leq \sum_{k=0}^{n-1}(\prod_{i=1}^{k-1}(l_{i} + 2r_{i})l_{i}\prod_{i=k+1}^{n-1}r_{i})}$ \hfill { } {\color{blue} (2.8)}\\

\n If $\mathbb{D}_{n} = det A$, where $A = (a_{ij})$ $1\leq i, j \leq n$, and the elements $a_{ij}$ satisfy {\bf {\color{blue}(2.5)}}, partition  $\mathbb{D}_{n}$ as follows :

\n $\displaystyle{\mathbb{D}_{n} = }$ $\left \vert \begin{array} {cc} A_{1}&\hat{a}_{2}\\
\quad\\
\hat{a}_{3}&l_{n} + r_{n}\\
\end{array} \right \vert $.\\
\quad\\

\n where \\

\n $\displaystyle{A_{1} = (a_{ij}) 1 \leq i, j \leq n-1}$, $\hat{a}_{2}$ is the column vector with components $\displaystyle{a_{in} 1 \leq i \leq n-1}$, $\hat{a}_{3}$ is the row vector with components $\displaystyle{a_{nj} 1 \leq j \leq n-1}$, and as in {\bf {\color{blue}(2.5)}} , $l_{n} + r_{n} = a_{nn}$, $\displaystyle{l_{n} \geq \sum_{j=1}^{n-1}\mid a_{nj}}$, $r_{n} \geq 0$.\\

\n Then we can write $\mathbb{D}_{n}$ as the sum of two determinants, i.e.,\\

\n $\displaystyle{\mathbb{D}_{n} = \Delta + r_{n}det A_{1}}$ \hfill { } {\color{blue}{\bf (2.9)}}\\

\n $\displaystyle{ \Delta = }$ $\left \vert \begin{array} {cc} A_{1}&\hat{a}_{2}\\
\quad\\
\hat{a}_{3}& l_{n} \\
\end{array} \right \vert $.\\

\quad\\
But the elements of $\Delta$ satisfy {\color{blue}{\bf (2.4)}}, hence, by {\color{blue}{\bf (2.5)}} and {\color{blue}{\bf (2.7)}} we deduce that \\

\n $\displaystyle{\Delta \geq \prod_{i=1}^{n}(a_{ii} - \hat{r}_{i}) \geq \prod_{i=1}^{n}(a_{ii} - r_{i}) = \prod_{i=1}^{n}l_{i}}$,\\

\n and  \hfill { } {\bf {\color{blue}(2.10)}}\\

\n $\displaystyle{\Delta \leq l_{n}\prod_{i=1}^{n-1}(a_{ii} + \hat{r}_{i}) \leq l_{n}\prod_{i=1}^{n-1}(l_{i} + 2r_{i})}$.\\

\n Also, by inductive assumption, since $A_{1}$, is of order $n - 1$, and,  by {\color{blue}{\bf (2.5)}},\\

\n $\displaystyle{ r_{i} \geq \sum_{ j= i+1}^{n}\mid a_{ij}\mid \geq \sum_{ j= i+1}^{n-1}\mid a_{ij}\mid}$,

\n We have , using {\color{blue}{\bf (2.10)}}, {\color{blue}{\bf (2.9)}} and {\color{blue}{\bf (2.8)}},\\

\n $\displaystyle{\mathbb{D}_{n} \geq \prod_{i=1}^{n}l_{i} + r_{n}\sum_{k=0}^{n-1}(\prod_{i=1}^{k}l_{i}\prod_{i=k+1}^{n-1}r_{i}) = \sum_{k=0}^{n}(\prod_{i=1}^{k}l_{i}\prod_{i=k+1}^{n-1}r_{i})}$,\\

\n and\\

\n $\displaystyle{\mathbb{D}_{n} \leq l_{n}\prod_{i=1}^{n-1}(l_{i} + 2r_{i}) + r_{n}\sum_{k=0}^{n-1}(\prod_{i=1}^{k-1}(l_{i} + 2r_{i})l_{k} \prod_{k+1}^{n-1}r_{i}}$\\

\n $\displaystyle{ = \sum_{k=0}^{n}(\prod_{i=1}^{k-1}(l_{i} + 2r_{i})l_{k}\prod_{i=k+1}^{n}r_{i}}$.\hfill { } $\blacklozenge$\\

\n {\bf{\color{red} Remark 2.5}}\\

\n G.B. Price {\bf {\color{blue}[ Price] }} , A. Ostrowski {\bf{ \color{blue}[Ostrowski1] }} , {\bf {\color{blue}[Ostrowski2] }}, J.L. Brenner {\bf{ \color{blue}[Brenner1] }}, {\bf {\color{blue}[Brenner2] }}and H. Schneider {\bf {\color{blue}[Schneider]}} have  given lower and upper bounds for the absolute value of determinants satisfying more general condition than {\bf {\color{blue}(2.4)}}. \\

\n However, the above theorem is not implied by any of their results.\hfill { } $\blacklozenge$\\

\n {\color{blue}$\bullet$} {\bf {\color{red}Bounds on norms of compound matrices}}\\

\n Let $A$ be a matrix in $M_{n}(\mathbb{C})$, For subsets $\alpha$ and $\beta$ of $\{1,..., n\}$ we denote by $A(\alpha \mid \beta)$ the sub-matrix of $A$ whose rows are indexed by $\alpha$ and whose columns are indexed by $\beta$ in their natural order.\\

\n  Let $k$ be a positive integer, $k \leq n$. we  denote by $C_{k}(A)$ the $k^{th}$ of the matrix $A$, that is, the $\left (\begin{array}{c} n\\\quad\\k\\\end{array}\right )$$\times$ $\left (\begin{array}{c} n\\\quad\\k\\\end{array}\right )$ matrix whose elements are {\color{red} the minors} det $A(\alpha \mid\beta)$ $\displaystyle{\alpha, \beta \subseteq \{1, ..., n\}}$, $\mid \alpha \mid = \mid \beta \mid = k$. We index $C_{k}(A)$ by $\displaystyle{\alpha \subseteq \{1, ..., n\}}$, $\mid \alpha \mid = k$ (ordered lexicographically).\\

\n {\color{blue} $\bullet$}  {\bf {\color{red}Example 2}} if $A \in M_{3}(\mathbb{R})$  and $k = 2$ then :\\

\n $C_{2}(A) = $ $\left ( \begin{array} {ccc} det A(\{1, 2\}|\{1, 2\})& det A(\{1, 2\}|\{1, 3\})& det A(\{1, 2\}|\{2, 3\})\\
\quad\\
det A(\{1, 3\}|\{1, 2\}) &det A(\{1, 3\}|\{1, 3\}) & det A(\{1, 3\}|\{2, 3\})\\
\quad\\
det A(\{2, 3\}|\{1, 2\}) &det A(\{2, 3\}|\{1, 3\}) &det A(\{2, 3\}|\{2, 3\})\\
\end{array}\right ) $\\
\quad\\

\n {\color{blue} $\bullet$}  {\bf {\color{red}Example  3}} if $A \in M_{4}(\mathbb{R})$  and $k = 3$ then :\\

\n $C_{3}(A) = $ $\left ( \begin{array} {cccc}det A({1, 2, 3}|{1, 2, 3})\, \,det A({1, 2, 3}|{1, 2, 4}) \,\, det A({1, 2, 3}|{1, 3, 4})\,\, det A({1, 2, 3}|{2, 3, 4})\\
\quad\\
det A({1, 2, 4}|{1, 2, 3}) \, \,det A({1, 2, 4}|{1, 2, 4}) \, \,det A({1, 2, 4}|{1, 3, 4}) \, \,det A({1, 2, 4}|{2, 3, 4})\\
\quad\\
det A({1, 3, 4}|{1, 2, 3})\, \, det A({1, 3, 4}|{1, 2, 4}) \, \,det A({1, 3, 4}|{1, 3, 4}) \, \,det A({1, 3, 4}|{2, 3, 4})\\
\quad\\
det A({2, 3, 4}|{1, 2, 3}) \, \,det A({2, 3, 4}|{1, 2, 4}) \, \,det A({2, 3, 4}|{1, 3, 4}) \, \,det A({2, 3, 4}|{2, 3, 4})\\
\end{array}\right ) $
\quad\\

\n The most important property of the compound mapping is that it is multiplicative.\\

\n {\bf {\color{red} Lemma 2.6}}. ({\bf{\color{blue}[Marshall et al]}}, Theorem 19.F.2)\\

\n  Let $A$ and $B$ be $n \times n$ matrices and let $1 \leq k \leq n$, then $C_{k}(AB) = C_{k}(A)C_{k}(B)$.\hfill { } $\blacklozenge$\\

\n This property is equivalent to the Binet-Cauchy theorem :\\

\n {\bf {\color{red} Theorem 2.7}} (Binet-Cauchy Theorem)\\

\n Let $A$  be  a $n\times m$ complex matrix, $B$ be  a $m\times l$ complex matrix and {\color{red}$p \leq min \{n, m, l\}$} then $C_{p}(AB) = C_{p}(A)C_{p}(B)$.\hfill { } $\blacklozenge$\\

\n Some  other principal properties of compound matrices are given  in [{\bf{\color{blue}[Aitken]}}, {\bf{\color{blue}[Marcus]}} , {\bf{\color{blue}[Mitrouli et al]}}, {\bf{\color{blue}[Kravvaritis et al]}}] for $A \in M_{n}(\mathbb{C})$  and $p$ an integer, $1 \leq  p \leq  n$:\\

\n in particular, let $A \in M_{n}(\mathbb{C})$ and $ k \leq n$ then we have :\\

\n {\color{red}$\bullet_{1}$} if $A$ is unitary, then $C_{k}(A)$ is unitary.\\

\n {\color{red}$\bullet_{2}$} if $A$ is diagonal, then $C_{k}(A)$ is digonal.\\

\n {\color{red}$\bullet_{3}$} if $A$ is upper (lower) triangular , then $C_{k}(A)$ is uppe (lower) triangular.\\

\n {\color{red}$\bullet_{4}$} $C_{k}(A^{T}) = C_{k}(A)^{T}$.\\

\n {\color{red}$\bullet_{5}$} $\displaystyle{ det (A + I) = 1 + det(A) + \sum_{i=1}^{n-1}tr(C_{i}(A)}$.\\

\n {\color{red}$\bullet_{6}$} if $\displaystyle{\{\lambda_{i},     i= 1 ....., n \}}$ are eigenvalues of $A$ then  the eigenvalues of $A^{[k]}$ are of the  following form :\\

\n $\displaystyle{\{\lambda_{i_{1}} + .... +\lambda_{i_{k}},  1\leq i_{1} <   ..... < i_{k}. \leq n \}}$.\\

\n {\color{red}$\bullet_{7}$} if $\displaystyle{\{\lambda_{i},     i= 1 ....., n \}}$ are eigenvalues of $A$ then  the eigenvalues of $C_{k}(A)$ are of the  following form :\\

\n $\displaystyle{\{\lambda_{i_{1}} .... \lambda_{i_{k}},  1\leq i_{1} <   ..... < i_{k}. \leq n \}}$.\\

\n The main use of compound matrices are their spectral properties which follow from the previous lemma together with the Jordan Canonical Form.\\

\n The compounds of {\color{red}{\bf companion matrices}} can be used to study products of roots of polynomials. An excellent example of this can be seen in  {\bf{\color{blue}[Hong et al]}} where there is an extensive description of the compounds of companion matrices as well as some applications.\\

\n Now, let $\nu$ be a vector norm on $\mathbb{C}^{n}$, and for a positive integer $k$, $k \leq n$, let $\mu$ be a norm on $M_{m}(\mathbb{C})$  where $m = $ $\left (\begin{array}{c} n\\\quad\\k\\\end{array}\right )$ then we have\\

\n {\bf {\color{red}Theorem 2.7}} ({\bf {\color{blue}[Elsner et al]}} Theorem 2.1)\\

\n $\displaystyle{\mu(C_{k}(A)) \leq \theta_{k}(\mu, \nu) max_{\alpha \subseteq \{1, ..., n\}; \mid\alpha\mid = k}\prod_{i \in \alpha}\nu(col_{i}(A))}$.\\

\n where\\

\n $\displaystyle{\theta_{k}(\mu, \nu) = max\{\mu(C_{k}(B)); B \in M_{n}(\mathbb{C}, \nu(col_{i}(B)) = 1, i = 1, ..., , \}}$ \\

\n  with $col_{i}(B)$ denotes the $i^{th}$ column of $B$.\hfill { } $\blacklozenge$\\

\n {\color{blue}$\bullet$} {\color{red}{\bf Some criteria  of stability on matrices given by Li and Wang using the compound matrix and Lozinskii measure.}}\\

\n {\bf{\color{red} Lemma 2.8}} {\bf{\color{blue}[Li-Wang]}}\\

\n {\color{red}(i)} Let $\mu$ be a Lozinski$\breve{i}$ measure. Then $s(\mathbb{A}) \leq \mu(\mathbb{A})$\\

\n {\color{red}(ii)} $\displaystyle{s(A) = inf\{\mu(A);  \mu \quad is \quad a\quad Lozinski\check{i}\quad measure\quad on\quad \mathbb{M}_{n}(\mathbb{K})\}}$ where $\mathbb{K} = \mathbb{R}$ or $\mathbb{C}$.  \hfill { } $\blacklozenge$\\

\n{\bf {\color{red}Proof}}\\

\n See {\bf{\color{blue}[Coopole]}}, p. 41 for {\color{red}(i)}\\

\n  and \\

\n {\bf{\color{blue}  [Li-Wang]}} , p. 252 for {\color{red}(ii)}\\

\n {\bf{\color{red} Proposition 2.9}} (see {\bf{\color{blue}[Li-Wang]}}) \\

\n $\displaystyle{ s(\mathbb{A}) < 0 \iff  s(\mathbb{A}^{[2]}) < 0}$ and $\displaystyle{(-1)^{n}det (\mathbb{A}) > 0}$. \hfill { } $\blacklozenge$\\

\n {\bf {\color{red} Theorem 2.10 }} (see {\bf{\color{blue}[Li-Wang]}})\\

\n  Assume that $\mathbb{A} \in \mathbb{M}_{n}(\mathbb{R})$ and $(-1)^{n}det(\mathbb{A}) > 0$. Then $\mathbb{A}$ is stable if and only if $\mu(A^{[2]}) < 0$ for some Lozinskii measure $\mu$ on $M_{\frac{n(n-1)}{2}}(\mathbb{R})$  \hfill { } $\blacklozenge$\\

\n  {\bf{\color {red} Corollary 2.11}}\\

\n  Assume that $\mathbb{A} \in \mathbb{M}_{n}(\mathbb{R})$ and $(-1)^{n}det(\mathbb{A}) > 0$. Then $\mathbb{A}$ is stable if the following conditions are verified:\\
\n {\bf {\color{red}$\bullet$}} $\displaystyle{\hat{a}_{j,j} + \sum_{i=1; i\neq j}^{\frac{n(n-1)}{2}}\mid\hat{a}_{i,j}\mid < 0 , \forall \quad j = 1, .... \frac{n(n-1)}{2}}$\\
\n where $(\hat{a}_{i,j})_{i,j = 1,...\frac{n(n-1)}{2}}$ are the entries of seconde additive compound matrix $A^{[2]}$. \hfill { } $\blacklozenge$\\

\n {\bf {\color{red} Proof}}\\

\n If we take as Lozinski$\breve{i}$'s measure $\displaystyle{\mu(A^{[2]}) = sup_{j}(\hat{a}_{j,j} + \sum_{i=1; i\neq j}^{\frac{n(n-1)}{2}}\mid \hat{a}_{i,j}\mid) , \forall \quad j = 1, .... \frac{n(n-1)}{2}}$\\

\n  then  by applying the above theorem $\displaystyle{\mu(A^{[2]}) < 0}$ and $\mathbb{A}$ is stable. \hfill { } $\blacklozenge$\\ 

\n {\color{red}{\bf Definition 2.12}}\\

\n A matrix $A = (a_{ij}); 1\leq i, j \leq n$ is said to have dominant principal diagonal if \\

\n $\displaystyle{\mid a_{ii} \mid > \sum_{k\neq i}^{n}\mid a_{ik} \mid}$ for  each $1 \leq i \leq n$. \hfill { } {\color{blue}($\star$)}\\

\n {\color{red} {\bf Lemma 2.13}}\\

\n  Let $A$ be a square real or complex matrix such that :\\

\n $\displaystyle{\mid a_{ii} \mid > \sum_{k\neq i}^{n}\mid a_{ik} \mid}$ for  each $1 \leq i \leq n$. \\

\n Then $A$ is invertible and the set of  its eigenvalues is included in  $\displaystyle{\bigcup_{i=1}^{n}\{z \in \mathbb{C}; \mid z - a_{ii} \mid \leq \sum_{k\neq i}^{n}\mid a_{ik} \mid\}}$.\\

\n {\color{red}{\bf Proof}}\\

\n Suppose that  $Ax = 0$ admit a solution $x \neq 0$ where $\displaystyle{x = (x_{1}, x_{2}, ...., x_{n})^{T}}$. Let $i_{0}$ such that $\displaystyle{ \mid x_{i_{0}} \mid = max_{1 \leq i \leq n}\mid x_{i} \mid }$.\\

\n  The $i_{0}^{th}$ equation of te system $Ax = 0$ can be written as follow :\\

\n $\displaystyle{\sum_{k=1}^{n} a_{i_{0}k}x_{k} = 0}$ or  $\displaystyle{a_{i_{0}i_{0}}x_{i_{0}} = \sum_{k \neq i_{0}}^{n}-a_{i_{0}k}x_{k}}$ \\

\n But $x_{i_{0}} \neq 0$ then  $\displaystyle{a_{i_{0}i_{0}} = \sum_{k \neq i_{0}}^{n}-a_{i_{0}k}\frac{x_{k}}{x_{i_{0}}}}$ and $\displaystyle{\mid a_{i_{0}i_{0}}\mid = \sum_{k \neq i_{0}}^{n}\mid a_{i_{0}k}\mid \frac{\mid x_{k}\mid}{\mid x_{i_{0}}}\mid  \leq \sum_{k \neq i_{0}}^{n}\mid a_{i_{0}k}\mid }$ which is impossible.\\

\n Now, let $z$ is an eigenvalue of $A$ then $ A- zI$ is not invertible.\\

\n  It follows that it is not  dominant principal diagonal in particular there exists $i$ such that  $\displaystyle{\mid a_{ii} - z\mid \leq \sum_{k = 1, k\neq i}^{n}\mid a_{ik}\mid}$ and  $\displaystyle{z \in \bigcup_{i=1}^{n}\{z \in \mathbb{C}; \mid z - a_{ii} \mid \leq \sum_{k\neq i}^{n}\mid a_{ik} \mid \}}$ \\

\n {\bf{\color{red} Remark 2.14}}\\

\n (i) If $A$ is a matrix with dominant principal diagonal, then $\displaystyle{\rho(I - D^{-1}A) < 1}$ where $D$ is the diagonal of  $A$ and $\displaystyle{\rho(I - D^{-1}A)}$ is the spectral radius of $\displaystyle{I - D^{-1}A }$ which is defined as the maximum of the moduli $\mid  \lambda \mid$ of eigenvalues $\lambda$ of $\displaystyle{I - D^{-1}A }$.\\

\n (ii) $\displaystyle{\mu(A^{[2]}) < 0}$ can be interpreted as $\hat{a}_{j,j} < 0 $ for $j = 1, . . . ,\frac{n(n-1)}{2}$, and $A^{[2]}$ is diagonally dominant in columns.  \hfill { } $\blacklozenge$\\

\n  {\color{blue}$\bullet$} {\color{red}{\bf Positive Definite Matrix}}\\

\n  {\color{red}{\bf Definition 2.15}}\\

\n An $n\times n$ complex matrix $A$ is called positive definite if\\

\n $\displaystyle{\Re e[x^{*}Ax] > 0 }$ \hfill { } {\bf {\color{blue}(2.4)}}\\	

\n for all nonzero complex vectors $x \in \mathbb{C}^{n}$, where $x^{*}$ denotes the conjugate transpose of the vector $x$.\\

\n  In the case of a {\color{red}real matrix} $A$, equation {\bf {\color{blue}(2.4)}} reduces to\\

\n $\displaystyle{\Re e[x^{T}Ax] > 0 }$ \hfill { } {\bf {\color{blue}(2.5)}}\\
	
\n where $x^{T}$ denotes the transpose. \hfill { } $\blacklozenge$\\

\n  {\color{red}$\bullet$} Positive definite matrices are of both theoretical and computational importance in a wide variety of applications. They are used, for example, in optimization algorithms and in the construction of various linear regression models (Johnson 1970).\\
\n  A positive definite matrix has at least one matrix square root. Furthermore, exactly one of its matrix square roots is itself positive definite.\\
\n A necessary and sufficient condition for a {\color{red}complex matrix} $A$ to be positive definite is that the {\color{red}Hermitian part} \\

\n $\displaystyle{ A_{H} = \frac{1}{2}(A+A^{H})}$ 	\hfill { } {\bf {\color{blue}(2.6)}}\\

\n where $A^{H}$ denotes the conjugate transpose, be positive definite.\\

\n  This means that a {\color{red} real matrix} $A$ is positive definite iff the symmetric part \\

 \n $\displaystyle{ A_{S}= \frac{1}{2}(A+A^{T})}$ 	\hfill { } {\bf {\color{blue}(2.7)}}\\

\n where $A^{T}$ is the transpose, is positive definite (Johnson 1970).\\

\n {\color{red}$\bullet$} Confusingly, the discussion of positive definite matrices is often restricted to only Hermitian matrices, or symmetric matrices .\\

\n in the case of real matrices (Pease 1965 {\bf{\color{blue}[ Pease]}}, Johnson 1970 {\bf{\color{blue}[Johnson]}}, Marcus and Minc 1988, p. 182 {\bf{\color{blue}[Marcus et al 1]}}; Marcus and Minc 1992, p. 69{\bf{\color{blue}[Marcus et al 2]}}; Golub and Van Loan 1996, p. 140 {\bf{\color{blue}[ Golub et al]}}). \\
\n A Hermitian (or symmetric) matrix is positive definite iff all its eigenvalues are positive. Therefore, a general complex (respectively, real) matrix is positive definite iff its Hermitian (or symmetric) part has all positive eigenvalues.\\
\n The determinant of a positive definite matrix is always positive, so a positive definite matrix is always nonsingular. \hfill { } $\blacklozenge$\\

\n {\bf {\color{red} Definition 2.16}}\\

\n (1) An  real square matrix $A$ is said Z-matrix if their of diagonal elements are all non-positive.\\

\n (1) An  real square matrix $A$ is said M-matrix if it is Z-matrix and fulfilling one of the conditions of the following theorem of Fiedler and Ptàk {\color{blue}{\bf [Fiedler-Ptàk]}} . \hfill { } $\blacklozenge$\\

\n {\bf {\color{red} Theorem 2.12}} (Fiedler-Ptàk)\\

\n Let $A$ be a Z-matrix. Then the following conditions are equivalent to each other :\\

\n $1^{o}$ There exists a vector $x \geq 0$ such that $Ax > 0$;\\

\n $2^{o}$ there exists a vector $x > 0$ such that $Ax > 0$;\\

\n $3^{o}$  there exists a diagonal matrix  $D$ with positive diagonal elements such that $ADe  > 0$ (here $e$ is the vector whose all coordinates are $1$);\\

\n $4^{o}$ there exists a diagonal matrix  $D$ with positive diagonal elements such that the matrix $W = AD$ is a matrix with dominant positive principal diagonal; \\

\n $5^{o}$ for each diagonal matrix $R$ such that $R \geq A$ the inverse $R^{-1}$ exists and $\rho(R^{-1}(P - A)) < 1$, where $P$ is the diagonal of $A$; \\

\n $6^{o}$ if $B$ is a Z-matrix and $B \geq A$, then $B^{-1}$ exists; \\

\n $7^{o}$ each eigenvalue of $A$ is positive; \\

\n $8^{o}$ all  principal minors of $A$ are positive;\\

\n $9^{o}$ there exists a strictly increasing sequence $\displaystyle{0 \neq M_{1} \subset   M_{2} \subset ......  M_{n} }$ such that the principal minors $det A(M_{i})$ are positive;\\

\n $10^{o}$ there exists  a permutation matrix $\mathcal{P}$ such that $\mathcal{P}A \mathcal{P}^{-1}$ may be written in the form $RS$ where $R$ is a lower triangular matrix positive diagonal elements such that $R$ is a Z-matrix and $S$ is an upper triangular matrix with positive diagonal elements such that $S$ is a Z-matrix;\\

\n $11^{o}$ the inverse $A^{-1}$ exists and $A^{-1} \geq 0$; \\

\n $12^{o}$ the real part of each eigenvalue of $A$ is positive;\\

\n $13^{o}$ for each vector $x \neq 0$ there exists an index $k$ that $x_{k}y_{k} > 0$  for $y = Ax$  \hfill { } $\blacklozenge$\\

\n {\color{red}{\bf Proof}}\\

\n See {\color{blue}{\bf [Fiedler-Ptàk]}} page 387-389.\\ 

\n {\color{red}{\bf Theorem 2.17}} (A caracterization of M-matrices in relation to the definite positive matrices)\\

\n Let $A $ be a square matrix of order $n$ satisfying:  $a_{ii} \geq 0; a_{ij} \leq 0$ if $i \neq j$.  Then the following conditions are equivalent:\\

\n (i) $A^{-1}$ exists and its elements are $\geq 0$.\\

\n (ii) There exists a diagonal matrix $D$ with elements $ > 0$ such that $DA$ is definite positive matrix.  \hfill { } $\blacklozenge$\\

\n {\color{red}{\bf Proof}}\\

\n {\color{red}$(i) \Longrightarrow (ii)$}\\

\n  Let $e$ be the vector having all its components equal to $1$. We define $x$ and $y$ by $Ax = e$ and $A^{* }y = e$. Then  we have $x_{i} > 0$ and $y_{i} > 0 $ for all $i$.\\

\n Let $B$ be the matrix defined by $\displaystyle{b_{ij} = y_{i}a_{ij}x_{j}}$ . Then $\displaystyle{b_{ii} - \sum_{j\neq i} \mid b_{ij}\mid = \sum_{j}b_{ijj} = y_{i} > 0}$ and  $\displaystyle{b_{ii} - \sum_{j\neq i} \mid b_{ji}\mid = \sum_{j}b_{jij} = x_{i} > 0}$.\\

\n Consequently $B$ and $B^{*}$  are strictly dominant diagonal matrices. Then $B + B^{*}$ is also a strictly dominant diagonal matrix and it is a definite positive matrix because it is symmetric. Then  there exists $\alpha > 0$ such that $\displaystyle{ < Bu , u > \geq \alpha \mid\mid u \mid\mid^{2}}$. \\ 

\n Now, let $\displaystyle{d_{i} = \frac{y_{i}}{x_{i}}}$ and $d_{i}$ are the elements the diagonal matrix $D$ then we have\\

\n $\displaystyle{ < DAu , u >  = \sum_{ij}d_{i}a_{ij}u_{i}u_{j} = \sum_{ij}b_{ij} \frac{u_{i}}{x_{i}}\frac{u_{j}}{x_{j}} \geq \alpha \sum_{i}\mid \frac{u_{i}}{x_{i}}\mid^{2} \beta \sum_{i}\mid u_{i}\mid^{2}}$.\\

\n {\color{red}$(ii) \Longrightarrow (i)$}\\

\n If $u \neq 0$ and  $Au = v $ then we have $\displaystyle{\sum_{i}d_{i}u_{i}v_{i} > 0}$.\\

\n Consequently,  there exists $i$  such that $u_{i}v_{i} > 0 $ which entails $(i) $ by applying the property $13^{o}$  of above theorem on the M-matrices.\\

\n {\color{blue}$\bullet$} {\bf {\color{red} The Schur stability criteria of matrices using the additive compound matrix}}\\

\n {\color{red}{\bf Definition 2.18}} (Shur stability) \\

\n A matrix $A$ is said to be Schur stable if  $\rho(A) < 1$, where $\displaystyle{\rho(A) = max\{\mid \lambda \mid ; \lambda \in \sigma(A)\}}$ (the spectral radius of  $A$).\\

\n Consider the $\mathcal{C}^{r} ; r \geq 1$ map:\\

\n $\displaystyle{ x \longrightarrow g(x) ; x \in \mathbb{R}^{n}}$ \hfill { } {\bf{\color{blue}(2.8)}}\\

\n If {\bf{\color{blue}(2.8)}} has a fixed point $x = x^{*}$, that is, $x^{*} = g(x^{*})$, then the {\color{red}linear map} corresponding to ({\bf{\color{blue}(2.8)}} {\color{red}is}\\

\n $\displaystyle{ y \longrightarrow Ay ; y \in \mathbb{R}^{n}}$ \hfill { } {\bf{\color{blue}(2.9)}}\\

\n where $\displaystyle{ A = Dg(x^{*})}$, the Jacobian matrix of  $g$ at $x^{*}$. \hfill { } $\blacklozenge$\\

\n {\color{red}{\bf Lemma 2.19}} (see {\bf {\color{blue}[Liao]}} Theorem 2.1). \\

\n If the matrix $A$ of the system {\bf{\color{blue}(2.9)}} is Schur stable, then the fixed point $x^{*}$ of the system {\bf{\color{blue}(2.9)}} is asymptotically stable. that is the eigenvalues of $A$ have strictly negative real part. \hfill { } $\blacklozenge$\\

\n  {\bf {\color{red}Theorem 2.20}} (see {\bf {\color{blue}[Zhang et al]}}). \\

\n {\color{blue}$\star_{1}$} Let $\displaystyle{B = I + \frac{2}{det(A - I)}(A - I)^{-1})}$, if  $\displaystyle{(-1)^{n} det(B) > 0}$,, then $A$ is Schur stable $\iff$ there exists some Lozinskii measure $\mu$ such that $\displaystyle{\mu(B^{[2]}) < 0}$.\\

\n {\color{blue}$\star_{2}$} Let $A \in M_{n}(\mathbb{R})$, then $\displaystyle{\rho(A) < 1 \iff  \rho(C_{2}(A)) < 1}$ and $det(I - A^{2}) > 0$.\\

\n {\color{blue}$\star_{3}$} Let $A \in M_{n}(\mathbb{R})$, then $\displaystyle{\rho(A) < 1 \iff  \sigma_{1}\sigma_{2} < 1}$ and $det(I - A^{2}) > 0$.\\

\n where  $\displaystyle{\{\sigma_{1}, \sigma_{2}, ....., \sigma_{n}\}}$ are the singular values of $A$, i.e the eigenvalues of the symmetric matrix $\sqrt{A^{*}A}$ such that $\displaystyle{\sigma_{1} \geq \sigma_{2} \geq ...... \geq \sigma_{n} \geq 0}$. \hfill { } $\blacklozenge$\\

\n  In next section we give some preliminary definitions and lemmas for linear stability of above system  and in section 4, we apply last corollary to stability of Covid 19 system.\\
\begin{center}
\n { \color{red} $\S \, 3$ {\bf Some preliminary definitions and lemmas}}\\
\end{center}

\n {\color{red} $\bullet$} Writing the above  {\color{red}five}-dimensional system as follow:\\

\n $\displaystyle{x' (t) = f(x(t))}$ where $\displaystyle{x(t) = (x_{1},x_{2}(t), x_{3}, x_{4}, x_{5})^{T} = (E(t), I(t), C(t), H(t), D(t))^{T}}$ and \\

\n $\displaystyle{f = (f_{1},f_{2}, f_{3}, f_{4}, f_{5})^{T}}$ such that: \\

\n  $\displaystyle{f_{1} : \mathbb{R}^{5} \longrightarrow \mathbb{R}}$ ; $\displaystyle{f_{1}[x(t)] = B - \mu E(t) + \beta_{9}H(t) + ( \beta_{10} - \beta_{1})E(t)I(t)  + \beta_{7}E(t)D(t) }$\\

\n $\displaystyle{f_{2} : \mathbb{R}^{5} \longrightarrow \mathbb{R}}$ ; $\displaystyle{f_{2}[x(t)] =  - (\beta_{2} + \beta_{6} + \beta_{8} + \mu)I(t) + (\beta_{1} -\beta_{10})E(t)I(t) }$\\

\n $\displaystyle{f_{3} : \mathbb{R}^{5} \longrightarrow \mathbb{R}}$ ; $\displaystyle{f_{3}[x(t)] = \beta_{2}I(t)  - (\beta_{5} + \beta_{3} + \mu)C(t) + \beta_{4}H(t)}$\\

\n $\displaystyle{f_{4} : \mathbb{R}^{5} \longrightarrow \mathbb{R}}$ ; $\displaystyle{f_{4}[x(t)] = \beta_{8}I(t) + \beta_{3}C(t)  - (\beta_{4} +  \beta_{9} +  \mu) H} $\\

\n $\displaystyle{f_{5} : \mathbb{R}^{5} \longrightarrow \mathbb{R}}$ ; $\displaystyle{f_{5}[x(t)] = \beta_{6}I(t) +  \beta_{5}C(t) -  \beta_{7}D(t)E(t)}$\\

\n where the <<$^{T}$>> denotes transpose.\\

\n {\color{blue}$\bullet$} {\bf {\color{red} Basic reproduction number}}\\

\n Mathematical modeling can play an important role in helping to quantify possible disease control strategies by focusing on the important aspects of a disease, determining threshold quantities for disease survival, and evaluating the effect of particular control strategies. \\

\n A very important threshold quantity is the basic reproduction number, sometimes called the basic reproductive number or basic reproductive ratio (Heffernan, Smith, $\&$ Wahl, 2005) \n {\color{blue}[Heffernan et al]}, which is usually denoted by $\mathcal{R}_{0}$. \\

\n The epidemiological definition of $\mathcal{R}_{0}$ is the average number of secondary cases produced by one infected individual introduced into a population of susceptible individuals, where an infected individual has acquired the disease, and susceptible individuals are healthy but can acquire the disease.\\
In reality, the value of $\mathcal{R}_{0}$ for a specific disease depends on many variables, such as location and density of population.\\

\n The study of the stability of jacobian matrices of order less than three of a dynamic system  yields a reasonable $\mathcal{R}_{0}$, but for more complex compartmental models, especially those with more infected compartments, the study of the stability is difficult as it relies on the algebraic Routh-Hurwitz conditions for stability of the Jacobian matrix.\\

\n  An alternative method proposed by Diekmann, Heesterbeek, and Metz (1990) in {\color{blue}[Diekmann et al]} and elaborated by van den Driessche and Watmough (2002) {\color{blue}[Driessche et al 1] } gives a way of determining $\mathcal{R}_{0}$ for a compartmental model by using the next generation matrix. \\

\n Here an outline of this method is given, the proofs and further details can be found in van den Driessche and Watmough (2002) {\color{blue}[Driessche et al 1] } and van den Driessche and Watmough (2008) {\color{blue}[Driessche et al 2] }.\\

\n Let $x = (x_{1} ,x _{2}, ...., x_{m}, ... , x_{n} )^{T}$ be the number of individuals in each compartment, where the first $ m < n$ compartments contain infected individuals.\\

\n  Assume that the equilibrium point $x^{*}$  exists and is stable in the absence of disease, {\color{red}and} that the linearized equations for $x_{1} , ..., x_{m}$ at the $x^{*}$ {\color{red}decouple} from the other equations. The assumptions are given in more details in the references cited above.\\

\n Consider these equations written in the form:\\

\n $\displaystyle{\frac{dx_{i}}{dt} = \mathfrak{F}_{i}(x) - \mathcal{V}_{i}(x) , \, 1 \leq i \leq {\color{red}m}}$\\

\n In this splitting,\\

\n $\mathfrak{F}_{i}(x)$ is the rate of appearance of {\color{red}new} infections in compartment $i$, \\

\n and\\

\n $\mathcal{V}_{i}(x)$ is the rate of other transitions between compartment i and other infected compartments.\\

\n It is assumed that  $\displaystyle{\mathfrak{F}_{i} ,  \mathcal{V}_{i} \in \mathcal{C}^{2}}$ and $\displaystyle{\mathfrak{F}_{i}  = 0, \,  {\color{red}m+1} \leq i \leq n}$ \hfill { } $\blacklozenge$\\

\n {\bf {\color{red} Remark 3.1}}\\

\n Let $n = 5$ and  $(x_{1} ,x _{2}, ...., x_{5} )^{T} = (E, I, C, H, D)^{T} $ the compments of our system {\color{blue}Covid-19} then we have   :\\

\n {\color{red}$\bullet$} $\displaystyle{m = {\color{red}2}}$\\

\n {\color{red}$\bullet$} $\left \{ \begin{array}{c}\displaystyle{\mathfrak{F}_{1}(E, I, C, H, D) = \beta_{7}ED + \beta_{10}EI}\\ 
\quad\\
\displaystyle{ \mathfrak{F}_{2}(E, I, C, H, D) = \beta_{1}EI }\quad\quad \quad \quad  \quad\\
\quad\\
\displaystyle{\mathfrak{F}_{i}(E, I, C, H, D) = {\color{red}0} ; 3 \leq i \leq 5}\quad \quad\\
\end{array} \right.$ \hfill { } {\bf {\color{blue}(3.1)}}\\
\quad\\

The Jacobian matrix associated to {\bf {\color{blue}(3.1)}} is  $\displaystyle{\mathbb{F} = }$ {\scriptsize $\left ( \begin{array} {ccccc} \displaystyle{\frac{\partial \mathfrak{F}_{1}}{\partial E}}&\displaystyle{\frac{\partial \mathfrak{F}_{1}}{\partial I}}&\displaystyle{\frac{\partial \mathfrak{F}_{1}}{\partial C}}&\displaystyle{\frac{\partial \mathfrak{F}_{1}}{\partial H}}&\displaystyle{\frac{\partial \mathfrak{F}_{1}}{\partial D}}\\
\quad\\
\displaystyle{\frac{\partial \mathfrak{F}_{2}}{\partial E}}&\displaystyle{\frac{\partial \mathfrak{F}_{2}}{\partial I}}&\displaystyle{\frac{\partial \mathfrak{F}_{2}}{\partial C}}&\displaystyle{\frac{\partial \mathfrak{F}_{2}}{\partial H}}&\displaystyle{\frac{\partial \mathfrak{F}_{2}}{\partial D}}\\
\quad\\
\displaystyle{\frac{\partial \mathfrak{F}_{3}}{\partial E}}&\displaystyle{\frac{\partial \mathfrak{F}_{3}}{\partial I}}&\displaystyle{\frac{\partial \mathfrak{F}_{3}}{\partial C}}&\displaystyle{\frac{\partial \mathfrak{F}_{3}}{\partial H}}&\displaystyle{\frac{\partial \mathfrak{F}_{3}}{\partial D}}\\
\quad\\
\displaystyle{\frac{\partial \mathfrak{F}_{4}}{\partial E}}&\displaystyle{\frac{\partial \mathfrak{F}_{4}}{\partial I}}&\displaystyle{\frac{\partial \mathfrak{F}_{4}}{\partial C}}&\displaystyle{\frac{\partial \mathfrak{F}_{4}}{\partial H}}&\displaystyle{\frac{\partial \mathfrak{F}_{4}}{\partial D}}\\
\quad\\
\displaystyle{\frac{\partial \mathfrak{F}_{5}}{\partial E}}&\displaystyle{\frac{\partial \mathfrak{F}_{5}}{\partial I}}&\displaystyle{\frac{\partial \mathfrak{F}_{5}}{\partial C}}&\displaystyle{\frac{\partial \mathfrak{F}_{5}}{\partial H}}&\displaystyle{\frac{\partial \mathfrak{F}_{5}}{\partial D}}\\
\end{array} \right )$ }\\

\n = $\left ( \begin{array} {ccccc} \displaystyle{\beta_{7}D + \beta_{10}I}&\displaystyle{\beta_{10}E}&\displaystyle{0}\quad &\displaystyle{0}&\displaystyle{\beta_{7}E}\\
\quad\\
\displaystyle{\beta_{1}I}&\displaystyle{\beta_{1}E}&\displaystyle{0}\quad &\displaystyle{0}&\displaystyle{0}\\
\quad\\
0&0&0\quad &0&0\\
\quad\\
0&0&0\quad &0&0\\
\quad\\
0&0&0\quad &0&0\\
\end{array} \right )$  \hfill { } {\bf {\color{blue}(3.2)}}\\
\quad\\

\n and for {\color{red}$m = 2$} we have\\

\n {\color{red}$\bullet$} $\mathbb{F}_{{\color{red}m}} = $ $\left ( \begin{array}{cc}\displaystyle{\beta_{7}D + \beta_{10}I}&\displaystyle{\beta_{10}E}\\
\quad\\
\displaystyle{\beta_{1}I}&\displaystyle{\beta_{1}E}\\
\end{array} \right )$   \hfill { } {\bf {\color{blue}(3.3)}}\\
\quad\\

\n {\color{red}$\bullet$} $\left \{ \begin{array}{c}\displaystyle{\mathcal{V}_{1}(E, I, C, H, D) = - B + \beta_{1}EI - \beta_{9}H + \mu E}\quad\quad\quad\quad\quad\\ 
\quad\\
\displaystyle{\mathcal{V}_{2}(E, I, C, H, D) = ( \beta_{2} +  \beta_{6} +  \beta_{8} + \mu)I +  \beta_{10}EI}\quad\quad\\ 
\quad\\
\displaystyle{\mathcal{V}_{3}(E, I, C, H, D) = -\beta_{2}I +  (\beta_{5} +  \beta_{3} + \mu)C -  \beta_{4}H}\\ 
\quad\\
\displaystyle{\mathcal{V}_{4}(E, I, C, H, D) = -\beta_{8}I - \beta_{3}C +  (\beta_{4} +  \beta_{9} + \mu)H }\\ 
\quad\\
\displaystyle{\mathcal{V}_{5}(E, I, C, H, D) = \beta_{6}I -  \beta_{5}C +  \beta_{7}DE}\quad\quad\quad\quad\quad\quad\\ 

\end{array} \right.$ \hfill { } {\bf {\color{blue}(3.4)}}\\

The Jacobian matrix associated to {\bf {\color{blue}(3.3)}} is  $\displaystyle{\mathbb{V} = }$ {\scriptsize $\left ( \begin{array} {ccccc} \displaystyle{\frac{\partial \mathcal{V}_{1}}{\partial E}}&\displaystyle{\frac{\partial \mathcal{V}_{1}}{\partial I}}&\displaystyle{\frac{\partial \mathcal{V}_{1}}{\partial C}}&\displaystyle{\frac{\partial\mathcal{V}_{1}}{\partial H}}&\displaystyle{\frac{\partial \mathcal{V}_{1}}{\partial D}}\\
\quad\\
\displaystyle{\frac{\partial  \mathcal{V}_{2}}{\partial E}}&\displaystyle{\frac{\partial  \mathcal{V}_{2}}{\partial I}}&\displaystyle{\frac{\partial \mathcal{V}_{2}}{\partial C}}&\displaystyle{\frac{\partial  \mathcal{V}_{2}}{\partial H}}&\displaystyle{\frac{\partial  \mathcal{V}_{2}}{\partial D}}\\
\quad\\
\displaystyle{\frac{\partial  \mathcal{V}_{3}}{\partial E}}&\displaystyle{\frac{\partial  \mathcal{V}_{3}}{\partial I}}&\displaystyle{\frac{\partial  \mathcal{V}_{3}}{\partial C}}&\displaystyle{\frac{\partial  \mathcal{V}_{3}}{\partial H}}&\displaystyle{\frac{\partial  \mathcal{V}_{3}}{\partial D}}\\
\quad\\
\displaystyle{\frac{\partial \ \mathcal{V}_{4}}{\partial E}}&\displaystyle{\frac{\partial  \mathcal{V}_{4}}{\partial I}}&\displaystyle{\frac{\partial  \mathcal{V}_{4}}{\partial C}}&\displaystyle{\frac{\partial  \mathcal{V}_{4}}{\partial H}}&\displaystyle{\frac{\partial  \mathcal{V}_{4}}{\partial D}}\\
\quad\\
\displaystyle{\frac{\partial  \mathcal{V}_{5}}{\partial E}}&\displaystyle{\frac{\partial  \mathcal{V}_{5}}{\partial I}}&\displaystyle{\frac{\partial  \mathcal{V}_{5}}{\partial C}}&\displaystyle{\frac{\partial  \mathcal{V}_{5}}{\partial H}}&\displaystyle{\frac{\partial  \mathcal{V}_{5}}{\partial D}}\\
\end{array} \right )$}\\ 

\n =  $\left ( \begin{array} {ccccc} \displaystyle{\beta_{1}I + \mu}&\displaystyle{\beta_{1}E}&\displaystyle{0}&\displaystyle{-\beta_{9}}&\displaystyle{0}\\
\quad\\
\displaystyle{\beta_{10}I}&\displaystyle{\beta_{2} + \beta_{6} + \beta_{8} + \mu + \beta_{10}E}&\displaystyle{0}&\displaystyle{0}&\displaystyle{0}\\
\quad\\
\displaystyle{0}&\displaystyle{-\beta_{2}}&\displaystyle{\beta_{5}+ \beta_{3} + \mu}&\displaystyle{-\beta_{4}}&\displaystyle{0}\\
\quad\\
\displaystyle{0}&\displaystyle{-\beta_{8}}&\displaystyle{-\beta_{3}}&\displaystyle{\beta_{4} + \beta_{9} + \mu}&\displaystyle{0}\\
\quad\\
\displaystyle{\beta_{7}D}& \displaystyle{-\beta_{6}}&-\displaystyle{\beta_{5}}&\displaystyle{0}&\displaystyle{\beta_{7}E}\\
\end{array} \right) $ \hfill { } {\bf {\color{blue}(3.5)}}\\

\n and for {\color{red}$m = 2$} we have\\

\n {\color{red}$\bullet$} $\mathcal{V}_{{\color{red}m}} = $ $\left ( \begin{array}{cc}\displaystyle{\beta_{1}I + \mu}&\displaystyle{\beta_{1}E}\\
\quad\\
\displaystyle{\beta_{10}I}&\displaystyle{ \beta_{2} + \beta_{6} + \beta_{8}  + \mu  + \beta_{10}E}\\
\end{array} \right )$   \hfill { } {\bf {\color{blue}(3.6)}}\\

\n {\color{blue}$\bullet$} {\bf{\color{red}Important case}}\\

\n Let {\color{red}$\displaystyle{E = \frac{B}{\mu}, I = C = H= D= 0}$ } and {\color{red}$\displaystyle{\alpha = \beta_{2} + \beta_{6} + \beta_{8}  + \mu}$} then  we have :\\

\n \n {\color{red}$\bullet$} $\mathbb{F}_{{\color{red}m}} = $ $\left ( \begin{array}{cc}\displaystyle{0}&\displaystyle{\frac{\beta_{10}B}{\mu}}\\
\quad\\
\displaystyle{0}&\displaystyle{\frac{\beta_{1}B}{\mu}}\\
\end{array} \right )$,  {\color{red}$\bullet$} $\mathcal{V}_{{\color{red}m}} = $ $\left ( \begin{array}{cc}\displaystyle{ \mu}&\displaystyle{\beta_{1}\frac{B}{\mu}}\\
\quad\\
\displaystyle{0}&\displaystyle{ \alpha  + \frac{\beta_{10}B}{\mu}}\\
\end{array} \right )$ ,  {\color{red}$\bullet$} $\mathcal{V}_{{\color{red}m}}^{-1} = $ $\left ( \begin{array}{cc}\displaystyle{\frac{1}{\mu}}&\displaystyle{\frac{-\beta_{1}B}{\mu[\beta_{10}B + \alpha \mu]}}\\
\quad\\
\displaystyle{0}&\displaystyle{ \frac{1}{\alpha  + \frac{\beta_{10}B}{\mu}}}\\
\end{array} \right )$ and \\  {\color{red}$\bullet$} $\mathbb{F}_{{\color{red}m}}\mathcal{V}_{{\color{red}m}}^{-1} = $ $\left ( \begin{array}{cc}\displaystyle{0}&\displaystyle{\frac{\beta_{10}B}{\alpha \mu + \beta_{10}B}}\\
\quad\\
\displaystyle{0}&\displaystyle{\frac{\beta_{1}B}{\alpha \mu + \beta_{10}B}}\\
\end{array} \right )$\\

\n The eigenvalues of  $\mathbb{F}_{{\color{red}m}}\mathcal{V}_{{\color{red}m}}^{-1}  $ are $\lambda_{0} = 0$ and {\color{blue}$\bullet$} {\color{red}$\displaystyle{\mathcal{R}_{0} = \frac{\beta_{1}B}{\alpha \mu + \beta_{10}B}}$} {\color{blue}$\displaystyle{ = \rho(\mathbb{F}_{{\color{red}m}}\mathcal{V}_{{\color{red}m}}^{-1})}$} which is called {\color{blue}effective basic reproduction number}. \\

\n These following figures give the curves of  $\mathcal{R}_{0}$-evolution  with respect $\mu$ as abscissa of step $\Delta \mu = 0.015$ and parameter $\beta_{10}$ but the other parameters are fixed as  in above table.\\
\begin{center}
 \includegraphics[scale=0.27]{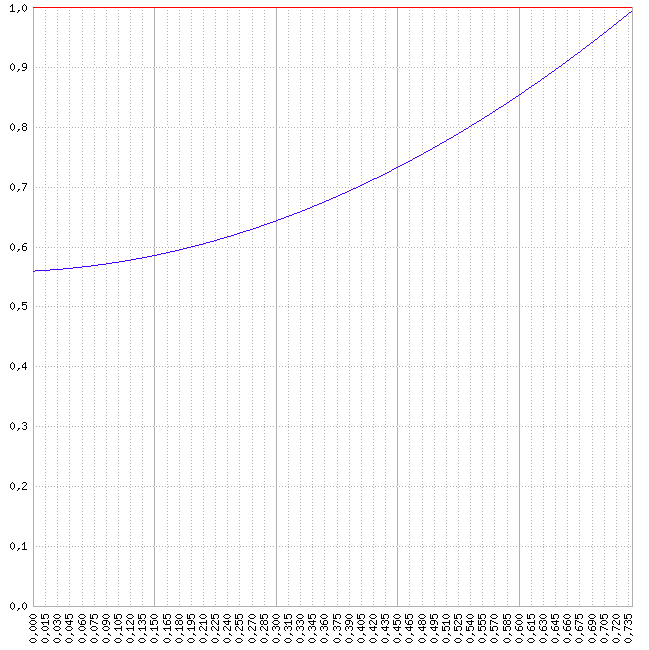} \quad \quad \quad \includegraphics[scale=0.27]{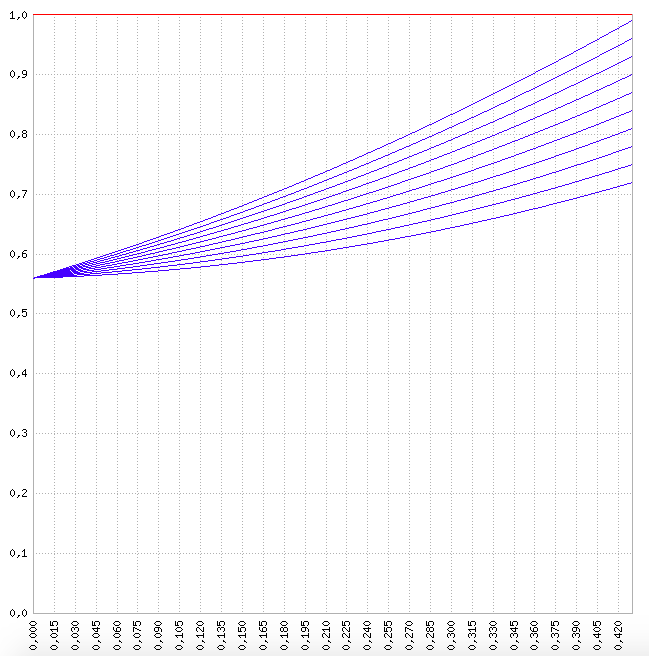} \\
  (a) $ \mu \in [0, 0.74]$; $\beta_{10} = 0.1$    \quad  \quad  \quad \quad  \quad  \quad  \quad \quad  \quad  (b) $ \mu \in [0, 043]$;  $ \beta_{10} \in [0.1, 1]$;  $\Delta \beta_{10} = 0.1$\\
   
  \includegraphics[scale=0.6]{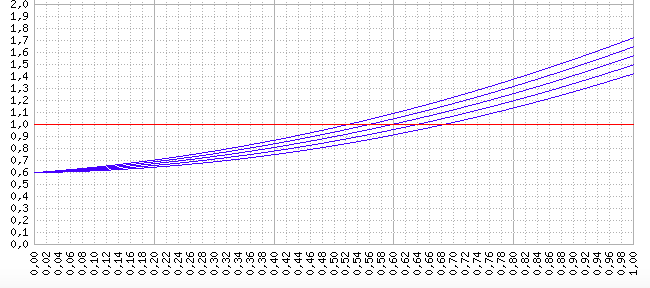}\\
(c) $ \beta_{10} \in [0.1, 1]$;  $\Delta \beta_{10} = 0.1$\\
\end{center}

\n {\bf{\color{red} Remark 3.2}}\\

\n Let $A  = (a_{ij})$ be a $n \times n$ real matrix such that $(a_{ij}) \leq 0$  for all $i \neq j, 1 \leq i,j \leq n$. Then matrix $A$ is also an M-matrix if it can be expressed in the form $A = sI - B$, where $B = (b_{ij}$) with $ b_{ij} \geq 0$, for all $1 \leq i,j \leq n$, where $s$  is at least as large as the maximum of the moduli of the eigenvalues of $B$, and $I$ is an identity matrix.\\

\n For the non-singularity of $A$, according to the Perron-Frobenius theorem, it must be the case that $s > \rho(B)$. Also, for a non-singular M-matrix, the diagonal elements $a_{ii}$ of $A$ must be positive. Here we will further characterize only the class of non-singular M-matrices. \hfill { } $\blacklozenge$\\

\n {\bf{\color{red}Definition 3.3}} (Metzler matrix)\\

\n In mathematics, especially linear algebra, a matrix is called Metzler, quasipositive (or quasi-positive) or essentially nonnegative if all of its elements are non-negative except for those on the main diagonal, which are unconstrained. That is, a Metzler matrix is any matrix A which satisfies $\displaystyle{ A=(a_{ij});\quad a_{ij}\geq 0,\quad i\neq j.}$ \hfill { } $\blacklozenge$\\

\n  $M$-matrices are very useful. We can found some of their applications to ecology, numerical analysis, probability, mathematical programming, game theory,
control theory, and matrix theory. \hfill { } $\blacklozenge$\\

\n {\color{blue}$\bullet$} {\color{red}{\bf Some fondamental properties of M-matrices}} \\
 
\n  An {\color{red}M-}matrix  $A \in M_{n}(\mathbb{R})$ is a matrix of the form $A = \alpha I - B$, where $B \geq 0$ ($B$ is elementwise nonnegative) and $\alpha \geq \rho(B)$. (By the Perron-Frobenius theorem e.g., {\color{blue}{\bf[Intissar arxiv]}}, p. 26, $\rho(B)$, the spectral radius of $B$, is an eigenvalue of $B$.)\\

\n If $A = \alpha I - B$, where $I$ is the identity matrix, $B$ is non-negative and $\alpha > \rho(B)$, then $ A$ is a non-singular M-matrix; if $\alpha = \rho(B)$, then $A$ is a singular M-matrix.\\ There are many definitions of M-matrices equivalent to the above. For example, if a matrix $A$ has the $\mathcal{Z}$ sign pattern and $\rho(A) > 0$, then $A$ is a non-singular M-matrix [4, p. 135 (G20)].\\

\n A matrix of the form $\alpha I - B, \, B \geq 0$ is called a $Z-$matrix.\\

\n {\color{red}$\bullet$}  Observe that a $Z-$matrix  $A $ is an $M-$matrix if and only if $A +\epsilon l $ is nonsingular for all $\epsilon  > 0$.\\

\n We said that a matrix $A =(a_{ij})$ of order $n$ has the $\mathcal{Z}$ sign pattern if $a_{ij}  \leq 0$ for all $i \neq j$.\\

\n if a matrix $A$ has the $\mathcal{Z}$ sign pattern and $\rho(A) > 0$, then $A$ is a non-singular M-matrix [{\bf{\color{blue}[Berman et al]}}, p. 135 (G20)].\\

\n From  Exercise $6b$ of Horn and Johnson [{\bf{\color{blue}[Horn-Johson]}}, p. 127] and  Berman and Plemmons [{\bf{\color{blue}[Berman et al]}}, p. 159 (5.2)], we get the following lemma :\\

\n {\color{red}{\bf Lemma 3.4}}\\

\n Let $A$ be a non-singular M-matrix and suppose $B$ and $BA^{-1}$ have the $\mathcal{Z}$ sign pattern.Then $B$ is a non-singular M-matrix if and only if $BA^{-1}$ is a non-singular M-matrix. \hfill { } $\blacklozenge$\\

\n In general, this lemma does not hold if $B$ a singular M-matrix. It can be shown to hold if $B$ is singular and irreducible. However, this is not sufficient for our needs in part II of this work. we shall need of the following lemma :\\

\n \n {\color{red}{\bf Lemma 3.5}}\\

\n Let $A$ be a non-singular M-matrix and suppose $B \geq 0$.Then,

\n (i) $(A - B)$ is a non-singular M-matrix if and only $\displaystyle{A - B)A^{-1}}$ is a non-singular M-matrix.\\

\n (ii)$(A - B)$ is a non-singular M-matrix if and only $\displaystyle{A - B)A^{-1}}$ is a non-singular M-matrix.\\

\n {\color{red}{\bf Proof}}\\

\n Let $ C = A - B)$. Then both $C$ and $CA^{-1} = I - BA^{-1}$ have the $\mathcal{Z}$ sign pattern. (Recall that $A^{-1} \geq 0$ since $A$ is a non-singular M-matrix).\\ 

\n Hence, the above lemma  implies statement (i). A separate continuity argument can be constructed for each implication in the singular case.\\

\n The following theorem collects conditions that characterize nonsingular $M-$matrices.\\

\n {\color{red}{\bf Theorem 3.6}}\\
 
 \n Let $A = \alpha l - B, \, B \geq 0$. Then the following statements are
equivalent:\\

\n a. $\alpha > \rho(B)$,\\

\n b. $A$ is positive stable: If $\lambda$ is an eigenvalue of $A$, then $\Re e \lambda > 0$,\\

\n c. $A$ is nonsingular and $A^{-1} \geq 0$,\\

\n d. $Ax$ is positive for some positive vector $x$,\\

\n e. The principal minors of $A$ are positive,\\

\n f. The leading principal minors of $A$ are positive. \\

\n {\bf {\color{red} Proof}}\\

\n Conditions (b), (c), and (e) are due to Ostrowski {\color{blue}[Ostrowski]}, who introduced the concept of $M-$matrices. Condition (e) is known in the economics literature as the Hawkins-Simon condition {\color{blue}[Hawkins-Simon]}\\

\n  Condition (d) is due to Schneider {\color{blue}[Schneider]}  and Ky-Fan in {\color{blue}[Fan]} and  the condition (f) to Fiedler and Ptak {\color{blue}[Fiedler-Ptàk]}. \\

\n Many additional characterizations of nonsingular (and of singular) M-matrices are given in {\color{blue}[Berman et al]} .\\ 

\n A subset of the set of all M-matrices that contains the nonsingular M-matrices and whose matrices share many of their properties is the set of group-invertible M-matrices (M-matrices with "property c"). \\

\n Basic reproduction number $\mathcal{R}_{0}$ for the model can be established using the next generation matrix method {\bf {\color{blue}[Cheng-Shan]}} and  {\bf {\color{blue}[Diekmann et al]}}.\\

\n {\bf {\color{red} Definition 3.7}}\\

\n The basic reproduction number $\mathcal{R}_{0}$ is obtained as the {\color{red}spectral radius} of matrix $\mathbb{F}\mathbb{V}^{-1}$ at disease free equilibrium point. Where $\mathbb{F}$ and $\mathbb{V}$ are constructed  as below:\\

\n $\displaystyle{\mathbb{F} = (\frac{\partial \mathfrak{F}_{i} (x^{*})}{\partial x_{j}})_{ij}}$ and $\displaystyle{\mathbb{V} = (\frac{\partial \mathcal{V}_{i} (x^{*})}{\partial x_{j}})_{ij}}$ for $1 \leq i, j \leq m$.\\
For our system the graph of $\mathcal{R}_{0}$ with respect $\frac{1}{\mu}$ is :\\
\begin{center}
\includegraphics[scale=0.25]{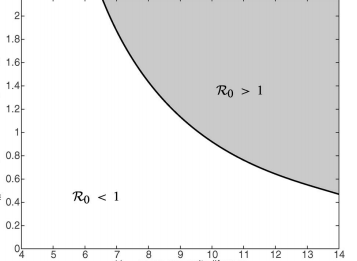}\\
\n Basic reproduction number of infections $\mathcal{R}_{0} $ as a function of $\frac{1}{\mu}$. All other parameters are fixed.\
\end{center}
\n From the above functions $(f_{i}), 1 \leq i \leq 5$ of our system, we consider the associated functions   $({\hat f}_{i}), 1 \leq i \leq 5$ where we delete the linear elements and the negative nonlinear elements, i.e :\\

\n $\displaystyle{f_{1}[x(t)] \longrightarrow {\hat f}_{1}[x(t)] = \beta_{10}E(t)I(t) + \beta_{7}E(t)D(t)}$\\

\n $\displaystyle{f_{2}[x(t)] \longrightarrow {\hat f}_{2}[x(t)] = \beta_{1}E(t)I(t)}$\\

\n $\displaystyle{f_{3}[x(t)] \longrightarrow {\hat f}_{3}[x(t)] = 0}$\\

\n $\displaystyle{f_{4}[x(t)] \longrightarrow {\hat f}_{4}[x(t)] = 0}$\\

\n $\displaystyle{f_{5}[x(t)] \longrightarrow {\hat f}_{5}[x(t)] = 0}$\\

\n and the associated functions   $({\hat g}_{i}), 1 \leq i \leq 5$ where we delete  the {\color{red}non negative nonlinear elements} and we  take the opposite of the obtained expression , i.e :\\

\n $\displaystyle{f_{1}[x(t)] \longrightarrow {\hat g}_{1}[x(t)] = -B + \mu E(t) - \beta_{9}H(t) + \beta_{1}E(t)I(t) }$\\

\n $\displaystyle{f_{2}[x(t)] \longrightarrow {\hat g}_{2}[x(t)] =  + (\beta_{2} + \beta_{6} + \beta_{8} + \mu)I(t) + \beta_{10}E(t)I(t) }$\\

\n $\displaystyle{f_{3}[x(t)] \longrightarrow {\hat g}_{3}[x(t)] =  -\beta_{2}I(t)  + (\beta_{5} + \beta_{3} + \mu)C(t) - \beta_{4}H(t)}$\\

\n $\displaystyle{f_{4}[x(t)] \longrightarrow {\hat g}_{4}[x(t)] = -\beta_{8}I(t) - \beta_{3}C(t)  + (\beta_{4} +  \beta_{9} +  \mu) H} $\\

\n $\displaystyle{f_{5}[x(t)] \longrightarrow {\hat g}_{5}[x(t)] = -\beta_{6}I(t) -  \beta_{5}C(t) + \beta_{7}D(t)E(t)}$\\

\n We define  the matrices $\mathbb{F}$ and  $\mathbb{V}$ as follow :\\

\n $\mathbb{F}$ = {\scriptsize $\left (\begin{array} {ccccc}\displaystyle{\frac{\partial {\hat f}_{1}}{\partial x_{1}}}& \displaystyle{\frac{\partial {\hat f}_{1}}{\partial x_{2}}}&\displaystyle{\frac{\partial {\hat f}_{1}}{\partial x_{3}}}&\displaystyle{\frac{\partial {\hat f}_{1}}{\partial x_{4}}}&\displaystyle{\frac{\partial {\hat f}_{1}}{\partial x_{5}}}\\
\quad\\
\displaystyle{\frac{\partial {\hat f}_{2}}{\partial x_{1}}}& \displaystyle{\frac{\partial {\hat f}_{2}}{\partial x_{2}}}&\displaystyle{\frac{\partial {\hat f}_{2}}{\partial x_{3}}}&\displaystyle{\frac{\partial {\hat f}_{2}}{\partial x_{4}}}&\displaystyle{\frac{\partial {\hat f}_{2}}{\partial x_{5}}}\\
\quad\\
\displaystyle{\frac{\partial {\hat f}_{3}}{\partial x_{1}}}& \displaystyle{\frac{\partial {\hat f}_{3}}{\partial x_{2}}}&\displaystyle{\frac{\partial {\hat f}_{3}}{\partial x_{3}}}&\displaystyle{\frac{\partial {\hat f}_{3}}{\partial x_{4}}}&\displaystyle{\frac{\partial {\hat f}_{3}}{\partial x_{5}}}\\
\quad\\
\displaystyle{\frac{\partial {\hat f}_{4}}{\partial x_{1}}}& \displaystyle{\frac{\partial {\hat f}_{4}}{\partial x_{2}}}&\displaystyle{\frac{\partial {\hat f}_{4}}{\partial x_{3}}}&\displaystyle{\frac{\partial {\hat f}_{4}}{\partial x_{4}}}&\displaystyle{\frac{\partial {\hat f}_{4}}{\partial x_{5}}}\\
\quad\\
\displaystyle{\frac{\partial {\hat f}_{5}}{\partial x_{1}}}& \displaystyle{\frac{\partial {\hat f}_{5}}{\partial x_{2}}}&\displaystyle{\frac{\partial {\hat f}_{5}}{\partial x_{3}}}&\displaystyle{\frac{\partial {\hat f}_{5}}{\partial x_{4}}}&\displaystyle{\frac{\partial {\hat f}_{5}}{\partial x_{5}}}\\
\quad\\
\end{array}\right )$} and  $\mathbb{V}$ = {\scriptsize$\left (\begin{array} {ccccc}\displaystyle{\frac{\partial {\hat g}_{1}}{\partial x_{1}}}& \displaystyle{\frac{\partial {\hat g}_{1}}{\partial x_{2}}}&\displaystyle{\frac{\partial {\hat g}_{1}}{\partial x_{3}}}&\displaystyle{\frac{\partial {\hat g}_{1}}{\partial x_{4}}}&\displaystyle{\frac{\partial {\hat g}_{1}}{\partial x_{5}}}\\
\quad\\
\displaystyle{\frac{\partial {\hat g}_{2}}{\partial x_{1}}}& \displaystyle{\frac{\partial {\hat g}_{2}}{\partial x_{2}}}&\displaystyle{\frac{\partial {\hat g}_{2}}{\partial x_{3}}}&\displaystyle{\frac{\partial {\hat g}_{2}}{\partial x_{4}}}&\displaystyle{\frac{\partial {\hat g}_{2}}{\partial x_{5}}}\\
\quad\\
\displaystyle{\frac{\partial {\hat g}_{3}}{\partial x_{1}}}& \displaystyle{\frac{\partial {\hat g}_{3}}{\partial x_{2}}}&\displaystyle{\frac{\partial {\hat g}_{3}}{\partial x_{3}}}&\displaystyle{\frac{\partial {\hat g}_{3}}{\partial x_{4}}}&\displaystyle{\frac{\partial {\hat g}_{3}}{\partial x_{5}}}\\
\quad\\
\displaystyle{\frac{\partial {\hat g}_{4}}{\partial x_{1}}}& \displaystyle{\frac{\partial {\hat g}_{4}}{\partial x_{2}}}&\displaystyle{\frac{\partial {\hat g}_{4}}{\partial x_{3}}}&\displaystyle{\frac{\partial {\hat g}_{4}}{\partial x_{4}}}&\displaystyle{\frac{\partial {\hat g}_{4}}{\partial x_{5}}}\\
\quad\\
\displaystyle{\frac{\partial {\hat g}_{5}}{\partial x_{1}}}& \displaystyle{\frac{\partial {\hat g}_{5}}{\partial x_{2}}}&\displaystyle{\frac{\partial {\hat g}_{5}}{\partial x_{3}}}&\displaystyle{\frac{\partial {\hat g}_{5}}{\partial x_{4}}}&\displaystyle{\frac{\partial {\hat g}_{5}}{\partial x_{5}}}\\
\quad\\
\end{array}\right )$}\\

\n {\bf {\color{red} Remark 3.8}}\\

\n (i) The explicit matrix $\mathbb{F}$ is :\\

\n $\mathbb{F}$ = {\scriptsize $\left (\begin{array} {ccccc}\displaystyle{\frac{\partial {\hat f}_{1}}{\partial x_{1}}}& \displaystyle{\frac{\partial {\hat f}_{1}}{\partial x_{2}}}&\displaystyle{\frac{\partial {\hat f}_{1}}{\partial x_{3}}}&\displaystyle{\frac{\partial {\hat f}_{1}}{\partial x_{4}}}&\displaystyle{\frac{\partial {\hat f}_{1}}{\partial x_{5}}}\\
\quad\\
\displaystyle{\frac{\partial {\hat f}_{2}}{\partial x_{1}}}& \displaystyle{\frac{\partial {\hat f}_{2}}{\partial x_{2}}}&\displaystyle{\frac{\partial {\hat f}_{2}}{\partial x_{3}}}&\displaystyle{\frac{\partial {\hat f}_{2}}{\partial x_{4}}}&\displaystyle{\frac{\partial {\hat f}_{2}}{\partial x_{5}}}\\
\quad\\
\displaystyle{\frac{\partial {\hat f}_{3}}{\partial x_{1}}}& \displaystyle{\frac{\partial {\hat f}_{3}}{\partial x_{2}}}&\displaystyle{\frac{\partial {\hat f}_{3}}{\partial x_{3}}}&\displaystyle{\frac{\partial {\hat f}_{3}}{\partial x_{4}}}&\displaystyle{\frac{\partial {\hat f}_{3}}{\partial x_{5}}}\\
\quad\\
\displaystyle{\frac{\partial {\hat f}_{4}}{\partial x_{1}}}& \displaystyle{\frac{\partial {\hat f}_{4}}{\partial x_{2}}}&\displaystyle{\frac{\partial {\hat f}_{4}}{\partial x_{3}}}&\displaystyle{\frac{\partial {\hat f}_{4}}{\partial x_{4}}}&\displaystyle{\frac{\partial {\hat f}_{4}}{\partial x_{5}}}\\
\quad\\
\displaystyle{\frac{\partial {\hat f}_{5}}{\partial x_{1}}}& \displaystyle{\frac{\partial {\hat f}_{5}}{\partial x_{2}}}&\displaystyle{\frac{\partial {\hat f}_{5}}{\partial x_{3}}}&\displaystyle{\frac{\partial {\hat f}_{5}}{\partial x_{4}}}&\displaystyle{\frac{\partial {\hat f}_{5}}{\partial x_{5}}}\\
\quad\\
\end{array}\right )$} {\color{red} =} $\left (\begin{array} {ccccc}\displaystyle{\beta_{7}D + \beta_{10}I }& \displaystyle{\beta_{10}E}&\displaystyle{0}&\displaystyle{0}&\displaystyle{\beta_{7}E}\\
\quad\\
\displaystyle{\beta_{1}I}& \displaystyle{\beta_{1}E}&\displaystyle{0}&\displaystyle{0}&\displaystyle{0}\\
\quad\\
0& 0&0&0&0\\
\quad\\
0& 0&0&0&0\\
\quad\\
0& 0&0&0&0\\
\quad\\
\end{array}\right )$\\

\n \n (ii) The explicit matrix $\mathbb{V}$ is :\\

\n $\mathbb{V}$ = {\scriptsize$\left (\begin{array} {ccccc}\displaystyle{\frac{\partial {\hat g}_{1}}{\partial x_{1}}}& \displaystyle{\frac{\partial {\hat g}_{1}}{\partial x_{2}}}&\displaystyle{\frac{\partial {\hat g}_{1}}{\partial x_{3}}}&\displaystyle{\frac{\partial {\hat g}_{1}}{\partial x_{4}}}&\displaystyle{\frac{\partial {\hat g}_{1}}{\partial x_{5}}}\\
\quad\\
\displaystyle{\frac{\partial {\hat g}_{2}}{\partial x_{1}}}& \displaystyle{\frac{\partial {\hat g}_{2}}{\partial x_{2}}}&\displaystyle{\frac{\partial {\hat g}_{2}}{\partial x_{3}}}&\displaystyle{\frac{\partial {\hat g}_{2}}{\partial x_{4}}}&\displaystyle{\frac{\partial {\hat g}_{2}}{\partial x_{5}}}\\
\quad\\
\displaystyle{\frac{\partial {\hat g}_{3}}{\partial x_{1}}}& \displaystyle{\frac{\partial {\hat g}_{3}}{\partial x_{2}}}&\displaystyle{\frac{\partial {\hat g}_{3}}{\partial x_{3}}}&\displaystyle{\frac{\partial {\hat g}_{3}}{\partial x_{4}}}&\displaystyle{\frac{\partial {\hat g}_{3}}{\partial x_{5}}}\\
\quad\\
\displaystyle{\frac{\partial {\hat g}_{4}}{\partial x_{1}}}& \displaystyle{\frac{\partial {\hat g}_{4}}{\partial x_{2}}}&\displaystyle{\frac{\partial {\hat g}_{4}}{\partial x_{3}}}&\displaystyle{\frac{\partial {\hat g}_{4}}{\partial x_{4}}}&\displaystyle{\frac{\partial {\hat g}_{4}}{\partial x_{5}}}\\
\quad\\
\displaystyle{\frac{\partial {\hat g}_{5}}{\partial x_{1}}}& \displaystyle{\frac{\partial {\hat g}_{5}}{\partial x_{2}}}&\displaystyle{\frac{\partial {\hat g}_{5}}{\partial x_{3}}}&\displaystyle{\frac{\partial {\hat g}_{5}}{\partial x_{4}}}&\displaystyle{\frac{\partial {\hat g}_{5}}{\partial x_{5}}}\\
\quad\\
\end{array}\right )$}\\

 {\color{red} =} $\left (\begin{array} {ccccc}\displaystyle{ \beta_{1}I  + \mu}& \displaystyle{\beta_{1}E}&\displaystyle{0}&\displaystyle{-\beta_{9}}&\displaystyle{0}\\
\quad\\
\displaystyle{\beta_{10}I}& \displaystyle{\beta_{10}E + \beta_{8} + \beta_{6}  + \beta_{2} + \mu}&\displaystyle{0}&\displaystyle{0}&\displaystyle{0}\\
\quad\\
0& -\beta_{2}&\beta_{3}+ \beta_{5} + \mu&-\beta_{4}&0\\
\quad\\
0& -\beta_{8}&- \beta_{3}&\beta_{9} + \beta_{4} + \mu&0\\
\quad\\
\beta_{7}D& -\beta_{6}&-\beta_{5}&0&\beta_{7}E\\
\quad\\
\end{array}\right )$\\
\quad\\

\n {\bf {\color{red} Lemma 3.9}}\\

\n (i) $\displaystyle{det \mathbb{V} }$ {\color{red}$\displaystyle{= \beta_{7}E}${\Huge{[}}$\displaystyle{(\beta_{1}I + \mu)\alpha(\beta \gamma - \beta_{3}\beta_{4}) - \beta_{10}I\beta_{1}E(\beta \gamma - \beta_{3}\beta_{4}) + \beta_{10}I \beta_{9}(\beta_{2}\beta_{3} + \beta \beta_{8})}${\Huge{]}}}\\

\n where \\

\n  $\displaystyle{ {\color{red}\alpha} = \beta_{10}E + \beta_{8} + \beta_{6} + \beta_{2} + \mu}$\\

\n  $\displaystyle{ {\color{red}\beta} = \beta_{3} + \beta_{5} + \mu}$\\

\n  $\displaystyle{ {\color{red}\gamma} = \beta_{9} + \beta_{4} + \mu}$\\

\n (ii) $\displaystyle{ \mathbb{V}^{-1} = \frac{1}{det \mathbb{V}} (-1)^{i + j} \mathbb{M}^{T}}$ \\

\n where $\mathbb{M}^{T}$ is the transpose of matrix of minor $(M_{ij})$ of $ \mathbb{V}$ \, $1 \leq i, j \leq 5$ \\

\n As the form of the matrix $\mathbb{F}$ is simple  $\mathbb{F}${\color{red} =} $\left (\begin{array} {ccccc}\displaystyle{\beta_{7}D + \beta_{10}I }& \displaystyle{\beta_{10}E}&\displaystyle{0}&\displaystyle{0}&\displaystyle{\beta_{7}E}\\
\quad\\
\displaystyle{\beta_{1}I}& \displaystyle{\beta_{1}E}&\displaystyle{0}&\displaystyle{0}&\displaystyle{0}\\
\quad\\
0& 0&0&0&0\\
\quad\\
0& 0&0&0&0\\
\quad\\
0& 0&0&0&0\\
\quad\\
\end{array}\right )$ then the matrix  $\mathbb{F}\mathbb{V}^{-1}$ has the following form:\\

\n $\mathbb{F}\mathbb{V}^{-1} = \displaystyle{\frac{1}{det\mathbb{V}}(-1)^{i + j}}$ $\left (\begin{array} {cc} A&B\\
\quad\\
C&D\\
\quad\\
\end{array} \right )$  where  $A$ is $2\times 2$ matrix ,  $B$ is $2\times 3 $ matrix , $C = 0$ is $3\times 2$ matrix and   $D = 0$ is $3\times 3$ matrix. \hfill { } $\blacklozenge$\\

\n {\bf {\color{red} Lemma 3.10}}\\

\n  $\displaystyle{\mathbb{F}\mathbb{V}^{-1}}$ $ = \displaystyle{\frac{1}{det\mathbb{V}}(-1)^{i + j}} \left (\begin{array} {ccccc}(\beta_{7}D + \beta_{10}I)M_{11} + \beta_{10}EM_{12} &(\beta_{7}D + \beta_{10}I)M_{21} + \beta_{10}EM_{22}&\displaystyle{*}&\displaystyle{*}&\displaystyle{*}\\
\quad\\
 \beta_{1}IM_{11} + \beta_{1}EM_{12} & \beta_{1}IM_{21} + \beta_{1}EM_{22}&\displaystyle{*}&\displaystyle{*}&\displaystyle{*}\\
\quad\\
0&0&0&0&0\\
\quad\\
0&0&0&0&0\\
\quad\\
0&0&0&0&0\\
\quad\\
\end{array} \right )$\\
\quad\\
\n where \\

\n $\displaystyle{M_{11} }$ {\color{red} $= \alpha \beta_{7}E(\beta \gamma -  \beta_{3} \beta_{4})$}\\

\n $\displaystyle{M_{12} }$ {\color{red}$ = \beta_{1}E(\beta \gamma - \beta_{3}\beta_{4}) - \beta_{9}(\beta_{2}\beta_{3})+ \beta \beta_{8})$}\\

\n $\displaystyle{M_{21}}$ {\color{red}$ = \beta_{1}E(\beta \gamma - \beta_{3}\beta_{4}) - \beta_{9}(\beta_{2}\beta_{3})+ \beta \beta_{8})$}\\

\n and\\

\n $\displaystyle{M_{22}}$ $= {\color{red} \beta_{7}E(\beta_{1}+ \mu)(\beta \gamma -  \beta_{3} \beta_{4})}$  \hfill { } $\blacklozenge$\\

\n In order to simplify the notations and avoid lengthy expressions, we define the parameters:\\

\n $\displaystyle{a = \frac{1}{det\mathbb{V}}[(\beta_{7}D + \beta_{10}I)M_{11} + \beta_{10}EM_{12}]}$ , \\

\n $\displaystyle{b = \frac{-1}{det\mathbb{V}}[(\beta_{7}D + \beta_{10}I)M_{21} + \beta_{10}EM_{22}]}$,\\

\n $\displaystyle{c = \frac{-1}{det\mathbb{V}}[\beta_{1}IM_{11} + \beta_{1}EM_{12}] }$,\\

\n and \\

\n $\displaystyle{ d =  \frac{1}{det\mathbb{V}}[\beta_{1}IM_{21} + \beta_{1}EM_{22}]}$ \\

\n then the eigenvalues of  $\mathbb{F}\mathbb{V}^{-1}$ are $\lambda_{i} ; 1 \leq i \leq 5$ where $\lambda_{1}$  and  $\lambda_{2}$  are the zeros of \\

\n $(a - \lambda)(d - \lambda) - bc = \lambda^{2} -(a + d)\lambda + (ad - bc) = 0$ \\

\n and \\

\n $\lambda_{3} =  \lambda_{4} = \lambda_{5} = 0$.\\

\n Consequently :\\

\n {\bf {\color{red} Lemma 3.11}}\\

\n  $\displaystyle{\mathcal{R}_{0} = }$ {\color{red}$\displaystyle{ \frac{a + d + \sqrt{\Delta}}{2}}$} where $\displaystyle{\Delta = (a + d)^{2} - 4(ad - bc)}$ \hfill { } $\blacklozenge$\\

\begin{center}
\n {\bf {\color{red} $\S \, 4.$ Determination of equilibrium points}}\\
\end{center}

\n {\color{red}{\bf Theorem 4.1}}\\

\n If the control reproduction number $\mathcal{R}_{0}$ is is less than $1$, model {\color{blue}{\bf (covid-19)}} has a unique equilibrium: the disease-free equilibrium (DFE) $\displaystyle{P_{0} = (\frac{B}{\mu}, 0, 0, 0, 0)}$.\\

\n Conversely, if $\mathcal{R}_{0} > 1$ , model {\color{blue}{\bf (covid-19)}} has two equilibria: the DFE and a unique {\color{red}endemic equilibrium}  $\displaystyle{P^{*} = (E^{*}, I^{*}, C^{*}, H^{*}, D^{*}) = (E^{*}, {\color{blue}\hat{\alpha}}H^{*} , {\color{blue}\hat{\beta}}H^{*}, H^{*},{\color{blue}\hat{\gamma}}H^{*})}$ where  $\displaystyle{H^{*}  = \frac{B - \mu E^{*}}{[(\beta_{1} - \beta_{10}){\color{blue}\hat{\alpha}} - \beta_{7}{\color{blue}\hat{\gamma}}]E^{*} - \beta_{9}}}$ and  $\displaystyle{E^{*} = {\color{red}\frac{\mu + \beta_{2} +  \beta_{6} +  \beta_{8} }{ \beta_{1}  -  \beta_{10}}} = \frac{{\color{red}\alpha}}{\beta_{1} - \beta_{10}}}$\\

\n with\\

\n $\displaystyle{{\color{red}\alpha} =  \beta_{2} +  \beta_{6} +  \beta_{8} + \mu}$;\\

\n $\displaystyle{{\color{blue}\hat{\alpha}} =   \frac{\beta_{3} \beta_{4} +(\beta_{4} +  \beta_{3}  + \mu) (\beta_{5} +  \beta_{3}  + \mu)} {\beta_{2}\beta_{3} + \beta_{8}(\beta_{5} + \beta_{3} + \mu)}}$;\\

\n $\displaystyle{{\color{blue}\hat{\beta}} =   \frac{\beta_{8} \beta_{4} + \beta_{2} (\beta_{4} +  \beta_{9}  + \mu)} {\beta_{2}\beta_{3} + \beta_{8}(\beta_{5} + \beta_{3} + \mu)}}$; \\

\n and\\

\n  $\displaystyle{ {\color{blue}\hat{\gamma}} = \frac{\beta_{6}}{\beta_{7}E^{*}}{\color{blue}\hat{\alpha}}+ \frac{\beta_{5}}{\beta_{7}E^{*}}{\color{blue}\hat{\beta}}}$. \hfill { } $\blacklozenge$\\

\n{\bf {\color{red} Proof}}\\

\n  $\left \{ \begin{array} {c} \displaystyle{0 = B - \beta_{1}EI  + \beta_{7}ED + \beta_{9}H  + \beta_{10}EI  - \mu E}\quad \quad {\color{red}(1)}\\
\quad\\
\displaystyle{0 =  \beta_{1}EI  - \beta_{2}I - \beta_{6}I - \beta_{8}I -\beta_{10}EI  - \mu I} \quad \quad \quad {\color{red}(2)}\\
\quad\\
\displaystyle{0 =  \beta_{2}I  - \beta_{5}C - \beta_{3}C + \beta_{4}H  - \mu C} \quad \quad \quad \quad \quad \quad {\color{red}(3)}\\
\quad\\
\displaystyle{0 =  \beta_{3}C  - \beta_{4}H + \beta_{8}I - \beta_{9}H  - \mu H} \quad \quad \quad \quad \quad \quad {\color{red}(4)}\\
\quad\\
\displaystyle{0 =  \beta_{5}C  + \beta_{6}I -  \beta_{7}DE} \quad \quad \quad \quad \quad \quad \quad \quad \quad \quad \quad {\color{red}(5)}\\
\end{array} \right.$ \hfill {} {\bf{\color{blue} (Equilibrium points)}}\\

\n (i) We observe that $\displaystyle{P_{0} = (\frac{B}{\mu}, 0, 0, 0, 0)}$ is an equilibrium point which is called {\color{red}disease free equilibrium point}.\\

\n (ii) A second equilibrium point  $\displaystyle{P^{*} } $ is given by  $\displaystyle{P^{*} = (E^{*}, I^{*}, C^{*}, H^{*}, D^{*}) = (E^{*}, {\color{blue}\hat{\alpha}}H^{*} , {\color{blue}\hat{\beta}}H^{*}, H^{*},{\color{blue}\hat{\gamma}}H^{*})}$\\

\n  where  $\displaystyle{H^{*}  = \frac{B - \mu E^{*}}{[(\beta_{1} - \beta_{10}){\color{blue}\hat{\alpha}} - \beta_{7}{\color{blue}\hat{\gamma}}]E^{*} - \beta_{9}}}$ and  $\displaystyle{E^{*} = {\color{red}\frac{\mu + \beta_{2} +  \beta_{6} +  \beta_{8} }{ \beta_{1}  -  \beta_{10}}}}$ which is called {\color{red}Endemic equilibrium point}.\\

\n In fact, If $I \neq 0$ then from equation  {\color{red}(2)},  we deduce  that $\displaystyle{E^{*} = {\color{red}\frac{\mu + \beta_{2} +  \beta_{6} +  \beta_{8} }{ \beta_{1}  -  \beta_{10}}}}$.\\ 

\n Writing the equations  {\color{red}(3)} and  {\color{red}(4)} in the following form \\

\n $\left \{ \begin{array} {c} \displaystyle{\beta_{2}I  - (\beta_{5} + \beta_{3} +  \mu) C = -  \beta_{4}H} \quad \quad \quad \quad \quad \quad {\color{red}(3)_{bis}}\\
\quad\\
 \displaystyle{ \beta_{8}I  +  \beta_{3}C =  (\beta_{4} +  \beta_{9}  + \mu) H} \quad \quad \quad \quad \quad \quad \quad {\color{red}(4)_{bis}}\\
  \end{array}\right. $\\
  
\n to deduce that  \\

\n $\displaystyle{I  =   \frac{[\beta_{3} \beta_{4} +(\beta_{4} +  \beta_{3}  + \mu) (\beta_{5} +  \beta_{3}  + \mu)]H} {\beta_{2}\beta_{3} + \beta_{8}(\beta_{5} + \beta_{3} + \mu)}}$ $= {\color{blue}\hat{\alpha}}H $ ,\\

\n where   $\displaystyle{{\color{blue}\hat{\alpha} =   \frac{\beta_{3} \beta_{4} +(\beta_{4} +  \beta_{3}  + \mu) (\beta_{5} +  \beta_{3}  + \mu)} {\beta_{2}\beta_{3} + \beta_{8}(\beta_{5} + \beta_{3} + \mu)}}}$\\

\n and\\

\n $\displaystyle{C  =   \frac{[\beta_{8} \beta_{4} + \beta_{2} (\beta_{4} +  \beta_{9}  + \mu)]H} {\beta_{2}\beta_{3} + \beta_{8}(\beta_{5} + \beta_{3} + \mu)}}$ $= {\color{blue}\hat{\beta}}H $\\

\n  where   $\displaystyle{{\color{blue}\hat{\beta} =   \frac{\beta_{8} \beta_{4} + \beta_{2} (\beta_{4} +  \beta_{9}  + \mu)} {\beta_{2}\beta_{3} + \beta_{8}(\beta_{5} + \beta_{3} + \mu)}}}$ \\

\n it follows from equation {\color{red}(5)} that $\displaystyle{D  = \frac{\beta_{6}}{\beta_{7}E}I + \frac{\beta_{5}}{\beta_{7}E}C = [\frac{\beta_{6}}{\beta_{7}E}{\color{blue}\hat{\alpha}}+ \frac{\beta_{5}}{\beta_{7}E}{\color{blue}\hat{\beta}}]H = {\color{blue}\hat{\gamma}} H}$\\

\n  where $\displaystyle{ {\color{blue}\hat{\gamma}} = \frac{\beta_{6}}{\beta_{7}E}{\color{blue}\hat{\alpha}}+ \frac{\beta_{5}}{\beta_{7}E}{\color{blue}\hat{\beta}}}$.\\

\n and from equation {\color{red}(1)}, we deduce that  $\displaystyle{H  = \frac{B - \mu E}{[(\beta_{1} - \beta_{10}){\color{blue}\hat{\alpha}} - \beta_{7}{\color{blue}\hat{\gamma}}]E - \beta_{9}}}$.\\

\n Then we get $\displaystyle{P^{*} = (E^{*}, I^{*}, C^{*}, H^{*}, D^{*}) = (E^{*}, {\color{blue}\hat{\alpha}}H^{*} , {\color{blue}\hat{\beta}}H^{*}, H^{*},{\color{blue}\hat{\gamma}}H^{*})}$ \\

\n where  $\displaystyle{H^{*}  = \frac{B - \mu E^{*}}{[(\beta_{1} - \beta_{10}){\color{blue}\hat{\alpha}} - \beta_{7}{\color{blue}\hat{\gamma}}]E^{*} - \beta_{9}}}$\\

\n  and\\

\n   $\displaystyle{E^{*} = {\color{red}\frac{\mu + \beta_{2} +  \beta_{6} +  \beta_{8} }{ \beta_{1}  -  \beta_{10}}}}$ which is called Endemic equilibrium point.\\

\n {\color{red}{\bf Corollary 4.2}}\\

\n If the parameters $\displaystyle{(\beta_{1}, \beta_{2}, ....., \beta_{10}) }$ satisfy one of the following conditions :\\

\n (i) $\displaystyle{ \beta_{1} < \beta_{10} }$;\\

\n (ii) $\displaystyle{ (\beta_{2} + \beta_{6} + \beta_{8} + \mu){\color{blue}\hat{\alpha}}  <  \beta_{5}{\color{blue}\hat{\beta}} + \beta_{9}  }$.\\

\n Then the model {\color{blue}{\bf (covid-19)}} has a unique equilibrium: \\

\n The disease-free equilibrium (DFE) $\displaystyle{P_{0} = (\frac{B}{\mu}, 0, 0, 0, 0)}$. \hfill { } $\blacklozenge$\\

\n {\bf {\color{red} Definition 4.3}}\\

\n The equilibrium $\displaystyle{P^{*} = (E^{*}, I^{*}, C^{*}, H^{*}, D^{*})}$ is called feasible if its components are positive. \hfill { } $\blacklozenge$\\

\n Thanks to {\color{blue}{\bf [Driessche et al 1]}}., the following result is straightforward.\\

\n {\bf {\color{red} Theorem 4.4}}\\

\n If $\mathcal{ R}_{0} < 1$, the DFE is locally asymptotically stable. If $\mathcal{ R}_{0} > 1$, the DFE is unstable. \hfill { } $\blacklozenge$\\

\n {\color{red}The epidemiological interpretation of Theorem 4.4} is that, {\bf {\color{blue}(covid-19)}} can be eliminated in the population when $\mathcal{ R}_{0} < 1$ if the initial conditions of the dynamical system  {\bf {\color{blue}(covid-19)}} are in the basin of attraction of the DFE $P_{0}$.\\ 

\n The theorem 4.4 shows also that, $\mathcal{R}_{0}$ is a threshold which can determine if the disease will be spread or not. Thus, reducing its value, is a means to mitigate or even eliminate the (covid-19) . It can be therefore important to determine among model parameters those who mostly influence its value.\\

\n Now let ${\color{red}e(t)} = E(t) - E^{*}$,  ${\color{red}i(t)} = I(t) - I^{*}$,  ${\color{red}c(t)} = C(t) - C^{*}$,  ${\color{red}h(t)} = H(t) - H^{*}$ and  ${\color{red}d(t)} = D(t) - D^{*}$ then it is easy  to verify that $e, i, c, h$ and $d$ satisfy the following system of differential equations:\\

\n {\scriptsize $\left \{ \begin{array} {c} \displaystyle{ \frac{d{\color{red}e(t)}}{dt} = [(\beta_{10} - \beta_{1})I^{*} + \beta_{7}D^{*} - \mu]{\color{red}e(t)} + (\beta_{10} - \beta_{1})E^{*}{\color{red}i(t)} + \beta_{7}E^{*}{\color{red}d(t) }+ \beta_{9}{\color{red}h(t)} + (\beta_{10} - \beta_{1}){\color{red}e(t)i(t)} + \beta_{7}{\color{red}e(t)d(t)}}\\
\quad\\
 \displaystyle{ \frac{d{\color{red}i(t)}}{dt} = [(\beta_{1} - \beta_{10} )I^{*} + \beta_{7}D^{*} - \mu]{\color{red}e(t)} + [(\beta_{1} - \beta_{10})E^{*} - (\beta_{2} + \beta_{6} + \beta_{8} + \mu)]{\color{red}i(t)} + (\beta_{1} - \beta_{10}){\color{red}e(t)i(t)}}\quad \quad \quad \quad \quad \quad\\
\quad\\
\displaystyle{\quad \frac{d{\color{red}c(t)}}{dt} = \beta_{2}{\color{red}i(t)} - (\beta_{2} + \beta_{5} + \mu){\color{red}c(t)} + \beta_{4}{\color{red}h(t)}} \quad \quad \quad \quad \quad \quad \quad \quad \quad \quad \quad  \quad \quad \quad \quad \quad  {\color{blue} (4.1)} \quad \quad \quad \quad \quad \quad \quad \quad \quad \quad \quad \quad \quad \quad \\ 
\quad\\
\displaystyle{ \frac{d{\color{red}h(t)}}{dt} = \beta_{8}{\color{red}i(t) }+ \beta_{3}{\color{red}c(t)} + (\beta_{8} - \beta_{4} - \beta_{9} -\mu){\color{red}h(t)}}\quad \quad \quad \quad \quad\quad \quad \quad \quad \quad\quad \quad \quad \quad \quad \quad \quad \quad \quad \quad \quad \quad \quad \quad \quad \quad \quad \quad \quad \\
\quad\\
\displaystyle{ \frac{d{\color{red}d(t)}}{dt} = -\beta_{7}D^{*}{\color{red}e(t) }+ \beta_{6}{\color{red}i(t)} + \beta_{5}{\color{red}c(t)} -\beta_{7}E^{*}{\color{red}d(t)} -\beta_{7}{\color{red}e(t)d(t)}}\quad \quad \quad \quad \quad\quad \quad \quad \quad \quad\quad \quad \quad \quad \quad \quad \quad \quad \quad \quad \quad \quad \quad \quad\\
\end{array} \right.$ }\\

\n with subject to the restriction  $\displaystyle{ e + i + c + h + d \leq \frac{B}{\mu} - [E^{*} + (1 + \alpha + \beta + \gamma)H^{*}]}$.\\

\n the point $p^{*} = (e^{*}, i^{*}, c^{*}, h^{*}, d^{*}) = (0, 0, 0, 0)$ is an equilibrium point of the system {\color{blue}(4.1)}.\\

\n The jacobian matrix of the system {\color{blue}(4.1)} is given by:\\

\n $\displaystyle{\mathbb{J}_{p^{*}}} = $ $\left ( \begin{array} {ccccc} a_{11}&a_{12}&a_{13}&a_{14}&a_{15}\\
\quad\\
a_{21}&a_{22}&a_{23}&a_{24}&a_{25}\\
\quad\\
a_{31}&a_{32}&a_{33}&a_{34}&a_{35}\\
\quad\\
a_{41}&a_{42}&a_{43}&a_{44}&a_{45}\\
\quad\\
a_{51}&a_{52}&a_{53}&a_{54}&a_{55}\\
\end{array} \right )$\\
\quad\\

\n where\\

\n {\color{red}$\bullet$} $a_{11} = (\beta_{10} - \beta_{1})I^{*} + \beta_{7}D^{*} - \mu$ , $a_{12} = (\beta_{10} - \beta_{1})E^{*}$ , $a_{13} = 0$ , $a_{14} = \beta_{9}$ , $a_{15} = \beta_{7}E^{*}$\\

\n {\color{red}$\bullet$} $a_{21} = (\beta_{1} - \beta_{10})I^{*} + \beta_{7}D^{*} - \mu$ , $a_{22} = (\beta_{1} - \beta_{10})E^{*} - (\beta_{2} + \beta_{6} + \beta_{8} + \mu)$ , $a_{23} = 0$ , $a_{24} = 0$ , $a_{25} = 0$\\

\n {\color{red}$\bullet$} $a_{31} = 0$ ,$a_{32} = \beta_{2}$ ,  $a_{33} = -(\beta_{2} + \beta_{5} + \mu)$ , $a_{34} = \beta_{4}$ , $a_{35} = 0$\\

\n {\color{red}$\bullet$} $a_{41} = 0$ ,$a_{42} = \beta_{8}$ , $a_{43} = \beta_{3}$ , $a_{44} = \beta_{8} - \beta_{4} - \beta_{9}  - \mu$ , $a_{45} = 0$\\

\n {\color{red}$\bullet$} $a_{51} = -\beta_{7}D^{*}$ , $a_{52} = \beta_{6}$,  $a_{53} = \beta_{5}$ , $a_{54} = 0$,  $a_{55} = - \beta_{7}E^{*}$\\

\n i.e. $\displaystyle{\mathbb{J}_{p^{*}}} = $ $\left ( \begin{array} {ccccc} a_{11}&a_{12}&0&a_{14}&a_{15}\\
\quad\\
a_{21}&a_{22}&0&0&0\\
\quad\\
0&a_{32}&a_{33}&a_{34}&0\\
\quad\\
0&a_{42}&a_{43}&a_{44}&0\\
\quad\\
a_{51}&a_{52}&a_{53}&0&a_{55}\\
\end{array} \right )$\\
\quad\\

\n In particular we deduce that\\

\n $\displaystyle{\mathbb{J}_{p^{0}}} = $ $\left ( \begin{array} {ccccc} -\mu& (\beta_{10} - \beta_{1})E&0&\beta_{9}&\beta_{7}E\\
\quad\\
-\mu&(\beta_{1} - \beta_{10})E - (\beta_{2} + \beta_{6} + \beta_{8} + \mu)&0&0&0\\
\quad\\
0&\beta_{2}&-(\beta_{2} + \beta_{5} + \mu)& \beta_{4}&0\\
\quad\\
0&\beta_{8}&\beta_{3}&\beta_{8} - \beta_{4} - \beta_{9}  - \mu&0\\
\quad\\
0&0&0&0&-\beta_{7}E\\
\end{array} \right )$\\

\n where $E = \frac{B}{\mu}$\\

\n Now we recall some technic calculations of determinant of a matrix  in the following form :\\

\n {\color{red}{\bf Lemma 4.5}}\\

\n let $\displaystyle{A = }$ $\left (\begin{array}{rrrr}a_{1,1}&a_{1,2}&a_{1,3}&a_{1,4}\\
\quad\\
a_{2,1}&a_{2,2}&a_{2,3}&a_{2,4}\\
\quad\\
a_{3,1}&a_{3,2}&a_{3,3}&a_{3,4}\\
\quad\\
a_{4,1}&a_{4,2}&a_{4,3}&a_{4,4}\\
\end{array}\right )$ then we have\\

\n {\scriptsize $det$ $\left (\begin{array}{rrrr}a_{1,1}&a_{1,2}&a_{1,3}&a_{1,4}\\
\quad\\
a_{2,1}&a_{2,2}&a_{2,3}&a_{2,4}\\
\quad\\
a_{3,1}&a_{3,2}&a_{3,3}&a_{3,4}\\
\quad\\
a_{4,1}&a_{4,2}&a_{4,3}&a_{4,4}\\
\end{array}\right )$ }$=$ {\scriptsize det $\left (\begin{array}{rr}a_{1,1}&a_{1,2}\\
\quad\\
a_{2,1}&a_{2,2}\\
\end{array} \right )$det $\left (\begin{array}{rr}a_{3,3}&a_{3,4}\\
\quad\\
a_{4,3}&a_{4,4}\\
\end{array} \right )$ - det $\left (\begin{array}{rr}a_{1,1}&a_{1,2}\\
\quad\\
a_{3,1}&a_{3,2}\\
\end{array} \right )$det $\left (\begin{array}{rr}a_{2,3}&a_{2,4}\\
\quad\\
a_{4,3}&a_{4,4}
\end{array} \right )$ \\
\quad\\

\n + det $\left (\begin{array}{rr}a_{1,1}&a_{1,2}\\
\quad\\
a_{4,1}&a_{4,2}\\
\end{array} \right )$det $\left (\begin{array}{rr}a_{2,3}&a_{2,4}\\
\quad\\
a_{3,3}&a_{3,4}\\
\end{array} \right )$}{\scriptsize {\color{red} +} det $\left (\begin{array}{rr}a_{2,1}&a_{2,2}\\
\quad\\
a_{3,1}&a_{3,2}\\
\end{array} \right )$det $\left (\begin{array}{rr}a_{1,3}&a_{1,4}\\
\quad\\
a_{4,3}&a_{4,4}\\
\end{array} \right )$\\
\quad\\

\n - det $\left (\begin{array}{rr}a_{2,1}&a_{2,2}\\
\quad\\
a_{4,1}&a_{4,2}\\
\end{array} \right )$det $\left (\begin{array}{rr}a_{1,3}&a_{1,4}\\
\quad\\
a_{3,3}&a_{3,4}\\
\end{array} \right )$ + det $\left (\begin{array}{rr}a_{3,1}&a_{3,2}\\
\quad\\
a_{4,1}&a_{4,2}\\
\end{array} \right )$det $\left (\begin{array}{rr}a_{1,3}&a_{1,4}\\
\quad\\
a_{2,3}&a_{2,4}
\end{array} \right )$}\\
\quad \\

\n i.e.\\

\includegraphics[scale=0.7]{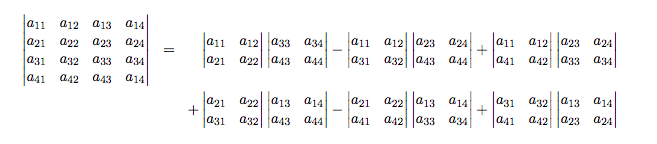}

\n In particular we have:\\

\n  $det$ $\left (\begin{array}{rrrr}a_{1,1}&a_{1,2}&{\color{red}0}&a_{1,4}\\
\quad\\
a_{2,1}&a_{2,2}&{\color{red}0}&{\color{red}0}\\
\quad\\
{\color{red}0}&a_{3,2}&a_{3,3}&{\color{red}0}\\
\quad\\
{\color{red}0}&a_{4,2}&a_{4,3}&a_{4,4}\\
\end{array}\right )$ $=$ \\
\quad\\

\n  det $\left (\begin{array}{rr}a_{1,1}&a_{1,2}\\
\quad\\
a_{2,1}&a_{2,2}\\
\end{array} \right )$det $\left (\begin{array}{rr}a_{3,3}&{\color{red}0}\\
\quad\\
a_{4,3}&a_{4,4}\\
\end{array} \right )$ - det $\left (\begin{array}{rr}a_{1,1}&a_{1,2}\\
\quad\\
{\color{red}0}&a_{3,2}\\
\end{array} \right )$det $\left (\begin{array}{rr}{\color{red}0}&{\color{red}0}\\
\quad\\
a_{4,3}&a_{4,4}
\end{array} \right )$ $+$ \\
 \quad\\
 \quad\\
 
 \n  det $\left (\begin{array}{rr}a_{1,1}&a_{1,2}\\
\quad\\
{\color{red}0}&a_{4,2}\\
\end{array} \right )$det $\left (\begin{array}{rr}{\color{red}0}&{\color{red}0}\\
\quad\\
a_{3,3}&{\color{red}0}\\
\end{array} \right )$ {\color{red} +} det $\left (\begin{array}{rr}a_{2,1}&a_{2,2}\\
\quad\\
{\color{red}0}&a_{3,2}\\
\end{array} \right )$det $\left (\begin{array}{rr}{\color{red}0}&a_{1,4}\\
\quad\\
a_{4,3}&a_{4,4}\\
\end{array} \right )$ $-$\\
\quad\\
\quad\\

\n det $\left (\begin{array}{rr}a_{2,1}&a_{2,2}\\
\quad\\
{\color{red}0}&a_{4,2}\\
\end{array} \right )$det $\left (\begin{array}{rr}{\color{red}0}&a_{1,4}\\
\quad\\
a_{3,3}&{\color{red}0}\\
\end{array} \right )$ + det $\left (\begin{array}{rr}{\color{red}0}&a_{3,2}\\
\quad\\
{\color{red}0}&a_{4,2}\\
\end{array} \right )$det $\left (\begin{array}{rr}{\color{red}0}&a_{1,4}\\
\quad\\
{\color{red}0}&{\color{red}0}\\
\end{array} \right )$\\
\quad\\

\n{\color{red} $\displaystyle{ = a_{11}a_{22}a_{33}a_{44} + a_{21}a_{42}a_{14}a_{33} - (a_{21}a_{12}a_{33}a_{44} + a_{21}a_{32}a_{14}a_{43}) }$}  \hfill { } $\blacklozenge$\\

\n {\color{red}{\bf Corollary 4.6}} \\

\n Let  $\displaystyle{\mathbb{J}_{p^{0}}} = $ $\left ( \begin{array} {ccccc} -\mu& (\beta_{10} - \beta_{1})E^{*}&0&\beta_{9}&\beta_{7}E^{*}\\
\quad\\
-\mu&(\beta_{1} - \beta_{10})E^{*} - (\beta_{2} + \beta_{6} + \beta_{8} + \mu)&0&0&0\\
\quad\\
0&\beta_{2}&-(\beta_{2} + \beta_{5} + \mu)& \beta_{4}&0\\
\quad\\
0&\beta_{8}&\beta_{3}&\beta_{8} - \beta_{4} - \beta_{9}  - \mu&0\\
\quad\\
0&0&0&0&-\beta_{7}E^{*}\\
\end{array} \right )$\\

\n where $E^{*} = \frac{B}{\mu}$\\

\n Then  we have\\

\n (i) $\displaystyle{det \mathbb{J}_{p^{0}} = -\beta_{7}E^{*} det}$ $\left ( \begin{array} {cccc} -\mu& (\beta_{10} - \beta_{1})E^{*}&0&\beta_{9}\\
\quad\\
-\mu&(\beta_{1} - \beta_{10})E^{*} - (\beta_{2} + \beta_{6} + \beta_{8} + \mu)&0&0\\
\quad\\
0&\beta_{2}&-(\beta_{2} + \beta_{5} + \mu)& \beta_{4}\\
\quad\\
0&\beta_{8}&\beta_{3}&\beta_{8} - \beta_{4} - \beta_{9}  - \mu\\
\end{array} \right )$\\

\n $\displaystyle{ =-\mu\beta_{7}E^{*}\{ {\color{red}\beta}(\beta_{8} - {\color{red}\gamma})({\color{red} a}E^{*} -{\color{red}\alpha}) +  {\color{red}\beta}\beta_{8}\beta_{9} + {\color{red}\beta}(\beta_{8}{\color{red} a}E^{*} - {\color{red}\gamma}) + \beta_{2}\beta_{3}\beta_{9}\}}$ \\

\n $\displaystyle{ =-\mu\beta_{7}{\color{red}\beta}E^{*}(2\beta_{8}{\color{red} a}E^{*} + {\color{red} \alpha}{\color{red} \gamma} + \beta_{8}\beta_{9} - \beta_{8}{\color{red} \alpha})}$.\\ 

\n  where  $\displaystyle{{\color{red} a }= \beta_{1} - \beta_{10}}$ , $\displaystyle{{\color{red}\alpha} = \beta_{2} + \beta_{6} + \beta_{8} + \mu}$,  $\displaystyle{{\color{red}\beta} = \beta_{2} + \beta_{5} + \mu}$ and   $\displaystyle{{\color{red}\gamma} = \beta_{4} + \beta_{9} + \mu}$ \\

\n (ii)  Let {\color{red}$\beta_{10} < \beta_{1}$} and {\color{blue}$\displaystyle{ 2\beta_{8}{\color{red} a}E^{*} + {\color{red} \alpha}{\color{red} \gamma} + \beta_{8}\beta_{9} > \beta_{8}{\color{red} \alpha}}$} then {\color{red}$\displaystyle{det \mathbb{J}_{p^{0}} < 0}$}. \hfill { } $\blacklozenge$\\

\n (iii) Under the conditions of (ii) we observe  that a assymption of Li-Wang criterion is satisfied.\\

\n Now, if {\color{red}$\beta_{10} < \beta_{1}$}, we write $\mathbb{J}_{p^{0}}$ in the following form \\

\n  $\displaystyle{\mathbb{J}_{p^{0}} = \mathbb{M} - \mathbb{E}}$ where \\

\n $\displaystyle{\mathbb{M} = }$ $\left ( \begin{array} {ccccc} 0&0&0&\beta_{9}&\beta_{7}E^{*}\\
\quad\\
0&(\beta_{1}- \beta_{10})E^{*}&0&0&0\\
\quad\\
0&\beta_{2}&0&\beta_{4}&0\\
\quad\\
0&\beta_{8}&\beta_{3}&\beta_{8}&0\\
\quad\\
 0&0&0&0&0\\
 \end{array} \right )$ and $\displaystyle{\mathbb{E} = }$ $\left ( \begin{array} {ccccc} \mu&(\beta_{1}- \beta_{10})E^{*}&0&0&0\\
\quad\\
\mu&\alpha&0&0&0\\
\quad\\
0&0&\beta&0&0\\
\quad\\
0&0&0&\gamma&0\\
\quad\\
 0&0&0&0&\beta_{7}E^{*}\\
 \end{array} \right )$ \\
 
 \n Then  if $\displaystyle{E^{*} \neq \frac{\alpha}{\beta_{1} - \beta_{10}}}$ we deduce that :\\
 
 $\displaystyle{\mathbb{E}^{-1} = }$ $\left ( \begin{array} {ccccc} \displaystyle{\frac{\alpha}{\mu(\alpha - {\color{red}a}})}&\displaystyle{\frac{-{\color{red}a}}{\mu(\alpha - {\color{red}a})}}&0&0&0\\
\quad\\
\displaystyle{\frac{-1}{\alpha - {\color{red}a}}}&\displaystyle{\frac{1}{\alpha - {\color{red}a}}}&0&0&0\\
\quad\\
0&0&\displaystyle{\frac{1}{\beta}}&0&0\\
\quad\\
0&0&0&\displaystyle{\frac{1}{\gamma}}&0\\
\quad\\
 0&0&0&0&\displaystyle{\frac{1}{\beta_{7}E^{*}}}\\
 \end{array} \right )$ \\
 
 \n where $\displaystyle{{\color{red}a} = (\beta_{1} - \beta_{10})E^{*}}$\\
 
 \n and\\
 
 \n $\displaystyle{\mathbb{M}\mathbb{E}^{-1} = }$ $\left (  \begin{array} {ccccc} 0&0&0&\displaystyle{\frac{\beta_{9}}{\gamma}}&1\\
 a{\color{red}u}&a{\color{red}v}&0&0&0\\
 \beta_{2}{\color{red}u}&\beta_{2}{\color{red}v}&0&\displaystyle{\frac{\beta_{4}}{\gamma}}&0\\
 \beta_{8}{\color{red}u}&\beta_{8}{\color{red}v}&\displaystyle{\frac{\beta_{3}}{\beta}}&\displaystyle{\frac{\beta_{8}}{\gamma}}&0\\
 0&0&0&0&0\\
 \end{array}\right)$ where $\displaystyle{{\color{red}u} = \frac{-1}{\alpha - a} }$ and   $\displaystyle{{\color{red}v} =  \frac{1}{\alpha - a}}$\\
 
 \n and\\
 
  \n $\displaystyle{\chi(\lambda): = det(\mathbb{M}\mathbb{E}^{-1} \lambda I) = }$ $\left \vert  \begin{array} {ccccc} -\lambda&0&0&\displaystyle{\frac{\beta_{9}}{\gamma}}&1\\
 a{\color{red}u}&a{\color{red}v} -\lambda &0&0&0\\
 \beta_{2}{\color{red}u}&\beta_{2}{\color{red}v}&-\lambda&\displaystyle{\frac{\beta_{4}}{\gamma}}&0\\
 \beta_{8}{\color{red}u}&\beta_{8}{\color{red}v}&\displaystyle{\frac{\beta_{3}}{\beta}}&\displaystyle{\frac{\beta_{8}}{\gamma}} - \lambda&0\\
 0&0&0&0&- \lambda\\
 \end{array}\right\vert$ $= - \lambda$ $\left \vert  \begin{array} {cccc} -\lambda&0&0&\displaystyle{\frac{\beta_{9}}{\gamma}}\\
 a{\color{red}u}&a{\color{red}v} -\lambda &0&0\\
 \beta_{2}{\color{red}u}&\beta_{2}{\color{red}v}&-\lambda&\displaystyle{\frac{\beta_{4}}{\gamma}}\\
 \beta_{8}{\color{red}u}&\beta_{8}{\color{red}v}&\displaystyle{\frac{\beta_{3}}{\beta}}&\displaystyle{\frac{\beta_{8}}{\gamma}} - \lambda\\
 \end{array}\right\vert$
 
 $=  \lambda^{2}$ $\left \vert  \begin{array} {ccc} 
a{\color{red}v} -\lambda &0&0\\
\beta_{2}{\color{red}v}&-\lambda&\displaystyle{\frac{\beta_{4}}{\gamma}}\\
\beta_{8}{\color{red}v}&\displaystyle{\frac{\beta_{3}}{\beta}}&\displaystyle{\frac{\beta_{8}}{\gamma}} - \lambda\\
 \end{array}\right\vert$ $+ \displaystyle{\frac{\beta_{9}}{\gamma}\lambda}$ $\left \vert  \begin{array} {ccc} 
 a{\color{red}u}&a{\color{red}v} -\lambda &0\\
 \beta_{2}{\color{red}u}&\beta_{2}{\color{red}v}&-\lambda\\
 \beta_{8}{\color{red}u}&\beta_{8}{\color{red}v}&\displaystyle{\frac{\beta_{3}}{\beta}}\\
 \end{array}\right\vert$\\
 
 \n $ = \lambda^{2} (a{\color{red}v} -\lambda )$$\left \vert  \begin{array} {cc} 
 - \lambda&\displaystyle{ \frac{\beta_{4}}{\gamma}}\\
\displaystyle{\frac{\beta_{3}}{\beta}}&\displaystyle{\frac{\beta_{8}}{\gamma} - \lambda}\\
\end{array}\right\vert$ $+ \displaystyle{\frac{\beta_{9}}{\gamma}\lambda a}{\color{red}u}$ $\left \vert  \begin{array} {ccc}
\beta_{2}{\color{red}v}&-\lambda\\
\beta_{8}{\color{red}v}& \displaystyle{\frac{\beta_{3}}{\beta}}\\
\end{array}\right \vert $ $- \displaystyle{\frac{\beta_{9}}{\gamma}\lambda (a{\color{red}v} - \lambda)}$ $\left \vert  \begin{array} {ccc}
 \beta_{2}{\color{red}u}&-\lambda\\
\beta_{8} {\color{red}u}& \displaystyle{\frac{\beta_{3}}{\beta}}\\
\end{array}\right\vert$\\ 

\n $\displaystyle{= - \lambda^{2} (a{\color{red}v} -\lambda )[\lambda (\frac{\beta_{8}}{\gamma} -\lambda) + \frac{\beta_{3}\beta_{4}}{\beta \gamma}]}$+ $\displaystyle{\frac{\beta_{9}}{\gamma}\lambda a{\color{red}u}{\color{red}v}(\frac{\beta_{2}\beta_{3}}{\beta} +\beta_{8}\lambda) }$ - $\displaystyle{\frac{\beta_{9}}{\gamma}\lambda (a{\color{red}v} -\lambda){\color{red}u}(\frac{\beta_{2}\beta_{3}}{\beta} + \beta_{8}\lambda)}$\\

\n $\displaystyle{= - \lambda^{2} (a{\color{red}v} -\lambda )(- \lambda^{2}   + \frac{\beta_{8}}{\gamma}\lambda + \frac{\beta_{3}\beta_{4}}{\beta \gamma})}$ + $\displaystyle{\frac{\beta_{9}\beta_{8}}{\gamma} a{\color{red}u}{\color{red}v}\lambda^{2} + \frac{\beta_{9}\beta_{8}}{\gamma} {\color{red}u}\lambda^{2} }$.\\

\n $\displaystyle{= - \lambda^{2} (- a{\color{red}v}\lambda^{2}   + \frac{\beta_{8}}{\gamma}a{\color{red}v}\lambda + \frac{\beta_{3}\beta_{4}}{\beta \gamma}a{\color{red}v}  +  \lambda^{3}   - \frac{\beta_{8}}{\gamma}\lambda^{2} - \frac{\beta_{3}\beta_{4}}{\beta \gamma}\lambda) + \frac{\beta_{9}\beta_{8}}{\gamma} a{\color{red}u}{\color{red}v}\lambda^{2} + \frac{\beta_{9}\beta_{8}}{\gamma} {\color{red}u}\lambda^{2}}$ \\

\n $\displaystyle{= - \lambda^{2}[\lambda^{3} + ( \frac{\beta_{9}\beta_{8}}{\gamma} {\color{red}u}(a{\color{red}v}+1) - a{\color{red}v}- \frac{\beta_{8}}{\gamma}) \lambda^{2} + ( \frac{\beta_{8}}{\gamma}a{\color{red}v}- \frac{\beta_{3}\beta_{4}}{\beta \gamma}) \lambda + \frac{\beta_{3}\beta_{4}}{\beta \gamma}a{\color{red}v}]}$.\\

\n Setting $\displaystyle{a_{1}=  \frac{\beta_{9}\beta_{8}}{\gamma} {\color{red}u}(a{\color{red}v}+1) - a{\color{red}v}- \frac{\beta_{8}}{\gamma}}$, $\displaystyle{a_{2}= \frac{\beta_{8}}{\gamma}a{\color{red}v}- \frac{\beta_{3}\beta_{4}}{\beta \gamma}}$ , $\displaystyle{a_{3}= \frac{\beta_{3}\beta_{4}}{\beta \gamma}a{\color{red}v}}$,\\

\n  and\\

\n  $\displaystyle{\chi(\lambda) = \lambda^{3} + a_{1}\lambda^{2} + a_{2}\lambda + a_{3}}$. If   $\lambda_{1}, \lambda_{1}$ and $\lambda_{3}$ are the zeros of $\displaystyle{\chi(\lambda) = 0}$ \\

\n Then we have\\

\n{\color{red}$\bullet$} $\displaystyle{\lambda_{1} + \lambda_{2} + \lambda_{3} = - a_{1}}$\\

\n{\color{red}$\bullet$} $\displaystyle{\lambda_{1}\lambda_{2} + \lambda_{1}\lambda_{3} +\lambda_{2}\lambda_{3} = a_{2}}$\\

\n{\color{red}$\bullet$} $\displaystyle{\lambda_{1}\lambda_{2}\lambda_{3} = - a_{3}}$.\\

\n and\\

\n {\color{red}{\bf Proposition 4.7}} \\

\n Let $\chi(\lambda) = \lambda^{3 } + a_{1}\lambda^{2 }+ a_{2}\lambda + a_{3}$, so that $\chi$ is uniformly asymptotically stable (uas), it is {\color{red}necessary} that it {\color{red}suffices} that $\Delta_{1} = a_{1} > 0, \Delta_{2 }= a_{1}a_{2 } - a_{3} > 0 $ and $\Delta_{3 }= a_{3}\Delta_{2 } > 0$.\\
\n  A {\color{red}necessary condition } for all the roots of the characteristic polynomial to admit a negative real part, all the coefficients must be positive, that is to say: $a_{1} > 0, a_{2} > 0, ..., a_{3} > 0.$ \\

\n  {\color{red}$\bullet$} As $\displaystyle{\frac{ a_{1} = \beta_{9}\beta_{8}}{\gamma} {\color{red}u}(a{\color{red}v}+1) - a{\color{red}v}- \frac{\beta_{8}}{\gamma} \leq 0}$  then  we can not apply this above proposition  for $\chi(\lambda)$.\\

\n Now, if we consider the discriminant of $\chi$ which is given by :\\

\n $\displaystyle{\Delta_{\chi} = a_{1}^{2}a_{2}^{2} + 18 a_{1}a_{2}a_{3} - 27a_{3}^{2}  - 4a_{2}^{3} - 4a_{1}^{3}a_{3}}$.\\

\n we observe that :\\

{\color{red}$\bullet _{1}$} If $\Delta_{\chi}  > 0$, 3 different  real roots of the equation $\chi(\lambda) = 0$.\\

{\color{red}$\bullet_{2}$} If $\Delta_{\chi}  =  0$, one double or triple root of the equation $\chi(\lambda) = 0$.\\

{\color{red}$\bullet_{3}$} If $\Delta_{\chi}  < 0$, one real root and two complex roots of the equation $\chi(\lambda) = 0$.\\ 

{\color{red}$\bullet_{4}$}  if $\Delta_{\chi}  > 0$, then a necessary and sufficient condition for an equilibrium point  to be \\

 locally asymptotically stable is $a_{1} > 0$, $a_{3} > 0$, $a_{1}a_{2} - a_{3} > 0$.\\

{\color{red}$\bullet_{5}$} if $\Delta_{\chi} < 0$, $a_{1} < 0$, $a_{2} < 0$, then all roots of $\chi(\lambda) = 0$ satisfy the condition $|arg(\lambda)| < \frac{\pi}{2}$.\\

{\color{red}$\bullet_{6}$} if $\Delta_{\chi} > 0$, $a_{1} > 0$, $a_{2} > 0$, $a_{1}a_{2} - a_{3} = 0$, then an equilibrium point is locally asymptotically stable.\\

{\color{red}$\bullet_{7}$} A {\color{red} necessary condition} for an equilibrium point to be locally asymptotically stable is $a_{3} > 0$.\\

{\color{red}$\bullet_{8}$} if the conditions $\Delta_{\chi} < 0$, $a_{1} > 0$, $a_{2} > 0$, $a_{1}a_{2} - a_{3} = 0$ are satisfied, then an equilibrium point is not locally asymptotically stable.\\

\n We remark that {\color{red}$\bullet_{4}$}, {\color{red}$\bullet_{6}$} and {\color{red}$\bullet_{4}$} are not satisfy by the coefficients of $\chi(\lambda)$. So we have to solve the cubic equation  $\chi(\lambda) = 0$ by the Cardan's method  which is ingenious and effective, but quite non-intuitive. \\

\n {\bf {\color{red}Theorem 4.8}} (solutions of cubic equation)\\

\n  Let $P$ the general cubic equation: \\

\n $\displaystyle{ax^{3} + bx^{2} + cx + d = 0  ; a  \neq 0 }$ \hfill { } {\color{red}{\Large $\star$}}\\

\n Then $P$   has solutions: \\

\n $\displaystyle{x_{1} = S + T - \frac{b}{3a}}$ \hfill { } {\color{red}{\Large $\star_{1}$}}\\

\n $\displaystyle{x_{2} = -\frac{S + T}{2} - \frac{b}{3a} + i\sqrt{3}(S - T)}$ \hfill { } {\color{red}{\Large $\star_{2}$}}\\

\n $\displaystyle{ x_{3} =  -\frac{S + T}{2} - \frac{b}{3a}  - i\sqrt{3}(S - T)}$ \hfill { } {\color{red}{\Large $\star_{3}$}}\\
 
\n where 
 
\n $\displaystyle{S = \sqrt[3\,]{R + \sqrt{R^{2} + Q^{3}}}}$\\
 
\n $\displaystyle{T =  \sqrt[3\,]{R - \sqrt{R^{2} + Q^{3}}}}$\\

\n and \\

\n $\displaystyle{Q = \frac{3c - b^{2}}{9a^{2}}}$\\

\n $\displaystyle{R = \frac{9abc - 27a^{2}d -2 b^{3}}{54a^{3}}}$\\

\n The expression $\displaystyle{\Delta = Q^{3} + R^{2}}$ is called the {\color{red} discriminant} of the equation\\

\n See for example Nickalls in  {\color{blue}{\bf [Nickalls]}} for a brief description of Cardan's method.\\

\n {\color{blue} $\bullet$} {\color{red}{\bf Substantial technical difficulties for explicit expression of} $\mathbb{J}_{p^{*}}^{[2]}$}\\

\n In Appendix of {\bf {\color{blue}[Wang-Li]}}, we found that for $n = 2, 3$, and $4$, an explicit expression of second additive compound matrices $A^{[2]}$ of $n\times n$ matrices $A= (a_{i j})_{1 \leq i, j \leq n}$ which are given respectively by:\\

\n {\color{red}$\bullet$} $n = 2$: $ \displaystyle{A^{[2]}  = a_{11} +  a_{22}  = tr(A)}$\\

\n {\color{red}$\bullet$} $n = 3$: $ \displaystyle{A^{[2]}  = }$ $\left ( \begin{array} {ccc} a_{11} +  a_{22}&a_{23}&-a_{13}\\
\quad\\
a_{32}& a_{11} +  a_{33}&a_{12}\\
\quad\\
- a_{31}& a_{21} &  a_{22}+a_{33}\\
\end{array} \right)$\\
\quad\\

\n in the same way that the section 5 of {\bf {\color{blue}[Wang-Li]}} where Li and Wang  studied the stability of  an epidemic model of SEIR type, we apply their criterion to the following epidemic model:\\

\n $\left \{ \begin{array} {c} \frac{dS}{dt} = \Lambda - (\beta_{1}I_{1}  +  \beta_{2}I_{2} )S - \mu S\quad \quad\\\quad\\
\frac{dI_{1}}{dt} = (\beta_{1}I_{1}  +  \beta_{2}I_{2} )S - (\mu  + \gamma)I_{1}\\
\quad\\
\frac{dI_{2}}{dt} = \gamma I_{1}  -(\mu + d)I_{2} \quad \quad \quad \quad\quad\quad\\
\end{array} \right.$ \hfill { } {\bf{\color{blue}(*$_{1}$) }}\\

\n where $(\Lambda, \beta_{1}, (\beta_{2}, \mu, \gamma, d)$ are given parameters.\\

\n {\color{red}$\bullet$} Determination of equilibrium points of the system {\bf{\color{blue}(*$_{1}$) }} and calculation of basic reproduction number $\mathcal{R}_{0}$ \\

\n Let $E = (S, I_{1}, I_{2})$ then the Jacobian matrix of above system is :\\

\n $\mathbb{J}_{E} = $ $\left ( \begin{array} {ccc} - (\beta_{1}I_{1}  +  \beta_{2}I_{2}  + \mu)& - \beta_{1}S&-\beta_{2}S\\
\quad\\
\beta_{1}I_{1}  +  \beta_{2}I_{2}&\beta_{1}S - \mu - \gamma&\beta_{2}S\\
\quad\\
0& \gamma& - \mu - d\\
\end{array} \right )$ \hfill { } {\bf{\color{blue}(*$_{2}$)}}\\
\quad\\

\n Now, we consider the following equations:\\

\n  $\left \{ \begin{array} {c} \frac{dS}{dt} = \Lambda - (\beta_{1}I_{1}  +  \beta_{2}I_{2} )S - \mu S = 0  \quad \quad {\color{red} (1)}\\\quad\\
\frac{dI_{1}}{dt} = (\beta_{1}I_{1}  +  \beta_{2}I_{2} )S - (\mu  + \gamma)I_{1} = 0 \quad {\color{red} (2)}\\
\quad\\
\frac{dI_{2}}{dt} = \gamma I_{1}  -(\mu + d)I_{2} = 0 \quad \quad \quad \quad\quad \quad {\color{red} (3)}\\
\end{array} \right.$ \hfill { } {\bf{\color{blue}(*$_{3}$) }}\\

\n then we observe that $E^{0} = (\frac{\Lambda}{\mu}, 0, 0)$ is a trivial equilibrium point  of {\bf{\color{blue}(*$_{3}$)}} (Disease free equilibrium point)\\

\n and so \\

\n $\mathbb{J}_{E^{0}} = $ $\left ( \begin{array} {ccc} -  \mu& - \frac{\beta_{1}\Lambda}{\mu}&- \frac{\beta_{2}\Lambda}{\mu}\\
\quad\\
0& \frac{\beta_{1}\Lambda}{\mu} - \mu - \gamma& \frac{\beta_{2}\Lambda}{\mu}\\
\quad\\
0& \gamma& - \mu - d\\
\end{array} \right )$ \hfill { }{\bf{\color{blue}(*$_{4}$)}}\\
\quad\\

\n By using the next generation matrix method,  the basic reproduction number $\mathcal{R}_{0}$ is obtained as the spectral radius of
matrix $(-\mathbb{F}\mathbb{V}^{-1})$ at disease free equilibrium point where  $\mathbb{F}$ and $\mathbb{V}$ are as below :\\

\n $\mathbb{F} =$  $\left ( \begin{array} {cc} \displaystyle{\frac{\beta_{1}\Lambda}{\mu}}& \displaystyle{\frac{\beta_{2}\Lambda}{\mu}}\\
\quad\\
0&0\\
\end{array} \right ) $ , $\mathbb{V} = $ $\left ( \begin{array} {cc} \displaystyle{-\mu - \gamma}& 0\\
\quad\\
\gamma&- \mu - d\\
\end{array} \right ) $ , $\mathbb{V}^{-1} = $ $\left ( \begin{array} {cc} \displaystyle{-\frac{1}{\mu + \gamma}}& 0\\
\quad\\
\displaystyle{-\frac{\gamma}{(\mu + d)(\mu + \gamma)}}&\displaystyle{- \frac{1}{\mu + d}}\\
\end{array} \right ) $\\

\n and $- \mathbb{F}\mathbb{V}^{-1} = \Lambda$ $\left ( \begin{array} {cc} \displaystyle{\frac{\beta_{1}(\mu + d) + \beta_{2}\gamma}{\mu(\mu + d)(\mu + \gamma)}}& \displaystyle{\frac{\beta_{2}}{\mu(\mu + d)}}\\
\quad\\
0&0\\
\end{array} \right ) $\\

\n It follows that :\\

\n {\color{red}$\mathcal{R}_{0} = \Lambda$ $\displaystyle{\frac{\beta_{1}(\mu + d) + \beta_{2}\gamma}{\mu(\mu + d)(\mu + \gamma)}}$} \hfill { } {\bf{\color{blue}(*$_{5}$)}}\\

\begin{center}
 Evolution of $\mathcal{R}_{0} $ with respect $\mu$ \\
 \includegraphics[scale=0.26]{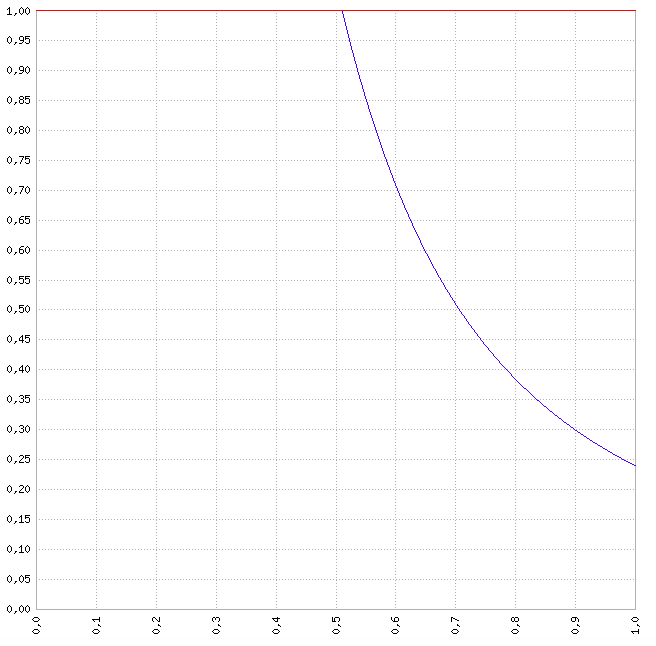} \\
 where $\beta_{1} = 0,3$, $\beta_{2} = 0,8$, $\gamma = 0,1$ , $\Lambda = 0,7$ and $d = 0,04$\\
 \end{center}

\n $\mathbb{J}_{E^{0}}^{[2]} = $ $\left ( \begin{array} {ccc} \frac{\beta_{1}\Lambda}{\mu} -  2\mu - \gamma&  \frac{\beta_{2}\Lambda}{\mu}& \frac{\beta_{2}\Lambda}{\mu}\\
\quad\\
\gamma&- 2\mu - d& -\frac{\beta_{1}\Lambda}{\mu}\\
\quad\\
0& 0& - 2\mu - \gamma - d\\
\end{array} \right )$ \hfill { } {\bf{\color{blue}(*$_{6}$)}}\\
\quad\\

\n Now let  $I_{1} \neq 0$ then from {\color{red}(3)} we deduce that :\\

\n $\displaystyle{ I_{2} = \delta I_{1} ; \delta = \frac{\gamma}{\mu + d}}$  \hfill { } {\bf{\color{blue}(*$_{7}$)}}\\

\n From  {\bf{\color{blue}(*$_{7}$)}} and {\color{red}(2)} we deduce that :\\

\n  $\displaystyle{ S = \frac{\mu + \gamma}{\beta_{1} + \beta_{2}\delta}}$  \hfill { } {\bf{\color{blue}(*$_{8}$)}}\\

\n Now from {\color{red}(1)} + {\color{red}(2)}, we deduce that :\\

\n  $\displaystyle{ I_{1} = \frac{\Lambda - \mu S}{\mu + d} = \frac{\Lambda(\beta_{1} + \beta_{2}\delta) - \mu(\mu + \gamma)}{\beta_{1} + \beta_{2}\delta)(\mu + d)}}$ \hfill { } {\bf{\color{blue}(*$_{8}$)}}\\

\n Let $\displaystyle{E^{*} = (S^{*}, I_{1}^{*}, I_{2}^{*})}$ where $\displaystyle{S^{*} = \frac{\mu + \gamma}{\beta_{1} + \beta_{2}\delta}, I_{1}^{*} = \frac{\Lambda - \mu S^{*}}{\mu + d}}$ and  $\displaystyle{I_{2}^{*} = \delta I_{1}^{*}}$ where $\displaystyle{\delta = \frac{\gamma}{\mu + d}}$.\\

\n then

\n {\color{blue}$\bullet$} The Jacobian matrix at the  endemic equilibrium point $E^{*} = (S^{*}, I_{1}^{*}, I_{2}^{*})$ of the system  {\bf{\color{blue}(*) }} is:\\

\n $J_{E^{*}} = $ $\left ( \begin{array} {ccc} -\beta_{1}I_{1}^{*}  -  \beta_{2}I_{2}^{*} - \mu& -\beta_{1}S^{*}& -\beta_{2}S^{*}\\
\quad\\
\beta_{1}I_{1}^{*}  +  \beta_{2}I_{2}^{*}&\beta_{1}S^{*} - \gamma - \mu & \beta_{2}S^{*}\\
\quad\\
0& \gamma& - d -\mu\\
\end{array} \right )$\\
\quad\\
 
\n  and\\

\n {\color{blue}$\bullet$} the second additive compound matrix associated to $J_{E^{*}}$  is:\\

\n $J_{E^{*}}^{[2]} = $ $\left ( \begin{array} {ccc} -\beta_{1}I_{1}^{*}  -  \beta_{2}I_{2}^{*} + \beta_{1}S^{*} - \gamma - 2\mu&\beta_{2}S^{*}&\beta_{2}S^{*}\\
\quad\\
\gamma& -\beta_{1}I_{1}^{*}  -  \beta_{2}I_{2}^{*} - 2\mu -d & - \beta_{1}S^{*}\\
\quad\\
0&\beta_{1}I_{1}^{*}  +  \beta_{2}I_{2}^{*} &\beta_{1}S^{*} - d - \gamma - 2\mu\\
\end{array} \right )$\\
\quad\\

\n {\bf {\color{red} Proposition 4.9}}\\

\n Let 

\n  {\color{red} $\bullet_{1}$} $\displaystyle{\beta_{2} < \frac{\gamma}{\delta^{2}}}$\\

\n  {\color{red} $\bullet_{2}$} $\displaystyle{\frac{ (\mu + \gamma)(\mu +d)(\beta_{1} + \beta_{2}\delta)}{\Lambda(\beta_{1} + \beta_{2}\delta) - \mu(\mu + \gamma)} + \frac{\beta_{1}(\mu + \gamma)}{\beta_{1} + \beta_{2}\delta} < d + \gamma + 2\mu}$\\

\n  {\color{red} $\bullet_{3}$} $\displaystyle{ \beta_{2}\delta \Lambda + \mu(\gamma +\mu) < \beta_{1}\Lambda }$.\\

\n  then the endemic equilibrium point of {\bf{\color{blue}(*) }} is asymptotically stable.\\

\n {\bf {\color{red} Proof}}\\

\n Let $\mathbb{P} = $ $\left (\begin{array}{ccc}I_{2}^{*}&0&0\\
\quad\\
0&I_{1}^{*}&0\\
\quad\\
0&0&S^{*}\\
\end{array} \right ) $\\

\n  then the matrix $J_{E^{*}}^{[2]}$  is similar to matrix $\displaystyle{\mathbb{A} = \mathbb{P} J_{E^{*}}^{[2]} \mathbb{P}^{-1} = (a_{ij})_{1\leq i,j \leq 3}}$ which is given by :\\

\n $\displaystyle{\mathbb{A}  = }$ $\left ( \begin{array} {ccc} -\beta_{1}I_{1}^{*}  -  \beta_{2}I_{2}^{*} + \beta_{1}S^{*} - \gamma - 2\mu&\beta_{2}S^{*}\frac{I_{2}^{*}}{I_{ 1}^{*}}&\beta_{2}S^{*}\frac{I_{2}^{*}}{S^{*}}\\
\quad\\
\gamma\frac{I_{1}^{*}}{I_{ 2}^{*}}& -\beta_{1}I_{1}^{*}  -  \beta_{2}I_{2}^{*} - 2\mu -d & - \beta_{1}S^{*}\frac{I_{1}^{*}}{S^{*}}\\
\quad\\
0&(\beta_{1}I_{1}^{*}  +  \beta_{2}I_{2}^{*} )\frac{S^{*}}{I_{ 1}^{*}}&\beta_{1}S^{*} - d - \gamma - 2\mu\\
\end{array} \right )$\\
\quad\\

\n Under the conditions {\color{red} $\bullet_{1}$} and {\color{red} $\bullet_{2}$} , we observe that the diagonal elements of $\mathbb{A}$ are negative and \\

\n (1) $\displaystyle{a_{11} + \mid a_{12} \mid + \mid a_{13} \mid < 0}$\\

\n (2) $\displaystyle{a_{22} + \mid a_{21} \mid + \mid a_{23} \mid < 0}$\\

\n (1) $\displaystyle{a_{33} + \mid a_{32} \mid < 0}$\\

\n i.e  $\mathbb{A}$ is {\color{red}diagonally dominant } in rows.\\

\n In order to apply the corollary of the Li-Wang criterion, it remains to calculate the determinant of $J_{E^{*}}$  \\

$detJ_{E^{*}} = $ $\left \vert \begin{array} {ccc} -\beta_{1}I_{1}^{*}  -  \beta_{2}I_{2}^{*} - \mu& -\beta_{1}S^{*}& -\beta_{2}S^{*}\\
\quad\\
\beta_{1}I_{1}^{*}  +  \beta_{2}I_{2}^{*}&\beta_{1}S^{*} - \gamma - \mu & \beta_{2}S^{*}\\
\quad\\
0& \gamma& - d -\mu\\
\end{array} \right \vert$\\
\quad\\

\n Under condition {\color{red} $\bullet_{3}$} we deduce that $detJ_{E^{*}} < 0 $\\

\n {\color{red}$\bullet$} $n = 4$: $ \displaystyle{A^{[2]}  = }$ $\left ( \begin{array} {cccccc} a_{11} +  a_{22}&a_{23}&a_{24}&-a_{13}&-a_{14}&0\\
\quad\\
a_{32}& a_{11} +  a_{33}&a_{34}&a_{12}&0&-a_{14}
\quad\\
a_{42}& a_{43} & a_{11}+a_{44}&0&a_{12}&a_{13}\\
\quad\\
- a_{13}& a_{21} &0&  a_{22}+a_{33}&a_{34}&- a_{24}\\
\quad\\
- a_{41}& 0 &a_{21}& a_{43}& a_{22}+a_{44}&a_{23}\\
\quad
0& -a_{41} &a_{31}& - a_{42}&a_{32}& a_{33}+a_{44}\\
\end{array} \right)$\\

\n In next lemma, we give the explicit  entries of second additive compound matrix of $n\times n$ matrix $A= (a_{i j})$  where {\color{red}$n = 5$}\\
\quad\\
\n {\color{red}{\bf Lemma 4.10}}\\

\n For $n = 5$, an explicit  expression of second additive compound matrix $A^{[2]}$ is given by: \\

\n \includegraphics[scale=0.35]{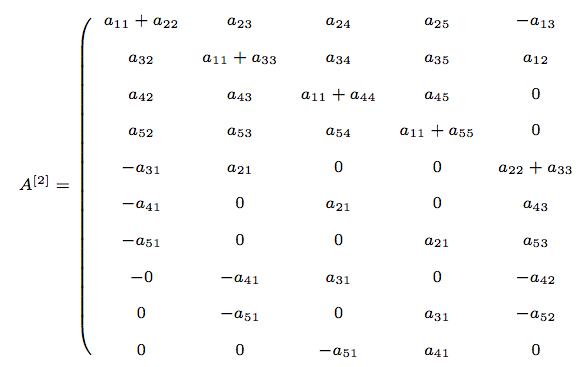}\includegraphics[scale=0.35]{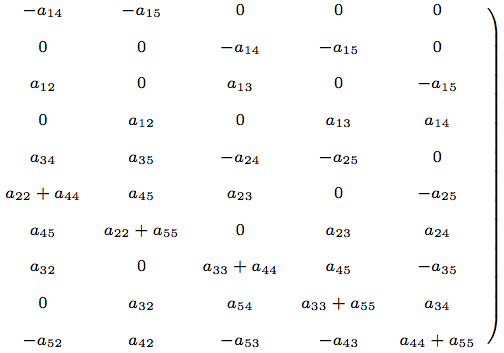}  \\

\n {\color{green}$\bullet$} As  {\color{red}$a_{13} = 0$, $a_{23} = 0$, $a_{24} = 0$, $a_{25} = 0$, $a_{31} = 0$, $a_{35} = 0$, $a_{41} = 0$, $a_{45} = 0$} and {\color{red}$a_{54} = 0$} in $\mathbb{J}_{p^{*}}^{[2]}$ , we deduce that the explicit expression of second additive compound matrix $\mathbb{J}_{p^{*}}^{[2]}$ where $\displaystyle{P^{*} = (\frac{B}{\mu}, 0, 0, 0, 0)}$ is :  \\
\quad\\

\n{\color{red}$\mathbb{J}_{p^{*}}^{[2]} = $}{\scriptsize $\left ( \begin{array} {cccccccccc}a_{11}+a_{22}&{\color{red}0}&{\color{red}0}&{\color{red}0}&{\color{red}0}&-a_{14}&-a_{15}&0&0&0\\
\quad\\
a_{32}&a_{11} + a_{33}&a_{34}&a_{35}&a_{12}&0&0&-a_{14}&-a_{15}&0\\
\quad\\
a_{42}&a_{43} &a_{11} + a_{44}&a_{45}&0&a_{12}&0 &{\color{red}0}&0&-a_{15} \\
\quad\\
a_{52}&a_{53} &{\color{red}0}&a_{11} + a_{55}&0&0&a_{12} &0&{\color{red}0}&a_{14} \\
\quad\\
{\color{red}0}&a_{21} &0&0&a_{22} + a_{33}&a_{34}&a_{35} &{\color{red}0}&-a_{25}&0 \\
\quad\\
{\color{red}0}&0 & a_{21}&0&a_{43}&a_{22} + a_{44}&{\color{red}0} &{\color{red}0}&0&{\color{red}0}  \\
\quad\\
-a_{51}&0 &0&a_{21}&a_{53}&{\color{red}0} &a_{22} + a_{55}&0&{\color{red}0}&{\color{red}0} \\
\quad\\
0&-a_{41} &{\color{red}0}&0&-a_{42}&a_{32}&0 &a_{33} + a_{44}&{\color{red}0}&{\color{red}0} \\
\quad\\
0&-a_{51} &0&{\color{red}0}&-a_{52}&0&a_{32}&a_{54}&a_{33} + a_{55}&a_{34}\\
\quad\\
0&0 &-a_{51}&{\color{red}0}&0&-a_{52}& a_{42}&-a_{53}&- a_{43}&a_{44} + a_{55} \\
\end{array} \right )$} $=$\\

\n \includegraphics[scale=0.33]{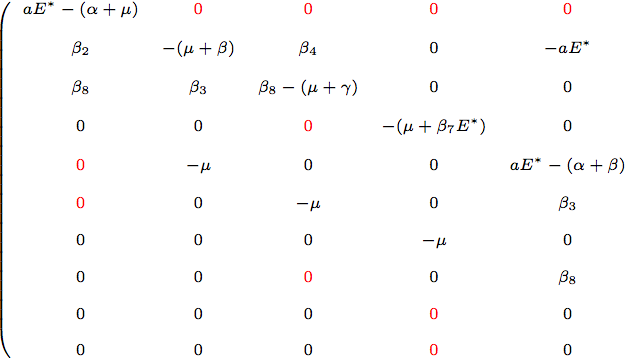}\includegraphics[scale=0.33]{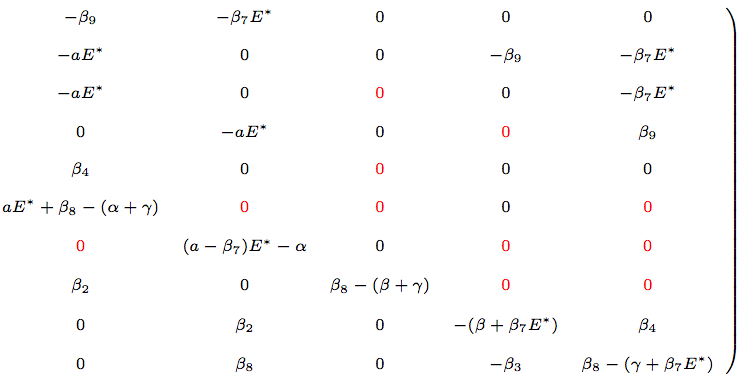}  \\

\n where\\

\n ${\color{red}a} = \beta_{1} - \beta_{10}$, ${\color{red}\alpha} = \beta_{2} + \beta_{6} + \beta_{8} + \mu$, ${\color{red}\beta} = \beta_{2} + \beta_{5} + \mu$ and ${\color{red}\gamma} = \beta_{4} + \beta_{9} + \mu$.\\

\n {\color{red}{\bf Theorem 4.11}} \\ 

\n  (i) If {\color{red}$\beta_{10} < \beta_{1}$} and {\color{blue}$\displaystyle{ 2\beta_{8}{\color{red} a}E^{*} + {\color{red} \alpha}{\color{red} \gamma} + \beta_{8}\beta_{9} > \beta_{8}{\color{red} \alpha}}$} then {\color{red}$\displaystyle{det \mathbb{J}_{p^{0}} < 0}$}.\\

\n (ii) if we have :\\

\n {\color{red}(a)} $\beta_{3} < \beta$\\

\n {\color{red}(b)} $\displaystyle{\frac{2(\beta_{1} - \beta_{10})B}{\mu}  <  \beta_{6} + \mu }$\\

\n {\color{red}(c)} $\displaystyle{\beta_{9} < \frac{\beta_{7}B}{\mu}}$\\

\n {\color{red}(d)} $\displaystyle{ \beta_{8} +  \frac{\beta_{7}B}{\mu} < \mu}$\\

\n then   $\mathbb{J}_{p^{0}}^{[2]}$ is {\color{red}diagonally dominant } in columns.\\

\n (iii) the equilibrium point of {\color{blue}(4.1)} is asymptotically stable. \hfill { } $\blacklozenge$\\

\n {\color{red}{\bf Conclusion}}\\

\n In second paper (Part II), In order to control the Covid-19 system, i.e., force the trajectories to go to the equilibria we will add some control parameters with uncertain parameters to stabilize  the five-dimensional  Covid-19 system studied in this paper. \\
\n Based on compound matrices theory, we have constructed in {\color{blue}[Intissar]} the controllers :\\

\n $\mathbb{U} = $ $\left (\begin{array} {ccccc} 0&{\color{red}u_{1}}&0&0&0\\
\quad\\
{\color{red}u_{2}}&0&0&0&0\\
\quad\\
0&0&0&0&0\\
\quad\\
0&0&0&0&0\\
\quad\\
{\color{red}u_{3}}&0&0&0&0\\
\end{array}\right ) $i.e. $\mathbb{U} = (u_{ij}) $ where  $u_{ij} = 0$ except $(u_{12}, u_{21}, u_{51}) \in \mathbb{R}^{3}$; $1 \leq i, j \leq 5$\\

\n to stabilize the system {\color{blue}(4.1)}, in particular to study the stability of following matrix :\\
\n $\displaystyle{\mathbb{J}_{p^{*}, u_{1},u_{2},u_{3}} =  \mathbb{J}_{p^{*}} + \mathbb{U}}$ and its second additive compound matrix $\displaystyle{ (\mathbb{J}_{p^{*}} + \mathbb{U})^{[2]}}$,  by applying  again the criterion of Li-Wang  on second compound matrix associated to the system {\color{blue}(4.1)} with these controllers. \\
\n We have constructed a Lyapunov function $\mathbb{L}$ of the system {\color{blue}(4.1)} for apply the classical Lyapunov theorem and to get :\\

\n {\bf{\color{red}Theorem 4.12}}\\

\n (i)  $(0, 0, 0, 0, 0)$ is a stable equilibrium point in the sense of Lyapunov.\\

\n  (ii) $\mathbb{V}(e, i, c, h, d) < 0, 0 < \mid\mid  (e, i, c, h, d ) \mid\mid < r_{1}$ for some $r_{1}$, i.e. if $\mathbb{L}$ is lnd.\\

 \n (iii)  $(0, 0, 0, 0, 0)$ is an asymptotically stable equilibrium point. \hfill { } $\blacklozenge$\\
 
 \begin{center}
 {\bf { \Large {\color{red} References}}}
 \end{center}
  
\n  {\color{blue}[Aimar et al] } M-T Aimar, A. Intissar and J-K Intissar :. Elementary Mathematical Analysis of Chaos: in Classical and a few New Dynamical Systems From Lorenz to Gribov-Intissar systems, Kindle Edition, ASIN: B07QV7JPLW.  (2019) \\

\n {\bf{\color{blue}[Aitken]}} Aitken, A.C. :. Determinants and matrices, Oliver $\&$ Boyd, Edinburgh (1967)\\

\n {\bf{\color{blue}[Arnold]}} V.I. Arnold, Geometric Methods in the Theorey of Ordinary Differential Equations, 2nd ed., Springer-Verlag, New York, 1998.\\

\n {\color{blue} [Banks et al]} Banks, J., Brooks, J., Cairns, G., Davis, G. and Stacey, P. , On Devaneys Defi- nition of Chaos. American Mathematical Monthly, 1992 ; 99, 332-334.\\

\n {\color{blue}[Beretta et al]} E. Beretta and V. Capasso :. On the general structure of epidemic systems . Global asymp- totic stability, Corn. $\&$ Maths. with Appls. Vol, 12A, No. 6. (1986), pp, 677-694\\

\n {\color{blue}[Berman et al]}  Berman, A. and R.J. Plemmons, Nonnegative Matrices in the Mathematical Sciences, Academic Press, New York, (1979), chapter 6.\\

\n {\bf{ \color{blue}[Brenner1]}} J.L. Brenner :. A bound for a determinant with dominant main diagonal, Proc. Amer. Math. Soc. vol. 5 (1954) pp. 631-634.\\

\n {\bf{ \color{blue}[Brenner2]}} J.L. Brenner :.  Bounds for determinants II , Proc. Amer. Math. Soc. vol. 8 (1957) pp. 532-534.\\

\n {\color{blue}{\bf [Bylov et al ]}} B. F. Bylov, R. E. Vinograd, D. M. Grobman, and V. V. Nemytskii, :.Theory of Lyapunov Exponents and Its Application to Problems of Stability (Moscow, Nauka, 1966)\\

\n {\color{blue}[Cheng-Shan]} Cheng, Z. J. and Shan, J. 2019 Novel coronavirus : where we are and what we know. Infection, (2020)1-9.\\

\n {\color{blue}[Coppel]} Coppel, W. A., :. Stability and Asymptotic Behavior of Differential Equations, Heath, Boston, (1965)\\

\n {\color{blue}{\bf [Dahlquist]}}  G. Dahlquist, :. Stability and Error Bounds in the Numerical Integration of Ordinary Differential Equations, Kungl. Tekn. Hogsk. Handl. Stockholm. (Stockholm, 1959), Vol. 130.\\

\n {\color{blue}{\bf [Daletskii-Krein]}} Yu. L. Daletskii and M. G. Krein, Stability of Solutions of Differential Equations in Banach Space, in Nonlinear Analysis and Its Applications (Nauka, Moscow, 1970).\\.

\n {\color{blue}[Diekmann et al]} Diekmann, O., Heesterbeek, J. A. P., and Metz, J. A. . : On the definition and the computation of the basic reproduction ratio $R_{0}$ in models for infectious diseases in hetero- geneous populations. Journal of Mathematical Biology, 28 (4), (1990) 365-382.\\

\n {\color{blue}[Driessche et al 1] } P. van den Driessche and James Watmough  Reproduction numbers and sub-threshold endemic equilibria for compartmental models of disease transmission, Mathematical Biosciences 180 (2002) 29-48\\

\n {\color{blue}[Driessche et al 2] } P. van den Driessche and James Watmough :. Further Notes on the Basic Reproduction Number, April (2008)
DOI: $10.1007/978-3-540-78911-6_6$ (In book: Mathematical Epidemiology, Springer Lecture Notes in Mathematics Vol. 1945).\\

\n {\bf {\color{blue}[Elsner et al]}} L. Elsner, D. Hershkowitz and H. Schneider :. Bounds  on norms of compound matrices and on products of eigenvalues, arXiv: math/9802106v1 [math. RA] 22 Fb 1998\\

\n  {\color{blue}[Fan]}, Fan, K., :.Topological Proofs for Certain Theorems on Matrices with Nonnegative Elements, Montach. Math., 62, 219-237, (1958). \\

\n {\color{blue}[Fiedle]} Fiedler, M. :. Additive compound matrices and inequality for eigenvalues of stochastic matrices, Czech. Math. J. 99, (1974), 392-402.\\

\n  {\color{blue}{\bf [Fiedler-Ptàk]}} M. FIEDLER et V. PTAK, On matrices with positive off diagonal éléments and positive principal minors. Czec. Math. J., 12 (87), 1962, p. 382-400.\\

\n {\bf{\color{blue}[ Golub et al]}} Golub, G. H. and Van Loan, C. F. "Positive Definite Systems." $\S$ 4.2 in Matrix Computations, 3rd ed. Baltimore, MD: Johns Hopkins University Press, pp. 140-141, 1996.\\

\n  {\bf{\color{blue}[Grobman]}} D. Grobman, Homeomorphisms of Systems of Differential Equations (Russian), Dokl. Akad., Nauk., Vol. 128, 1959, pp 880-881.\\
 
\n {\bf{\color{blue}[Hartman1]}}  P. Hartman, A Lemma in the Theory of Structural Stability of Differential Equations, Proc. Amer. Math. Soc., Vol. 11, 1960, pp 610-620.\\
 
\n {\bf{\color{blue}[Hartman2]}}  P. Hartman, Ordinary Differential Equations, Wiley, New York, 1964.\\

\n {\color{blue}[Hawkins-Simon]} Hawkins, D. and H.A. Simon, :. Note: Some Conditions of Macroeconomic Stability, Econometrica 17 (1949), 245-248\\

\n {\color{blue}[Heffernan et al]} J. M. Heffernan, R. J. Smith, and L. M. Wahl, Perspectives on the basic reproductive ratio, J. R. Soc. Interface, 2 (2005), pp. 281-293.\\

\n {\bf{\color{blue}[Horn-Johson]}} R.A. Horn, C.R. Johnson, Topics in Matrix Analysis, Cambridge University, Cambridge, (1991).\\

\n {\color{blue}[Intissar]} Construction of controllers to stabilize five-dimensional Covid-19 mathematical system with uncertain parameters (preprint 2020).\\

\n {\color{blue}[Intissar et al]} A. Intissar and J.-K. Intissar :. Differential Calculus, Principles and Applications, CEPADUES Editions (2017) I.S.B.N. : 9782364935884.\\

\n \n {\color{blue}[Johnson]} Johnson, C. R. "Positive Definite Matrices." Amer. Math. Monthly 77, 259-264 1970.\\

\n {\color{blue}[Khanh}  Nguyen Huu Khanh : Stability analysis of an influenza virus model with disease resistance, Journal of the Egyptian Mathematical Society, (2015) .\\

\n {\bf{\color{blue}[Kravvaritis et al]}} Christos Kravvaritis and  Marilena Mitrouli :.Compound matrices: properties, numerical issues and analytical computations: (2008) \\ https://www.researchgate.net/publication/220393932.\\

\n {\bf {\color{blue}[Liao]}} X. Liao, :. On asymptotic behavior of solutions to several classes of discrete dynamical systems, Science in China A, vol. 45, no. 4, pp. 432-442, 2002.\\

\n {\color{blue}[Li et al]}, Michael Y. Li, John R. Graef, Liancheng Wang, János Karsai :. Global dynamics of a SEIR model with varying total population size, Mathematical Biosciences, Volume 160, Issue 2, Au- gust 1999, Pages 191-213 https ://doi.org/10.1016/S0025-5564(99)00030-9.\\

\n {\color{blue}[Li-Wang]} Michael Y. Li and Liancheng Wang :. A Criterion for Stability of Matrices Journal of Ma- thematical Analysis and Applications, 225, (1998), 249- 264.\\

\n {\color{blue}{\bf [Lozinski]}}  S. M. Lozinski :. Error estimate for numerical integration of ordinary differential equations. I, Izv. Vyssh. Uchebn. Zaved. Mat., No. 5, 52?90 (1958).\\

\n {\bf{\color{blue}[Marcus]}}  Marcus, M.: Finite Dimensional Multilinear Algebra, Two Volumes. Marcel Dekker, New York (1973-1975)\\

\n \n {\color{blue}[Marcus et al 1]} Marcus, M. and Minc, H. Introduction to Linear Algebra. New York: Dover, p. 182, 1988.\\

\n {\color{blue}[Marcus et al 2]} Marcus, M. and Minc, H. "Positive Definite Matrices." $\S$ 4.12 in A Survey of Matrix Theory and Matrix Inequalities. New York: Dover, p. 69, 1992.\\

\n {\bf{\color{blue}[Marshall et al]}}, A.W. Marshall, I. Olkin, and B.C. Arnold. Inequalities: Theory of Majorization and Its Applications, second edition. Springer, New York, 2011.\\

\n  {\bf{\color{blue}[Mitrouli et al]}} Mitrouli, M., Koukouvinos, C.: On the computation of the Smith normal form of compound matrices. Numer. Algorithms 16, 95-105 (1997)\\

\n {\color{blue}[Muldowney]} J. S. Muldowney :.Compound matrices and ordinary differential equations, Rocky Mountain Journal of Mathematics, Volume 20, Number 4, 1990.\\

\n {\color{blue}{\bf Nickalls}} R W D Nickalls, A new approachto solving the cubic: Cardan?s solution revealed, Mathematical Gazette, Vol.77, 1993.\\

\n {\color{blue}[Ostrowski]} A.M. Ostrowski, Über die Determinaten mit uberweigender Hauptdiagonale, Comment. Math. Heiv., 10, 69- 96,1937.\\

\n {\bf{ \color{blue}[Ostrowski1]}} A.  Ostrowski :. Note  on bounds for determinants with dominant principal diagonal, Proc. Amer. Math. Soc. vol. 3 (1952) pp. 26-30.\\

\n {\bf {\color{blue}[Ostrowski2]}}:. Note  on bounds for determinants, Duke Math. J, vol. 22 (1955) pp. 95-102.\\

\n \n {\color{blue}[Pease]} Pease, M. C. Methods of Matrix Algebra. New York: Academic Press, 1965.\\

\n {\bf{\color{blue}[Perko,]}} L. Perko, Differential Equations and Dynamical Systems, 3rd ed., Springer-Verlag, New York, 2001.\\
 
\n {\color{blue}{\bf [Perov]}} A. I. Perov, :. New stability criteria for constant-coefficient linear systems, Avtomat. Telemekh., No. 2, 22-33 (2002) [Automat. Remote Control 63 (2),189-199 (2002)].\\

\n {\bf{\color{blue}[Poincaré]}} H. Poincaré, Sur le problème des trois corps et les équations, Dynamique Acta Math., Vol.13, 1890, pp 1-270.\\
 
\n {\bf {\color{blue}[Price] }} G. B. Price:. Determinants with dominant principal diagonal , proc, Amer, Math. Soc. vol. 2 (1951) pp. 497-502.\\

\n  {\color{blue}[Schneider]} H. Schneider, :. An Inequality for Latent Roots Applied to Determinants with Dominant Principal Diagonal, J. London Math. Society, 28, 8-20, (1953),\\ 

\n {\color{blue} [Shah et al]} Nita H. Shah, Ankush H. Suthar, Ekta N. Jayswal :. Control Strategies to Curtail Transmission of COVID-19, https ://doi.org 10.1101/2020.04.04.20053173.\\

\n  {\color{blue}[Sun et al]}. Chengjun Sun, Ying-Hen Hsieh :. Global analysis of an SEIR model with varying population size and vaccination, Applied Mathematical Modelling 34 (2010) 2685-2697.\\

\n {\bf{\color{blue}[X-Wang]}} A  global Hartman-Grobman theorem , arXiv:2002.06094v1 [math.DS] 14 Feb 2020.\\

\end{document}